\documentclass{gtart}

\def\ifplaintex{\expandafter\ifx\csname documentclass\endcsname\relax}

\ifplaintex 
\else
\fi

\expandafter\ifx\csname beginpicture\endcsname\relax
\expandafter\ifx\csname documentclass\endcsname\relax
\input pictex \else
\input prepictex \input pictex \input postpictex \fi\fi

\def\gt{{\mathsurround=0pt\it $\cal G\mskip-2mu$eometry \&\ 
$\cal T\!\!$opology}}        

\def\gtp{{\mathsurround=0pt\it $\cal G\mskip-2mu$eometry \&\ 
$\cal T\!\!$opology $\cal P\!$ublications}}  

\def\volumenumber#1{\def\thevolumenumber{#1}}
\def\papernumber#1{\def\thepapernumber{#1}}
\def\volumeyear#1{\def\thevolumeyear{#1}}

\def\pagenumbers#1#2{\def\startpage{#1}\def\finishpage{#2}}
\def\published#1{\def\publishdate{#1}}
\def\proposed#1{\def\theproposer{#1}}
\def\seconded#1{\def\theseconders{#1}}
\def\received#1{\def\receiveddate{#1}}

\def\accepted#1{\def\accepteddate{#1}}

\long\def\asciiabstract#1{\long\def\theasciiabstract{#1}}

\let\\\par\let\thevolumenumber\relax
\let\thepapernumber\relax
\let\thevolumeyear\relax\let\thesamplenumber\relax\let\startpage\relax
\let\finishpage\relax\let\publishdate\relax\let\receiveddate\relax
\let\reviseddate\relax\let\accepteddate\relax\let\theasciititle\relax
\let\theasciiauthors\relax
\let\theasciiabstract\relax
\let\theasciiemail\relax\let\theshortauthors\relax\let\theshorttitle\relax

\long\def\maketitlep{   

\count0=\startpage

\gt\hfill      
\beginpicture
\setcoordinatesystem units <0.33truein, 0.33truein> point at 2.2 0.9
\setplotsymbol ({$\cal G$})
\plotsymbolspacing=9truept
\circulararc 315 degrees from 0 1 center at 0 0
\setplotsymbol ({$\cal T$})
\circulararc 315 degrees from 1 -1 center at 1 0
\endpicture
\break
{\small\ifx\thesamplenumber\relax 
Volume \else Sample
\fi\thevolumenumber\ (\thevolumeyear)
\startpage--\finishpage\nl
Published: \publishdate}
\vglue 0.5truein plus 0.4fil minus 0.1truein

{\parskip=0pt\leftskip 0pt plus 1fil\def\\{\par\smallskip}{\ifplaintex\large
\else\Large\fi\bf\thetitle}\par\medskip}   

\vglue 0pt plus 0.1fil 

{\parskip=0pt\leftskip 0pt plus 1fil\def\\{\par}{\sc\theauthors}
\par\medskip}

\vglue 0pt plus 0.1fil 

{\small\parskip=0pt\let\newline\\
{\leftskip 0pt plus 1fil\def\\{\par}{\sl\theaddress}\par}
\expandafter\ifx\theemail\relax    
\relax\else\vglue 5pt plus 0.02fil minus 2pt\def\\{\stdspace{\rm 
and}\stdspace} 
\cl{Email:\stdspace\tt\theemail}\fi
\ifx\theurl\relax                  
\relax\else\vglue 5pt plus 0.02fil minus 2pt\def\\{\stdspace{\rm 
and}\stdspace}
\cl{URL:\stdspace\tt\theurl}\fi\par}

\vglue 7pt plus 0.3fil minus 3pt

{\bf Abstract}
\vglue 5pt plus 0.1fil minus 2pt

\theabstract

\vglue 7pt plus 0.3fil minus 3pt

{\bf AMS Classification numbers}\quad Primary:\quad \theprimaryclass

Secondary:\quad \thesecondaryclass

\vglue 5pt plus 0.3fil minus 2pt

{\bf Keywords:}\quad \thekeywords

\vglue 10pt plus 0.5fil minus 5pt

{\small  Proposed: \theproposer\hfill Received: \receiveddate\nl
Seconded: \theseconders\hfill 
\ifx\reviseddate\relax                         
Accepted: \accepteddate                        
\else
Revised: \reviseddate                          
\fi}
\eject
}       

\let\maketitlepage\maketitlep
\let\maketitle\maketitlepage

\font\phead=cmsl9 scaled 950
\font=cmsl9 scaled 1050
\font\pnum=cmbx10 scaled 913
\font\lnum=cmbx10 
\font\pfoot=cmsl9 scaled 950
\font=cmsl9 scaled 1050
\ifplaintex
\headline{\vbox to 0pt{\vskip -4.5mm\line{\small\phead\ifnum
\count0=\startpage ISSN 1364-0380 (on line)
1465-3060 (printed) \hfill {\pnum\folio}\else\ifodd\count0\def\\{ }%
\ifx\theshorttitle\relax\thetitle\else\theshorttitle\fi\hfill{\pnum\folio}
\else\def\\{ and }{\pnum\folio}\hfill\ifx\theshortauthors\relax\theauthors
\else\theshortauthors\fi\fi\fi}\vss}}
\footline{\vbox to 0pt{\vglue 0mm\line{\small\pfoot\ifnum\count0=\startpage
\copyright\ \gtp\hfill\else
\gt, Volume \thevolumenumber\ (\thevolumeyear)\hfill\fi}\vss
}}
\else
\makeatletter
\def\@oddhead{{\small\lhead\ifnum\count0=\startpage ISSN 1364-0380 (on line)
1465-3060 (printed) \hfill {\lnum\number\count0}\else\ifodd\count0
\def\\{ }\ifx\theshorttitle\relax \thetitle \else\theshorttitle\fi\hfill
{\lnum\number\count0}\else\def\\{ and }{\lnum\number\count0}
\hfill\ifx\theshortauthors\relax 
\theauthors\else\theshortauthors\fi\fi\fi}}\def\@evenhead{@oddhead}
\def\@oddfoot{\small\lfoot\ifnum\count0=\startpage\copyright\ \gtp\hfill\else
\gt, Volume \thevolumenumber\ (\thevolumeyear)\hfill\fi}
\def\@evenfoot{@oddfoot}
\makeatother
\fi

\newwrite\gtoutfile
\long\gdef\makeheadfile{  
{\def\\{, }\def\s{ }
\immediate\openout\gtoutfile head.xxx
\immediate\write\gtoutfile{To: math@arxiv.org}
\immediate\write\gtoutfile{Subject: put OR rep NNNNN:ppppp}
\immediate\write\gtoutfile{--text follows this line--}
\immediate\write\gtoutfile{Proxy-for: \ifx\theasciiauthors\relax
\theauthors\else\theasciiauthors\fi\s<\ifx\theasciiemail\relax\theemail\else\theasciiemail\fi>}
\immediate\write\gtoutfile{\noexpand\\}
\immediate\write\gtoutfile{Authors: \ifx\theasciiauthors\relax
\theauthors\else\theasciiauthors\fi}
{\def\\{ }\immediate\write\gtoutfile{Title: \ifx\theasciititle\relax
\thetitle\else\theasciititle\fi}}
\immediate\write\gtoutfile{Subj-class: GT or SG, GR etc}
\immediate\write\gtoutfile{MSC-class: \theprimaryclass\ifx\thesecondaryclass\relax\else, \thesecondaryclass\fi}
\immediate\write\gtoutfile{Journal-ref: Geom. Topol. \thevolumenumber\s
(\thevolumeyear) \startpage-\finishpage}
\immediate\write\gtoutfile{Comments: Published in Geometry and Topology at}
\immediate\write\gtoutfile{\s\s\s  http://www.maths.warwick.ac.uk/gt/GTVol\thevolumenumber/paper\thepapernumber.abs.html}
\immediate\write\gtoutfile{\noexpand\\}
\immediate\write\gtoutfile{}
\ifx\theasciiabstract\relax
\immediate\write\gtoutfile{\theabstract}\else
\immediate\write\gtoutfile{\theasciiabstract}\fi
\immediate\write\gtoutfile{}
\immediate\write\gtoutfile{\noexpand\\}
\immediate\write\gtoutfile{}
\immediate\write\gtoutfile{}
\immediate\closeout\gtoutfile}}  

\def\maketitlepage{\maketitlep\makeheadfile}
\let\maketitle\maketitlepage

\def\ifplaintex{\expandafter\ifx\csname documentclass\endcsname\relax}

\ifplaintex 
\else
\fi

\expandafter\ifx\csname beginpicture\endcsname\relax
\expandafter\ifx\csname documentclass\endcsname\relax
\input pictex \else
\input prepictex \input pictex \input postpictex \fi\fi

\def\gt{{\mathsurround=0pt\it $\cal G\mskip-2mu$eometry \&\ 
$\cal T\!\!$opology}}        

\def\gtp{{\mathsurround=0pt\it $\cal G\mskip-2mu$eometry \&\ 
$\cal T\!\!$opology $\cal P\!$ublications}}  

\def\volumenumber#1{\def\thevolumenumber{#1}}
\def\papernumber#1{\def\thepapernumber{#1}}
\def\volumeyear#1{\def\thevolumeyear{#1}}

\def\pagenumbers#1#2{\def\startpage{#1}\def\finishpage{#2}}
\def\published#1{\def\publishdate{#1}}
\def\proposed#1{\def\theproposer{#1}}
\def\seconded#1{\def\theseconders{#1}}
\def\received#1{\def\receiveddate{#1}}

\def\accepted#1{\def\accepteddate{#1}}

\long\def\asciiabstract#1{\long\def\theasciiabstract{#1}}

\let\\\par\let\thevolumenumber\relax
\let\thepapernumber\relax
\let\thevolumeyear\relax\let\thesamplenumber\relax\let\startpage\relax
\let\finishpage\relax\let\publishdate\relax\let\receiveddate\relax
\let\reviseddate\relax\let\accepteddate\relax\let\theasciititle\relax
\let\theasciiauthors\relax
\let\theasciiabstract\relax
\let\theasciiemail\relax\let\theshortauthors\relax\let\theshorttitle\relax

\long\def\maketitlep{   

\count0=\startpage

\gt\hfill      
\beginpicture
\setcoordinatesystem units <0.33truein, 0.33truein> point at 2.2 0.9
\setplotsymbol ({$\cal G$})
\plotsymbolspacing=9truept
\circulararc 315 degrees from 0 1 center at 0 0
\setplotsymbol ({$\cal T$})
\circulararc 315 degrees from 1 -1 center at 1 0
\endpicture
\break
{\small\ifx\thesamplenumber\relax 
Volume \else Sample
\fi\thevolumenumber\ (\thevolumeyear)
\startpage--\finishpage\nl
Published: \publishdate}
\vglue 0.5truein plus 0.4fil minus 0.1truein

{\parskip=0pt\leftskip 0pt plus 1fil\def\\{\par\smallskip}{\ifplaintex\large
\else\Large\fi\bf\thetitle}\par\medskip}   

\vglue 0pt plus 0.1fil 

{\parskip=0pt\leftskip 0pt plus 1fil\def\\{\par}{\sc\theauthors}
\par\medskip}

\vglue 0pt plus 0.1fil 

{\small\parskip=0pt\let\newline\\
{\leftskip 0pt plus 1fil\def\\{\par}{\sl\theaddress}\par}
\expandafter\ifx\theemail\relax    
\relax\else\vglue 5pt plus 0.02fil minus 2pt\def\\{\stdspace{\rm 
and}\stdspace} 
\cl{Email:\stdspace\tt\theemail}\fi
\ifx\theurl\relax                  
\relax\else\vglue 5pt plus 0.02fil minus 2pt\def\\{\stdspace{\rm 
and}\stdspace}
\cl{URL:\stdspace\tt\theurl}\fi\par}

\vglue 7pt plus 0.3fil minus 3pt

{\bf Abstract}
\vglue 5pt plus 0.1fil minus 2pt

\theabstract

\vglue 7pt plus 0.3fil minus 3pt

{\bf AMS Classification numbers}\quad Primary:\quad \theprimaryclass

Secondary:\quad \thesecondaryclass

\vglue 5pt plus 0.3fil minus 2pt

{\bf Keywords:}\quad \thekeywords

\vglue 10pt plus 0.5fil minus 5pt

{\small  Proposed: \theproposer\hfill Received: \receiveddate\nl
Seconded: \theseconders\hfill 
\ifx\reviseddate\relax                         
Accepted: \accepteddate                        
\else
Revised: \reviseddate                          
\fi}
\eject
}       

\let\maketitlepage\maketitlep
\let\maketitle\maketitlepage

\font\phead=cmsl9 scaled 950
\font=cmsl9 scaled 1050
\font\pnum=cmbx10 scaled 913
\font\lnum=cmbx10 
\font\pfoot=cmsl9 scaled 950
\font=cmsl9 scaled 1050
\ifplaintex
\headline{\vbox to 0pt{\vskip -4.5mm\line{\small\phead\ifnum
\count0=\startpage ISSN 1364-0380 (on line)
1465-3060 (printed) \hfill {\pnum\folio}\else\ifodd\count0\def\\{ }%
\ifx\theshorttitle\relax\thetitle\else\theshorttitle\fi\hfill{\pnum\folio}
\else\def\\{ and }{\pnum\folio}\hfill\ifx\theshortauthors\relax\theauthors
\else\theshortauthors\fi\fi\fi}\vss}}
\footline{\vbox to 0pt{\vglue 0mm\line{\small\pfoot\ifnum\count0=\startpage
\copyright\ \gtp\hfill\else
\gt, Volume \thevolumenumber\ (\thevolumeyear)\hfill\fi}\vss
}}
\else
\makeatletter
\def\@oddhead{{\small\lhead\ifnum\count0=\startpage ISSN 1364-0380 (on line)
1465-3060 (printed) \hfill {\lnum\number\count0}\else\ifodd\count0
\def\\{ }\ifx\theshorttitle\relax \thetitle \else\theshorttitle\fi\hfill
{\lnum\number\count0}\else\def\\{ and }{\lnum\number\count0}
\hfill\ifx\theshortauthors\relax 
\theauthors\else\theshortauthors\fi\fi\fi}}\def\@evenhead{@oddhead}
\def\@oddfoot{\small\lfoot\ifnum\count0=\startpage\copyright\ \gtp\hfill\else
\gt, Volume \thevolumenumber\ (\thevolumeyear)\hfill\fi}
\def\@evenfoot{@oddfoot}
\makeatother
\fi

\newwrite\gtoutfile
\long\gdef\makeheadfile{  
{\def\\{, }\def\s{ }
\immediate\openout\gtoutfile head.xxx
\immediate\write\gtoutfile{To: math@arxiv.org}
\immediate\write\gtoutfile{Subject: put OR rep NNNNN:ppppp}
\immediate\write\gtoutfile{--text follows this line--}
\immediate\write\gtoutfile{Proxy-for: \ifx\theasciiauthors\relax
\theauthors\else\theasciiauthors\fi\s<\ifx\theasciiemail\relax\theemail\else\theasciiemail\fi>}
\immediate\write\gtoutfile{\noexpand\\}
\immediate\write\gtoutfile{Authors: \ifx\theasciiauthors\relax
\theauthors\else\theasciiauthors\fi}
{\def\\{ }\immediate\write\gtoutfile{Title: \ifx\theasciititle\relax
\thetitle\else\theasciititle\fi}}
\immediate\write\gtoutfile{Subj-class: GT or SG, GR etc}
\immediate\write\gtoutfile{MSC-class: \theprimaryclass\ifx\thesecondaryclass\relax\else, \thesecondaryclass\fi}
\immediate\write\gtoutfile{Journal-ref: Geom. Topol. \thevolumenumber\s
(\thevolumeyear) \startpage-\finishpage}
\immediate\write\gtoutfile{Comments: Published in Geometry and Topology at}
\immediate\write\gtoutfile{\s\s\s  http://www.maths.warwick.ac.uk/gt/GTVol\thevolumenumber/paper\thepapernumber.abs.html}
\immediate\write\gtoutfile{\noexpand\\}
\immediate\write\gtoutfile{}
\ifx\theasciiabstract\relax
\immediate\write\gtoutfile{\theabstract}\else
\immediate\write\gtoutfile{\theasciiabstract}\fi
\immediate\write\gtoutfile{}
\immediate\write\gtoutfile{\noexpand\\}
\immediate\write\gtoutfile{}
\immediate\write\gtoutfile{}
\immediate\closeout\gtoutfile}}  

\def\maketitlepage{\maketitlep\makeheadfile}
\let\maketitle\maketitlepage

\volumenumber{4}\papernumber{11}\volumeyear{2000}
\pagenumbers{309}{368}
\proposed{Yasha Eliashberg}
\seconded{Tomasz Mrowka, Joan Birman}
\received{6 October 2000}
\accepted{14 October 2000}
\published{14 October 2000}

\usepackage{epsf,amssymb,amsmath,psfrag}

\newtheorem{thm}{Theorem}[section]
\newtheorem{prop}[thm]{Proposition}

\newtheorem{lemma}[thm]{Lemma}
\newtheorem{cor}[thm]{Corollary}

\newtheorem{quest}{Question}

\newtheorem*{claim}{Claim}
\newtheorem{p}[thm]{Problem}

\newcommand{\C}{\mbox{\bf{C}}}
\newcommand{\R}{\mbox{\bf{R}}}
\newcommand{\Z}{\mbox{\bf{Z}}}
\newcommand{\Q}{\mbox{\bf{Q}}}
\renewcommand{\H}{\mathbb{H}}
\renewcommand{\P}{\mbox{\bf{P}}}

\newcommand{\bdry}{\partial}

\newcommand{\be}{\begin{enumerate}}
\newcommand{\ee}{\end{enumerate}}
\newcommand{\sa}{\rightsquigarrow}

\begin{document}

\title{On the classification of tight contact structures I}

\author{Ko Honda}
\address{Mathematics Department, University of Georgia\\Athens, GA 30602, USA}
\email{honda@math.uga.edu}
\url{http://www.math.uga.edu/\char126 honda}

\begin{abstract}
We develop new techniques in the theory of convex surfaces to prove
complete classification results for tight contact structures on lens
spaces, solid tori, and $T^2\times I$.
\end{abstract}

\asciiabstract{We develop new techniques in the theory of convex
surfaces to prove complete classification results for tight contact
structures on lens spaces, solid tori, and T^2 X I$.}

\keywords{Tight, contact structure, lens spaces, solid tori}

\primaryclass{57M50}

\secondaryclass{53C15}

\maketitlepage

\section{Introduction}

It has been known for some time that, in dimension 3, contact structures fall into one of two classes:
tight or overtwisted.  A contact structure $\xi$ is said to be {\it overtwisted} if there exists an embedded disk
$D$ which is tangent to $\xi$ everywhere along $\bdry D$, and a contact structure is {\it tight} if it is not overtwisted.
This dichotomy was first discovered by Bennequin in his seminal paper \cite{B}, and
further elucidated by Eliashberg \cite{E}.  In \cite{E89}, Eliashberg classified overtwisted contact structures on
closed 3--manifolds, effectively reducing the overtwisted classification to a homotopy classification of
2--plane fields on 3--manifolds.  Eliashberg \cite{E} then proceeded to classify tight contact structures on
the 3--ball $B^3$, the 3--sphere $S^3$, $S^2\times S^1$, and $\R^3$.  In particular, he proved that there exists a unique
tight contact structure on $B^3$, given a fixed boundary characteristic foliation --- this theorem of Eliashberg
comprises the foundational building block in the study of tight contact structures on
3--manifolds.  Subsequent results on the classification
of tight contact structures were:  a complete classification on the 3--torus by Kanda \cite{K} and Giroux (obtained
independently), a complete classification on some lens spaces by Etnyre \cite{Et}, and some partial results
on solid tori $S^1\times D^2$ by Makar--Limanov \cite{ML} and circle bundles over Riemann surfaces by
Giroux.  One remarkable discovery by Makar--Limanov \cite{ML} was that there exist tight contact structures which
become overtwisted when pulled back to the universal cover $\widetilde{M}$ via the covering map
$\pi\co \widetilde{M}\rightarrow M$.  This prompts us to define a {\it universally tight} contact structure to be one
which remains tight when pulled back to $\widetilde{M}$ via $\pi$.  We call a tight contact structure $\xi$
{\it virtually overtwisted} if $\xi$ becomes overtwisted when pulled back to a {\it finite} cover.  It is not
known whether every tight contact structure is either universally tight or virtually overtwisted, although this
dichotomy holds when $\pi_1(M)$ is residually finite.

The goal of this paper is to give a complete classification of tight contact structures
on lens spaces, as well as a complete classification of tight contact structures on
solid tori $S^1\times D^2$ and toric annuli $T^2\times I$ with convex
boundary.  This completes
the classification of tight contact structures on lens spaces, initiated by
Etnyre in \cite{Et}, as well as the classification of tight contact structures on solid
tori (at least for convex boundary), initiated by Makar--Limanov
\cite{ML}.   We will also determine precisely which tight contact structures are universally tight and which are
virtually overtwisted  --- all the manifolds we consider this paper will have residually finite
$\pi_1(M)$, hence tight contact structures on these manifolds
will either be universally tight or virtually overtwisted.   Our method is
a systematic application of the methods developed by Kanda \cite{K}, which in turn use Giroux's
theory of convex surfaces \cite{Giroux2}.
In essence, we use Kanda's methods and apply them in Etnyre's setting: we decompose the 3--manifold
$M$ in a series of steps, along {\it closed convex surfaces} or {\it convex surfaces with Legendrian boundary}.
The difference between Etnyre's approach and ours
is that we require that the cutting surfaces have boundary consisting of  {\it Legendrian curves},
whereas Etnyre used cutting surfaces which had
{\it transverse curves} on the boundary.
The Legendrian curve approach appears to be more efficient and yields fewer
possible configurations than the transverse curve approach, although the author is not quite sure why this is the
case.

The classification theorems will reveal a closer connection between contact structures and 3--dimensional
topology than was previously expected.  In particular, the geometry of $\pi_0(\mbox{Diff}^+(T^2))=SL(2,\Z)$ (including the standard
Farey tessellation) plays a significant role for the 3--manifolds studied in this paper --- lens spaces have
Heegaard decompositions into solid tori, and the toric annulus contains incompressible $T^2$.
Unlike foliation theory (which is related to contact topology by the work of Eliashberg and Thurston \cite{ET}),
contact topology has a built-in `handedness', and we will see that the contact topology is determined in large
part by {\it positive Dehn twists} in $\pi_0({\mbox{Diff}}^+(T^2))=SL(2,\Z)$.   We believe the results in this paper
represent a tiny fraction of a large and emerging theory of contact structures applied to three--manifold topology.
The techniques developed in this paper are applied to other classes of 3--manifolds (circle bundles which
fiber over closed oriented surfaces and torus bundles over $S^1$) in the sequel \cite{H2}, and
in \cite{EH99} J. Etnyre and the author prove the non-existence of positive tight contact structures on the
Poincar\'e homology sphere for one of its orientations, thereby producing the first example of a
closed 3--manifold
which does not carry a tight contact structure.

\medskip
{\bf Note}\qua  E Giroux has independently obtained similar
classification results.   His approach and ours are surprisingly dissimilar, and the interested reader
will certainly increase his understanding by reading his account \cite{Giroux3} as well.

The first version of this paper was written in April 1999.  This
version was written on August 1, 2000.

\section{Statements of results}

In this paper all the 3--manifolds $M$ are oriented and compact,
and all the contact structures $\xi$ are positive, ie, given by a global 1--form $\alpha$
with $\alpha\wedge d\alpha>0$, and oriented.   We will simply write `contact structure', when we
mean `positive, oriented
contact structure'.

\subsection{Lens spaces}
Consider the lens space $L(p,q)$, where $p>q>0$ and $(p,q) = 1$. Assume
$-{p\over
q}$ has the continued fraction expansion
$$-{p\over q}=r_0-{1\over r_1-{1\over r_2-\cdots {1\over r_k}}},$$
with all $r_i<-1$.
 Then we have the following classification theorem for tight contact
structures on lens spaces $L(p,q)$.

\begin{thm} \label{1}
There exist exactly $|(r_0+1)(r_1+1)\cdots (r_k+1)|$
tight contact structures on the lens space $L(p,q)$ up to isotopy, where
$r_0,...,r_k$ are
the coefficients of the continued fraction expansion of $-{p\over q}$.
Moreover, all the tight contact structures on $L(p,q)$ can be obtained from
Legendrian surgery on links in $S^3$, and are therefore holomorphically
fillable.
\end{thm}

{\it Legendrian surgery} is a contact surgery technique due to Eliashberg \cite{E90}.  It
produces contact structures which are holomorphically fillable, and are therefore tight, by
a result of Eliashberg and Gromov
\cite{E91,Gromov}.


\subsection{The thickened torus $T^2\times I$}

When we study contact structures on manifolds
with boundary, we need to impose a boundary condition ---
a natural condition would be to ask that the boundary be {\it convex}.
A closed, oriented, embedded surface $\Sigma$ in a contact manifold
$(M,\xi)$ is said to be {\it convex}
if there is a vector field $v$ transverse to $\Sigma$ whose
flow preserves $\xi$.  A generic surface $\Sigma$ inside a contact 3--manifold is convex \cite{Giroux2},
so demanding
that the boundary be convex presents no loss of generality.

A convex surface $\Sigma\subset (M,\xi)$
has a naturally associated family of disjoint embedded
curves $\Gamma_\Sigma$, well-defined up to isotopy and called the {\it dividing curves} (for more details see
Section \ref{section:dividing}).
The dividing curves $\Gamma_\Sigma$ separate the surface $\Sigma$ into two
subsurfaces $R_+$ and $R_-$.  If $\xi$ is tight and $\Sigma\not=S^2$, then
the dividing curves $\Gamma_\Sigma$ are homotopically essential, in the sense that
none of them bounds an embedded disk in $\Sigma$.
In particular, if $\Sigma$ is a torus,
$\Gamma_\Sigma$ will consist of an even number of parallel essential
curves.

Consider a tight contact structure $\xi$ on $T^2\times I = T^2 \times [0, 1]$ with convex boundary.
Fix an oriented identification between the torus $T^2$ and $\R^2/\Z^2$.
Given a convex torus $T$ in $T^2\times I$, its set of dividing curves is, up to isotopy,
determined by the following data: (1) the {\it number} $\#\Gamma_T$ of these dividing curves
and (2) their {\it slope} $s(T)$, defined by the property that each curve is isotopic to a linear
curve of slope $s(T)$ in $T\simeq \R^2/\Z^2$.

\subsubsection{Twisting}    \label{section:twisting}

In order to state the classification theorem for $T^2\times I$ it is necessary to define the notions of
{\it twisting in the $I$--direction}, {\it minimal twisting in the $I$--direction}, and
{\it nonrotativity in the $I$--direction}.

Given a slope $s$ of a line in $\R^2$ (or $\R^2/\Z^2$), associate to it its standard angle
$\overline\alpha(s)\in \R\P^1=\R/\pi\Z$. For $\overline{\alpha}_1$, $\overline{\alpha}_2\in
\R\P^1$, let $[\overline{\alpha}_1, \overline{\alpha}_2]$ be the image of the
interval $[\alpha_1,\alpha_2]\subset \R$, where $\alpha_i\in\R$ are
representatives of $\overline{\alpha}_i$ and $\alpha_1\leq \alpha_2<\alpha_1+\pi$.
A slope $s$ is said to be {\it between} $s_1$ and $s_0$ if $\overline{\alpha}(s)
\in [\overline{\alpha}(s_1),\overline{\alpha}(s_0)]$.

Consider a tight contact structure $\xi$ on $T^2\times I$ with convex boundary
and boundary slopes $s_i=s(T_i)$,
$i=0,1$, where   $T_i=T^2\times \{i\}$.
We say $\xi$ is {\it minimally twisting} (in the $I$--direction)
if every convex torus parallel to the boundary
has slope $s$ between $s_1$ and $s_0$.  In particular, $\xi$ is {\it nonrotative} (in the $I$--direction) if
$s_1=s_0$ and $\xi$ is minimally twisting.
Define the {\it $I$--twisting} of a tight $\xi$ to be
$\beta_I=\alpha(s_0)-\alpha(s_1)=
\sum_{k=1}^l(\alpha(s_{k-1\over l})-\alpha(s_{k\over l}))$, where (i)  $s_{k\over l}=s(T_{k\over l})$,
$k=0,\cdots, l$, (ii) $T_0=T^2\times\{0\}$, $T_1=T^2\times\{1\}$,
and $T_{k\over l}$, $k=1,\cdots, l-1$ are mutually disjoint  convex tori parallel to the boundary,
arranged in order
from closest to $T_0$ to farthest from $T_0$, (iii) $\xi$ is minimally twisting
between $T_{k-1\over l}$ and $T_{k\over l}$, and (iv) $\alpha(s_{k\over l})\leq
\alpha(s_{k-1\over l})<\alpha(s_{k\over l})+\pi$.

The following will be shown in Proposition \ref{hiho}:
\be
\item The $I$--twisting of $\xi$ is  well-defined, finite,
and independent of the choices of $l$ and the $T_{k\over l}$.
\item The $I$--twisting of $\xi$ is always non-negative.
\ee

Notice that the $I$--twisting $\beta_I$ is dependent
on the particular identification $T^2=\R^2/\Z^2$.  We therefore introduce
$\phi_I(\xi)=\pi \lfloor{\beta_I \over \pi}\rfloor$,
which is independent of the identification.  Here $\lfloor \cdot\rfloor$ is the greatest integer function. Also,
$\phi_I=0$ is equivalent to minimal  twisting.

\subsubsection{Statement of theorem}
After normalizing via $\pi_0(\mbox{Diff}^+(T^2))=SL(2,\Z)$, we may assume that $T_1$
has dividing curves with slope $-{p\over q}$, where $p\geq q>0$, $(p, q) =
1$,  and $T_0$ has slope $-1$.     Denote $T_a=T^2\times \{a\}$.
For this boundary data, we have the following:

\begin{thm} \label{2} Consider $T^2\times I$ with convex boundary,
and assume, after normalizing via $SL(2,\Z)$,
that $\Gamma_{T_1}$ has slope $-{p\over q}$, and $\Gamma_{T_0}$ has slope
$-1$.   Assume we fix a characteristic foliation on $T_0$ and $T_1$ with these
dividing curves.  Then, up to an isotopy which fixes the boundary, we have the following classification:
\be
\item Assume either {\rm(a)} $-{p\over q}<-1$ or {\rm(b)} $-{p\over q}=-1$ and $\phi_I>0$.
Then there exists a unique factorization $T^2\times I =
(T^2\times[0,{1\over 3}])\cup (T^2\times[{1\over 3},{2\over 3}])\cup (T^2\times[{2\over 3},1])$, where
(i) $T_{i\over 3}$,  $i=0,1,2,3$, are convex, (2) $(T^2\times[0,{1\over 3}])$ and $(T^2\times[{2\over 3},1])$
are nonrotative, (3)  $\#\Gamma_{T_{1\over 3}}=\#\Gamma_{T_{2\over 3}}=2$, and (4) $T_{1\over 3}$ and
$T_{2\over 3}$ have fixed characteristic foliations which are adapted to $\Gamma_{T_{1\over 3}}$
and $\Gamma_{T_{2\over 3}}$.
\item Assume $-{p\over q}<-1$ and $\#\Gamma_{T_0}=\#\Gamma_{T_1}=2$.
\be
\item[\rm(a)] There exist exactly $|(r_0+1)(r_1+1)\cdots (r_{k-1}+1)(r_k)|$
tight contact structures with $\phi_I=0$.
Here, $r_0,...,r_k$ are the coefficients of the continued fraction
expansion of $-{p\over q}$, and $-{p\over q}<-1$.
\item[\rm(b)]   There exist exactly 2 tight contact
structures with $\phi_I=n$, for each $n\in \Z^+$.
\ee
\item Assume $-{p\over q}=-1$ and $\#\Gamma_{T_0}=\#\Gamma_{T_1}=2$.
Then there exist exactly 2 tight contact
structures with $\phi_I=n$, for each $n\in \Z^+$.
\item Assume $-{p\over q}=-1$ and $\#\Gamma_{T_0}=2n_0$, $\#\Gamma_{T_1}=2n_1$.
Then the nonrotative tight contact structures are in 1--1 correspondence with $\mathcal{G}$, the
set of all possible (isotopy classes of) configurations of arcs on an annulus $A=S^1\times I$ with markings
$\sigma_i\subset S^1\times \{i\}$, $i=0,1$, which
satisfy the following:
\be
\item[\rm(a)] $|\sigma_i|=2n_i$, $i=0,1$, where $|\cdot|$ denotes cardinality.
\item[\rm(b)] Every point of $\sigma_0\cup \sigma_1$ is precisely one endpoint of one arc.
\item[\rm(c)] There exist at least two arcs which begin on $\sigma_0$ and end on $\sigma_1$.
\item[\rm(d)] There are no closed curves.
\ee
\ee
\end{thm}

\subsection{Solid tori}
Finally, we have the analogous theorem for solid tori.
Fix an oriented identification of $T^2=\bdry (S^1\times D^2)$ with $\R^2/\Z^2$,
where $\pm(1,0)^T$ corresponds to the meridian of the solid torus, and $\pm (0,1)^T$ corresponds
the longitudinal direction determined by a chosen framing.
We consider tight contact structures $\xi$ on $S^1\times D^2$ with convex boundary $T^2$.
Let the
{\it slope} $s(T^2)$ of $T^2$ be the slope under the identification $T^2\simeq \R^2/\Z^2$.

\begin{thm} \label{3}   Consider the tight contact structures on $S^1\times D^2$ with convex
boundary $T^2$, for which $\#\Gamma_{T^2}=2$ and $s(T^2)=-{p\over q}$, $p\geq q>0,
(p, q) = 1$.   Fix a characteristic foliation $\mathcal{F}$ which is adapted to $\Gamma_{T^2}$.
There exist exactly $|(r_0+1)(r_1+1)\cdots (r_{k-1}+1)(r_k)|$  tight
contact structures on $S^1\times D^2$ with this boundary condition, up to
isotopy fixing $T^2$.
Here, $r_0,\cdots, r_k$ are the coefficients of the continued fraction
expansion of $-{p\over q}$.
\end{thm}

In other words, the number of tight contact structures
for the solid torus with (a fixed) convex boundary
with $\#\Gamma_{T^2}=2$ and $s(T^2)=-{p\over q}$ is the same as the number of tight contact
structures on
$T^2\times I$ with (fixed) convex boundary, $\#\Gamma_{T_i}=2$, $i=0,1$,
slopes $s(T_1)=-{p\over q}$ and $s(T_0)=-1$, and minimal twisting.

Via a multiplication by
$\left (\begin{array}{cc}
1 & m\\ 0 & 1 \end{array}\right )\in SL(2,\Z)$, $m\in \Z$, which is equivalent to a
Dehn twist which induces a change of framing,
all the boundaries of $S^1\times D^2$ can be put in the form described in the
theorem above.   In addition, the choice of slope $-{p\over q}$ with $p \ge
q > 0$ is unique.

\subsection{Strategy of proof}
First consider $T^2\times I$.  We fix a boundary condition by prescribing dividing sets
$\Gamma_i=\Gamma_{T_i}$, $i=0,1$.   Also fix a boundary characteristic foliation which is compatible
with $\Gamma_i$.  Giroux's Flexibility Theorem, described in Section \ref{convexity}, roughly states that it is the
{\it isotopy type} of the dividing set $\Gamma$ which dictates the geometry of $\Sigma$, not the
precise characteristic foliation which is compatible with $\Gamma$.  This allows us to reduce the classification
to one particular characteristic foliation compatible with $\Gamma_i$, and we choose a (rather non-generic)
realization of a convex surface --- one that is in {\it standard form} (see Section \ref{section:tori}).

In Section \ref{by} we introduce the notion of a {\it bypass}, which is the crucial new ingredient which allows us
to successively peel off `thin' $T^2\times I$ layers which we call {\it basic slices}.  We eventually obtain a
factorization of a $(T^2\times I,\xi)$ into basic $T^2\times I$ slices, if $\xi$ is tight and
minimally twisting.  This decomposition gives
a possible upper bound for the number of tight contact structures on $T^2\times I$ with
given boundary conditions.  These candidate tight contact structures are easily distinguished by the
relative Euler class.  We then successively embed $T^2\times I\subset S^1\times D^2\subset L(p,q)$,
and find that the upper bound is exact, since all of the candidate tight contact structures can be
realized by Legendrian surgery.  The remaining cases of Theorem \ref{2} when the $I$--twisting is not minimal
and when $\#\Gamma_i>2$ are treated in Section \ref{whatever}.

\section{Preliminaries}

\subsection{Convexity}   \label{convexity}

In this section only $(M,\xi)$ is a compact, oriented 3--manifold with
a contact structure, tight or overtwisted.

An oriented properly embedded surface $\Sigma$ in $(M,\xi)$ is called {\it convex} if
there is a vector field $v$ transverse to $\Sigma$ whose flow preserves $\xi.$  This {\it contact
vector field} $v$ allow us to find an $I$--invariant neighborhood $\Sigma\times I\subset M$ of $\Sigma$,
where $\Sigma=\Sigma\times \{0\}$.
In most applications, our convex surface $\Sigma$ will either be closed or compact with {\it Legendrian boundary}.
The theory of closed convex surfaces appears in detail in Giroux's paper \cite{Giroux2}.
However, the same results for the Legendrian boundary case have not appeared in the literature,
and we will rederive Giroux's results in this case.

\subsubsection{Twisting number of  a Legendrian curve}
A curve $\gamma$ which is everywhere tangent to $\xi$ is called {\it Legendrian}.
We define the {\it twisting number}
$t(\gamma,Fr)$ of a closed Legendrian curve $\gamma$
with respect to a given framing $Fr$ to be the number of counterclockwise
(right) $2\pi$ twists
of $\xi$ along $\gamma$, relative to $Fr$.
In particular, if $\gamma$ is a connected component of the boundary of a compact surface $\Sigma$,
$T\Sigma$ gives
a natural framing $Fr_\Sigma$, and if $\Sigma$ is a Seifert surface of $\gamma$, then
$t(\gamma,Fr_\Sigma)$ is the Thurston--Bennequin invariant
$tb(\gamma)$.  We will
often suppress $Fr$ when the framing is understood.
Notice that it is easy to
decrease $t(\gamma,Fr)$ by locally adding zigzags in a front projection, but not always possible to increase
$t(\gamma,Fr)$.

\subsubsection{Perturbation into a convex surface with Legendrian boundary}
Giroux \cite{Giroux2}
proved that a closed oriented embedded surface $\Sigma$ can be deformed
by a $C^\infty$--small isotopy so that the resulting embedded surface is convex.
We will prove the following proposition:

\begin{prop}  \label{prop:makeconvex}
Let $\Sigma\subset M$ be a compact, oriented, properly embedded surface with Legendrian boundary,
and assume $t(\gamma,Fr_\Sigma)\leq 0$ for all components $\gamma$ of $\bdry \Sigma$.  There exists a
$C^0$--small perturbation near the boundary (fixing $\bdry\Sigma$)
which puts an annular neighborhood $A$ of $\bdry\Sigma$ into
a standard form, and a subsequent $C^\infty$--small perturbation of the perturbed surface (fixing the annular
neighborhood of $\bdry\Sigma$), which makes $\Sigma$ convex.  Moreover, if $v$ is a contact vector field
defined on a neighborhood of $A$ and transverse to $A\subset \Sigma$, then $v$ can be extended to a
contact vector field transverse to all of $\Sigma$.
\end{prop}

\proof
Assume  that $t(\gamma,Fr_\Sigma) < 0$, for all boundary components
$\gamma$. After a $C^0$--small perturbation near the boundary (fixing the boundary), we may assume that
$\gamma$ has a {\it standard annular collar} $A$.  Here  $A=S^1\times [0,1] =({\bf R}/{\bf Z}) \times [0,1]$
with coordinates $(x,y)$ and  $\gamma=S^1\times\{0\}$.
Its neighborhood $A\times [-1,1]$ has coordinates $(x,y,t)$,  and the contact 1--form on $A\times[-1,1]$ is
$\alpha=\sin(2\pi nx)dy+\cos(2\pi nx)dt$.
The Legendrian curves $S^1\times \{\mbox{pt}\}\subset A$ are called the {\it Legendrian rulings} and
and $\{{k \over 2n}\}\times [0,1]\subset A$, $k=1,2,\cdots, 2n$ are called the {\it Legendrian divides}.

\begin{figure}[ht!]
{\epsfysize=2in\centerline{\epsfbox{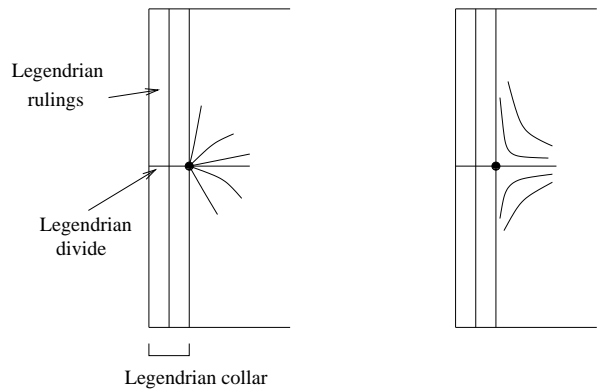}}}
\caption{Half-elliptic point and half-hyperbolic point}
\label{half}
\end{figure}

Once we have standard annular neighborhoods of $\bdry\Sigma$, we
use the following perturbation lemma, due to Fraser \cite{F} --- refer to Figure \ref{half}
for an illustration of {\it half-elliptic} and {\it half-hyperbolic} singular
points.

\begin{lemma}\label{perturb} It is possible to perturb $\Sigma$, while
fixing the
Legendrian collar, to make any tangency $({k\over 2n}, 1) \in A=({\bf
R}/{\bf Z})
\times [0,1]\subset \Sigma$ half-elliptic and any tangency half-hyperbolic.
\end{lemma}

\begin{proof}
It suffices, by a Darboux-type argument,
to extend the contact structure on $S^1\times [0,1]\times[-1,1]$  above
to $S^1\times [0,2]\times[-1,1]$, such that the characteristic foliation on $S^1\times[0,2]\times\{0\}$
has a half-elliptic or a half-hyperbolic singularity.  It therefore also suffices to treat the neighborhood of
a Legendrian divide.  Without loss of generality, let the Legendrian divide be $\{0\}\times[0,1]\times\{0\}
\subset [-\varepsilon,\varepsilon]\times
[0,1]\times\{0\}$, with contact 1--form $\alpha'=dt+xdy$.  Now extend to $\alpha'
=dt-f(y)dx+xdy$ for a half-elliptic singularity and $\alpha'=
dt+f(y)dx+xdy$ for a half-hyperbolic singularity, on $[-\varepsilon,\varepsilon]\times
[0,2]\times\{0\}$, where $f(y)=0$ on $[0,1]$ and ${df\over dy}>0$ on $[1,2]$.
\end{proof}

{\bf Note}\qua  M Fraser \cite{F} obtained normal forms near the boundary, for $\Sigma$ with
Legendrian boundary, even when $t(\gamma) > 0$ for some boundary component
$\gamma$.  In this case, Lemma \ref{perturb} is no longer applicable.  Instead, all the singularities
must be half-hyperbolic, after appropriate cancellations.  If $t(\gamma) > 0$,
$\Sigma$ cannot be made convex.
\medskip

When $t(\gamma)=0$, then perturb $\Sigma$, fixing $\gamma$, so that
the contact structure is given by $\alpha=dt-ydx$ on $A\times [-1,1]$, where $A$ is as before.

If $\Sigma$ is compact with Legendrian boundary, and all the boundary
components have $t \le 0$,
we use Lemma \ref{perturb} if $t< 0$, to make all the boundary tangencies of $\Sigma$
half-elliptic (if $t=0$ use the paragraph above), and perturb to
obtain $\Sigma$ with characteristic foliation $\mathcal{F}$ which is Morse--Smale on the
interior.  This means
that we have isolated singularities (which are `hyperbolic', in the dynamical systems sense, not to be confused
with elliptic vs. hyperbolic singular points, which will be written without quotes),
no saddle--saddle
connections, and all the
sources or sinks are elliptic singularities or closed orbits which are
Morse--Smale in the usual
sense.  This guarantees the convexity of $\Sigma$.  The actual construction
of the transverse
contact vector field follows from Giroux's argument in \cite{Giroux2} (Proposition
II.2.6), where it is shown that $\Sigma$ is convex if $\Sigma$ is closed
and the characteristic foliation is Morse--Smale.

The goal is to find some $I$--invariant contact structure $\xi'$
(given by a 1--form $\alpha'$)
which induces this characteristic foliation $\mathcal{F}$ on $\Sigma$.  Orient the characteristic foliation so that
the positive elliptic points are the sources and the negative elliptic points are the sinks.  This will naturally identify
which closed orbits are positive (sources) and which closed orbits are negative (sinks).
Let $X$ be a vector field which directs $\mathcal{F}$ and is nonzero away from the singularities of $\mathcal{F}$.
Consider the neighborhood $N(\Sigma)=\Sigma\times I$, where $I$ has coordinate $t$. The
`hyperbolicity' of the singularities implies that if $\xi$ is given by $\alpha=dt+\beta$ (here $\beta$ has no
$dt$--terms, but may be $t$--dependent), then $d\beta$ is nonzero near the singularity on $\Sigma\times\{0\}$.
(This means $X$ has positive {\it divergence} near the singularities.)
Now let $U\subset \Sigma$ be the union of small neighborhoods
of the half-elliptic or half-hyperbolic singularities, elliptic and hyperbolic singularities, the closed orbits,
and neighborhoods of connecting orbits which connect between singularities of the same sign.
Without loss of generality, restrict attention to $U_+$, the components of $U$ with positive singularities.
Let $\beta'$ be a 1--form on $\Sigma$ given by $\beta'=i_X\omega$, where $\omega$ is an area form on $\Sigma$.
The positive divergence ensures that $d\beta'$ is positive near the singular points.  In a neighborhood
$B=S^1\times [-1,1]$ of  a positive closed orbit $S^1\times \{0\}$, with coordinates $(x,y)$, let
$X={\partial \over \partial x}+\phi(x,y){\partial \over\partial y}$, and $\omega=dxdy$.  Then $\beta'=i_X\omega$
satisfies $d\beta'>0$ on $B$, since the Morse--Smale condition implies ${\partial \phi\over \partial y}>0$.
(However, away from the
singularities and closed orbits, we do not know whether $d\beta'$ is positive.) We now take a positive function
$f$ for which $f$ grows rapidly along $X$, ie, $df(X)>>0$, and form $\beta''=f\beta'$.
Since $d\beta''=df\wedge \beta'+fd\beta'$, we obtain  $d\beta''>0$.  Now let $\alpha'=dt+\beta''$.

Now, $\Sigma\backslash U$ consists of annuli $A'=(\R/\Z)\times I$, with coordinates
$(x,y)$ and $\mathcal{F}|_{A'}$ given by $x=\mbox{const.}$,   and $A''=I\times I$, with
coordinates $(x,y)$ and $\mathcal{F}|_{A''}$ also given by $x=\mbox{const.}$  Consider $A'$.
The $I$--invariant contact structure $\xi'$ is defined along $(\R/\Z)\times \{0\}$ by
$f(x,y)dt-dx$ for some positive function $f(x,y)$ satisfying ${\partial f\over \partial y}<0$, and
is defined along $(\R/\Z)\times \{1\}$ by $f(x,y)dt-dx$ for some negative function $f(y)$ satisfying
${\partial f\over \partial y}<0$.  We simply interpolate $f$ between $y=0$ and $y=1$, while keeping
${\partial f\over \partial y}<0$.  $A''$ is similar, but we need to remember that $f$ is already specified
along $\{0,1\}\times I$.

We have therefore constructed an $I$--invariant contact structure $\xi'$  such that
$\xi'|_\Sigma=\mathcal{F}$ and $\xi=\xi'$ on a neighborhood of $A$.  The proof of the proposition is complete
once we have the following lemma.

\begin{lemma}
Let $\Sigma$ be closed or with collared Legendrian boundary.  If $\xi$ and $\xi'$ are
contact structures defined on a neighborhood of $\Sigma$, inducing the same characteristic
foliation $\mathcal{F}$, then there exists a 1--parameter family
of diffeomorphisms $\phi_s$, $s\in[0,1]$, where $\phi_0=id$, $\phi_1^*(\xi')=\xi$, and
$\phi_s$ preserve $\mathcal{F}$.    Moreover, if $\xi$ and $\xi'$ agree on the collared Legendrian
boundary $A$, then $\phi_s$ can be made to have support away from $A$.
\end{lemma}

The proof of this lemma uses Moser's method, and is proven exactly as in Proposition 1.2 of \cite{Giroux2}.
\qed

\subsubsection{Dividing curves}    \label{section:dividing}

A convex surface $\Sigma$ which is closed or compact with Legendrian boundary
has a {\it dividing set} $\Gamma_\Sigma$.
We define a {\it dividing set} $\Gamma_\Sigma$ for $v$ to be the set of points $x$ where
$v(x)\in \xi(x)$.  We will write $\Gamma$ if there is no ambiguity of $\Sigma$.
$\Gamma$ is a union of smooth curves and arcs which are transverse
to the {\it characteristic foliation} $\xi|_\Sigma$.  If $\Sigma$ is closed, there
will only be closed curves $\gamma\subset \Gamma$;
if $\Sigma$ has Legendrian boundary, $\gamma\subset\Sigma$ may be an arc with endpoints on the boundary.
The isotopy type of $\Gamma$ is independent of the choice
of $v$ --- hence we will slightly abuse notation and call $\Gamma$ {\it the dividing set}  of $\Sigma$.
Denote the number of connected components of $\Gamma_\Sigma$ by $\#\Gamma_\Sigma$.
$\Sigma\backslash\Gamma_\Sigma = R_+-R_-$, where $R_+$ is the subsurface where
the orientations of
$v$ (coming from the normal orientation of $\Sigma$) and the normal orientation of $\xi$
coincide, and $R_-$ is the subsurface where they are opposite.

\subsubsection{Giroux's Flexibility Theorem}
The following informal principle highlights the importance of the dividing set:

\medskip

{\bf Key Principle}\qua It is the dividing set $\Gamma_\Sigma$ ({\it not the exact characteristic
foliation}) which encodes the essential
contact topology information in a neighborhood of $\Sigma$.

\medskip

To make this idea more precise, we will now present Giroux's Flexibility Theorem.
If ${\mathcal{F}}$ is a singular foliation on $\Sigma$, then a disjoint union of properly embedded
curves $\Gamma$ is said to {\it divide} ${\mathcal{F}}$  if there exists some $I$--invariant
contact structure $\xi$ on
$\Sigma\times I$ such that ${\mathcal{F}}=\xi|_{\Sigma\times \{0\}}$ and $\Gamma$ is the dividing
set for $\Sigma\times \{0\}$.

\begin{thm}[Giroux \cite{Giroux2}]\label{flexibility}
	Let $\Sigma$ be a closed convex surface or a compact convex surface with
	Legendrian boundary,  with characteristic foliation $\xi|_\Sigma$,
	contact vector field $v$, and dividing set $\Gamma$. If  $\mathcal{F}$ is another
	singular foliation on $\Sigma$ divided by $\Gamma$, then there is an
	isotopy $\phi_s$, $s\in[0,1]$, of $\Sigma$ such that $\phi_0(\Sigma)=\Sigma,$
	$\xi|_{\phi_1(\Sigma)}=\mathcal{F}$, the isotopy is fixed on $\Gamma$,
	and $\phi_s(\Sigma)$ is transverse to $v$ for all $s$.
\end{thm}

An isotopy $\phi_s$, $s\in[0,1]$, for which $\phi_s(\Sigma)\pitchfork v$ for all $s$ is
called {\it admissible}.

\proof
Consider two $I$--invariant contact structures $\xi_0$ and $\xi_1$ on $\Sigma\times I$ which induce
the same dividing set $\Gamma$ on $\Sigma$.  We may assume that $\xi_0=\xi_1$ on
$(N(\Gamma)\cup N(\bdry \Sigma))\times I$.   Here $N(\Gamma)$ and $N(\bdry\Sigma)$ are neighborhoods
of $\Gamma$ and $\bdry\Sigma$ in $\Sigma$.
Consider $\Sigma_0\times I$, where $\Sigma_0$ is a connected component
of $\Sigma\backslash N(\Gamma)$.
Here $\xi_s$, $s=0,1$, will be given by $\alpha_s=dt+\beta_s$, $s=0,1$, where $t$ is the variable in
the $I$--direction, $\beta_s$ is a 1--form on $\Sigma$ which is independent of $t$,
and $d\beta_s>0$.  We interpolate
$\beta_0$ and $\beta_1$ through $\beta_s=(1-s)\beta_0+s\beta_1$, $s\in[0,1]$.
Then $\alpha_s=dt+\beta_s$, $s\in[0,1]$ are all contact and $I$--invariant.  Also note that $\beta_s$ is independent
of $s$ on $N(\bdry \Sigma_0)\times I$.  We use a Moser-type argument to obtain a 1--parameter family
$\{\phi_s\}$ of diffeomorphisms satisfying
\begin{equation}  \label{eqn1}
\phi_s^*(\alpha_s)=f_s\alpha_0,
\end{equation} where $f_s$ is some function.
Differentiating this equation, we obtain:
\begin{equation}    \label{eqn2}
\phi_s^*\left(\mathcal{L}_{{X_s}}\alpha_s +{d\alpha_s\over ds}\right)={df_s\over ds}\alpha_0,
\end{equation}
where $X_s$ is the $s$--dependent vector field $\displaystyle
{d\phi_s\over ds}$, and $\mathcal{L}$ is the Lie derivative.
Substituting Equation \ref{eqn1}  into Equation \ref{eqn2} and removing $\phi_s^*$, we obtain
\begin{equation}
\mathcal{L}_{X_s}\alpha_s=-{d\alpha_s\over ds}+g_s\alpha_s,
\end{equation}
where $g_s$ is some function. We may set $g_s=0$, and solve the pair:
\begin{eqnarray}
i_{X_s}(d\alpha_s)&=&-{d\beta_s\over ds},\\
i_{X_s}(dt+\beta_s)& = &0.
\end{eqnarray}
It is important to note that, since $\beta_s$ is constant along $N(\bdry\Sigma_0)\cup N(\Gamma)$, $X_s=0$ and
$\phi_s$ leaves $(N(\bdry \Sigma_0)\cup N(\Gamma))\times I$ fixed.
By construction, $\phi_s(\Sigma\times\{0\})$ is transverse to
$v$.
\qed

\subsection{Convex surfaces in tight contact manifolds}
From now on let $(M,\xi)$ be a compact, oriented 3--manifold with a tight contact structure
$\xi$.
The following is Giroux's criterion for determining which convex surfaces have neighborhoods which
are tight:

\begin{thm}[Giroux's criterion] \label{thm:Girouxscriterion}
If $\Sigma\not = S^2$ is a convex surface (closed or compact with
Legendrian boundary) in a contact manifold
$(M,\xi)$, then $\Sigma$ has a tight neighborhood if and only if $\Gamma_\Sigma$ has no homotopically
trivial curves.  If $\Sigma=S^2$, $\Sigma$ has a tight neighborhood if and only if
$\#\Gamma_\Sigma=1$.
\end{thm}

We will prove the easy half of the theorem in Section  \ref{section:lerp}.

\medskip

\begin{figure}[ht!]
	{\epsfysize=1.5in\centerline{\epsfbox{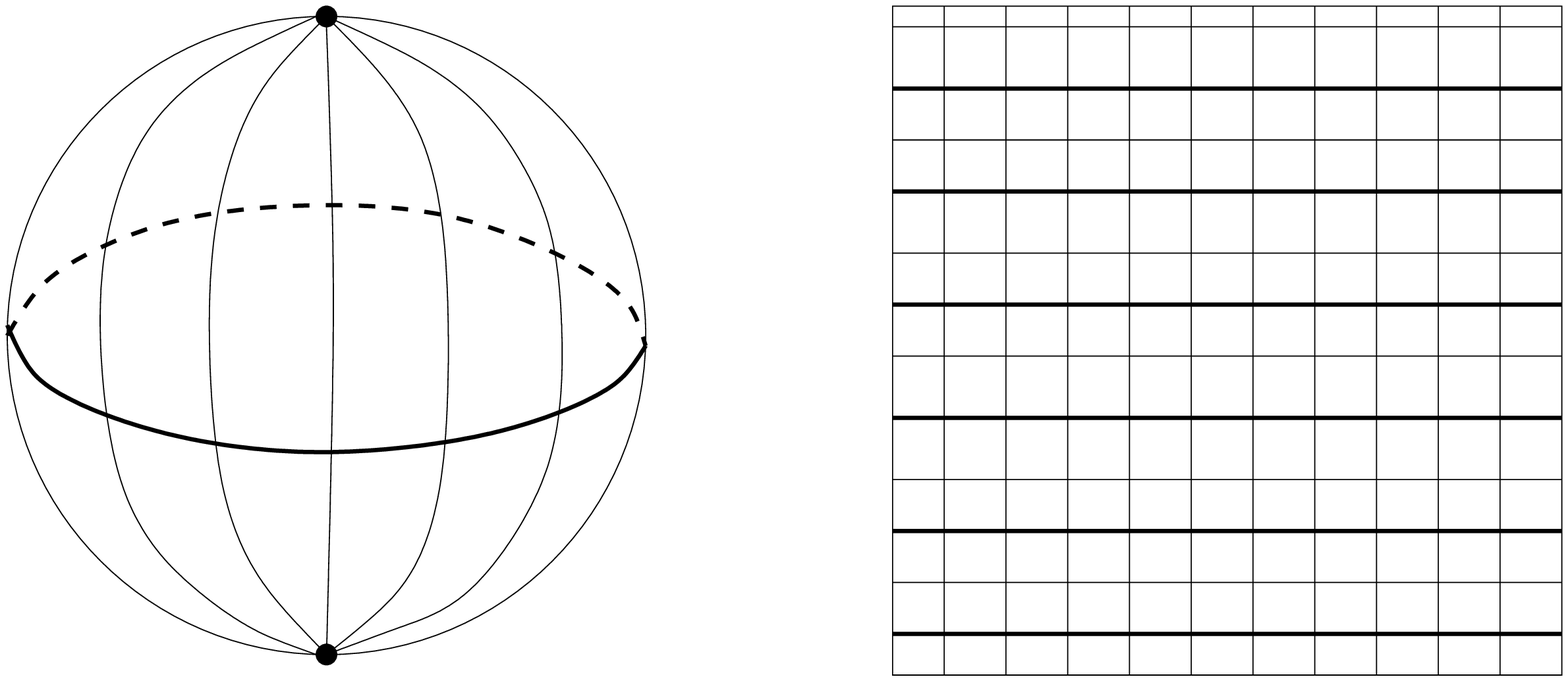}}}
	\caption{Dividing curves for $S^2$ and $T^2$}
	\label{fig18}
\end{figure}

{\bf Examples}\qua  The following are some examples of convex surfaces that can exist inside
tight contact manifolds.
\be
\item $\Sigma=S^2$.  Since $\#\Gamma_\Sigma=1$, there is only one possibility.  See Figure
\ref{fig18}. Note that any time there is more than one dividing curve the contact structure is
overtwisted.   In Figure \ref{fig18}, the thicker lines are the dividing curves and the thin lines represent the characteristic
foliation.

\item $\Sigma=T^2$.  Since there cannot be any homotopically trivial curves, $\Gamma_\Sigma$
consists of an even number $2n>0$ of parallel homotopically essential curves.  Depending on
the identification with $\R^2/\Z^2$ the dividing curves may look like as in Figure \ref{fig18}.
Note that in our planar representation of $T^2$ the sides are identified and the top and bottom are identified.
\ee

\subsubsection{Convex tori in standard form}    \label{section:tori}
One of the main ingredients in our study is the convex torus $\Sigma\subset M$
in {\it standard form}.  Assume $\Sigma$ is a convex torus in a tight contact manifold $M$.
Then, after
some identification of $\Sigma$ with $\R^2/\Z^2$, we may assume $\Gamma_\Sigma$
consists of $2n$ parallel homotopically essential curves of slope $0$.
The {\it torus division number} is given by $n={1\over 2}(\#\Gamma_\Sigma)$. Using Giroux's
Flexibility Theorem, we can
deform $\Sigma$
inside a neighborhood of $\Sigma\subset M$ into a torus which we still call $\Sigma$ and has the
same dividing set as the old $\Sigma$.
The characteristic foliation on this new $\Sigma=\R^2/\Z^2$ with coordinates $(x,y)$ is
given by
$y = rx + b$, where $r \ne 0$ is fixed, and $b$ varies in a family,
with tangencies $y={k\over 2n}$, $k=1,...,2n$. ($r = \infty$ will also be
allowed, in which case
we have the family $x = b$.)  We say such a $\Sigma$ is a {\it convex torus in standard form} (or simply
{\it in standard form}).  The horizontal
Legendrian curves $y={k\over 2n}$ are isolated and rather inflexible from
the point of view of
$\Sigma$ (as well as nearby convex tori), and will be called {\it Legendrian
divides}. The Legendrian curves that are in a family are much more
flexible, and will be called {\it Legendrian rulings}.  In particular, a consequence of Giroux's Flexibility Theorem
is the following:

\begin{cor}[Flexibility of Legendrian rulings]\label{rulings} 
Let $(\Sigma, \xi_{\Sigma})$ be a
torus in the above form, with coordinates $(x,y)\in\R^2/\Z^2$, Legendrian rulings $y=rx+b$
(or $x = b$), and Legendrian divides $y={k\over 2n}$. Then, via a $C^0$--small perturbation
near the Legendrian divides, we can modify the slopes of the rulings from $r \ne 0$ to any
other number $r' \ne 0$ ($r =\infty$ included).
\end{cor}

We will also say that a convex annulus $\Sigma=S^1\times I$ is in {\it
standard form} if,
after a diffeomorphism,
$S^1\times \{pt\}$ are Legendrian (ie, they are the Legendrian rulings),
with tangencies
$z={k\over 2n}$ (Legendrian divides), where $S^1=\R/\Z$ has coordinate $z$.

\subsection{Convex decompositions}

Let $(M,\xi)$ be a compact, oriented, tight contact 3--manifold with nonempty convex boundary
$\bdry M$.  Suppose $\Sigma$ is a properly embedded oriented surface with $\bdry \Sigma\subset
\bdry M$.  In this section we describe how to perturb $\Sigma$ into a convex surface with
Legendrian boundary (after possible modification of the characteristic foliation on $\bdry M$), and
perform a {\it convex decomposition}.

\subsubsection{Legendrian realization principle}  \label{section:lerp}

In this section we present the {\it Legendrian realization principle} --- a criterion for determining
whether a given curve or a collection of curves and arcs can be made Legendrian after a
perturbation of a convex surface $\Sigma$.  The result is surprisingly strong --- we can realize
almost any curve as a Legendrian one.   Our formulation of Legendrian realization is a
generalization of Kanda's \cite{K98}.
Call a union of disjoint properly embedded closed curves and arcs $C$ on a convex surface
$\Sigma$ with Legendrian boundary {\it nonisolating}
if (1) $C$ is transverse to $\Gamma_\Sigma$, and every arc in $C$ begins and ends on
$\Gamma_\Sigma$, and (2) every component of $\Sigma\backslash (\Gamma_\Sigma \cup
C)$ has a boundary component
which intersects $\Gamma_\Sigma$.  Here, $C\pitchfork \Gamma_\Sigma$, strictly speaking, makes
sense only after we have fixed a
contact vector field $v$.  For the Legendrian realization principle and its corollary, the contact
structure $\xi$ does not need to be tight.

\begin{thm}[Legendrian realization] \label{lerp}
Consider $C$, a nonisolating collection of disjoint properly embedded closed curves and arcs,
on a convex surface $\Sigma$
with Legendrian boundary.
Then there exists an
admissible isotopy $\phi_s$, $s\in[0,1]$ so that
\be
\item $\phi_0=id$,
\item $\phi_s(\Sigma)$ are all convex,
\item $\phi_1(\Gamma_\Sigma)=\Gamma_{\phi_1(\Sigma)}$,
\item $\phi_1(C)$ is Legendrian.
\ee
\end{thm}
Therefore, in particular, a nonisolating collection $C$ can be realized by a Legendrian collection
$C'$ with the same number of geometric intersections.    A corollary of this theorem, observed by
Kanda, is the following:

\begin{cor}[Kanda] A closed curve $C$ on $\Sigma$ which is transverse to $\Gamma_\Sigma$
can be realized as a Legendrian curve (in the
sense of Theorem \ref{lerp}), if $C\cap \Gamma_\Sigma\not=\emptyset$.
\end{cor}

Observe that  if $C$ is a Legendrian curve on a convex surface $\Sigma$, then its twisting
number $t(C,Fr_\Sigma)={1\over 2}\#(C\cap \Gamma_\Sigma)$, where $\#(C\cap \Gamma_\Sigma)$
is the geometric intersection number (signs ignored).

\proof
By Giroux's Flexibility Theorem, it suffices to find a
characteristic foliation ${\mathcal{F}}$ on $\Sigma$ with (an isotopic copy of) $C$ which is represented
by Legendrian curves and arcs.  We remark here that these Legendrian curves and arcs constructed
will always pass through singular points of ${\mathcal{F}}$.
Consider a component $\Sigma_0$ of  $\Sigma\backslash (\Gamma_\Sigma \cup
C)$ --- let us assume $\Sigma_0\subset R_+$, so all the elliptic singular points are
sources.
Denote $\bdry \Sigma_0 = \gamma^-- \gamma^+$, where $\gamma^-$ consists of
closed curves $\gamma$ which intersect $\Gamma_\Sigma$, and $\gamma^+$ consists of
closed curves $\gamma\subset C$.
This means that for $\gamma\subset \gamma^-$,
either $\gamma\subset \Gamma_\Sigma$ or $\gamma=
\delta_1\cup \delta_2\cup\cdots \cup \delta_{2k}$, where
$\delta_{2i-1}$, $i=1,\cdots,k$, are subarcs of $C$, $\delta_{2i}$, $i=1,\cdots, k$,
are subarcs of $\Gamma_\Sigma$, and the endpoint of $\delta_j$ is the initial point of $\delta_{j+1}$.
Since $C$ is nonisolating, $\gamma^-$ is nonempty.
What the $\gamma^-$ provide are `escape routes' for the flows whose
sources are $\gamma^+$ or the singular set of $\Sigma_0$  ---
in other words, the flow would be exiting  along
$\Gamma_\Sigma$.

\begin{figure}[ht!]
\centerline{\small
\psfrag {d1}{$\delta_1$}
	\psfrag {d2}{$\delta_2$}
	\psfrag {d3}{$\delta_3$}
	\psfrag {d4}{$\delta_4$}
	\psfrag {dn}{$\delta_n$}
	\psfrag {dn1}{$\delta_{2n-1}$}
\epsfysize=1.5in\epsfbox{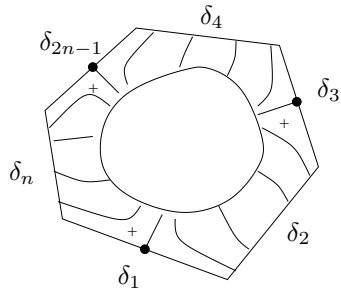}} 
\caption{Characteristic foliation on $\gamma\times I$}
\label{fig2}
\end{figure}

Construct ${\mathcal{F}}$ so that (1) the subarcs of $\gamma^-$ coming from $C$
are now Legendrian, with a single positive half-hyperbolic point in the interior of the arc,
(2) the curves of
$\bdry \Sigma_0$ contained in $C$ are Legendrian curves, with one positive half-elliptic point
and one positive half-hyperbolic point.   If $\gamma\subset \gamma^-$ intersects $C$,
then we give a neighborhood $\gamma\times I$ a characteristic foliation as in Figure \ref{fig2}.
After filling in this collar, we may assume that ${\mathcal{F}}$ is transverse to and flows out
of $\gamma^-$.   If $\gamma^+$ is empty, then we introduce a positive elliptic singular
point on the interior of $\Sigma_0$, and let $\gamma^+$ be a small closed loop around
the singular point, transverse to the flow.  At any rate, we may assume the flow enters through
$\gamma^+$ and exits through $\gamma^-$ --- by filling in appropriate positive hyperbolic
points we may extend ${\mathcal{F}}$ to all of $\Sigma_0$.
\qed

\medskip

Actually, Kanda observes the following slightly stronger statement.  The proof is identical --- instead
of single Legendrian curves, we insert a collar neighborhood.

\begin{cor}[Kanda] Let $C\pitchfork \Gamma_\Sigma$ be a closed curve on $\Sigma$ which satisfies
$|C\cap \Gamma_\Sigma|\geq 2$.  Then $C$ can be realized as a Legendrian curve, and, moreover,
$C$ can be made to have a standard annular collar neighborhood $A\subset\Sigma$ consisting of
a 1--parameter family of
Legendrian ruling curves which are translates of $C$.
\end{cor}

We will now give a proof of one-half of Theorem \ref{thm:Girouxscriterion}, as a corollary of
the Legendrian realization principle.  The converse is more involved, and will be omitted (it will not
be used in this paper).

\medskip

{\bf Proof of Giroux's Criterion}\qua
Assume $\Gamma_\Sigma$ has a homotopically trivial curve $\gamma$ which bounds a disk
$D$.  Then there exists a curve
$\gamma'\subset \Sigma\backslash D$ parallel to $\gamma$, such that $\gamma'\cap \Gamma_\Sigma=\emptyset$.
Provided
$\Gamma_\Sigma$ does not consist solely of the homotopically trivial curve $\gamma$,
$\gamma'$ is nonisolating, and
we may use Legendrian realization and assume, after modifying $\Sigma$ inside an $I$--invariant
neighborhood, that $\gamma'$ is Legendrian, and $t(\gamma')=0$ with
respect to $\Sigma$.  This implies that $\gamma'$ bounds an overtwisted disk.
The case $\#\Gamma_\Sigma=1$ requires a bit more work and one operation which
is introduced later.  We may assume $\Sigma$ is not a disk, since the boundary
Legendrian curve would then bound an overtwisted disk.
Take a closed curve $\delta\subset \Sigma$ which is homotopically essential,
has no intersection with $\#\Gamma_\Sigma$, and does not separate $\Sigma$ (note
that $\delta$ may be a boundary Legendrian curve).    Use Legendrian realization
to realize $\delta$ as a Legendrian curve with $t(\delta)=0$.  At this point, we will need to
apply the `folding' method for increasing the dividing curves described in Section \ref{section:increase}.
Each fold will introduce a pair of dividing curves parallel to $\delta$.  Now $\gamma'$ is
Legendrian-realizable.
\qed

\medskip

\subsubsection{Cutting and rounding}
Suppose $\Sigma\subset M$ is a properly embedded oriented surface with $\bdry \Sigma\subset
\bdry M$, where $\bdry M$ is convex.    Make $\bdry \Sigma\pitchfork \Gamma_{\bdry M}$, and modify
$\bdry \Sigma$ (by adding extraneous intersections) if necessary,
so that $|\bdry \Sigma\cap \Gamma_{\bdry M}|>0$.
Using the Legendrian realization principle, we may arrange $C$ to be Legendrian on
$\bdry M$, with a standard annular collar, after perturbation.

$C$ has a neighborhood $N(C)$ which is locally isomorphic to the neighborhood $\{x^2
+ y^2  \le \varepsilon\}$ of $M = {\bf R}^2 \times ({\bf R}/{\bf Z})$ with
coordinates $(x, y, z)$ and contact 1--form $\alpha = \sin(2\pi n z) dx +
\cos(2\pi n z) dy$, where $n={1\over 2}|C\cap \Gamma_{\bdry M}| \in {\bf Z}^+$.     Here
$C=\{x=y=0\}$ and $\bdry M\cap N(C)=\{x=0\}$.   Also let $\Sigma\cap N(C)=\{y=0\}$ and perturb
the rest (fixing $\Sigma\cap N(C)$) so $\Sigma$ is convex with Legendrian boundary.

\begin{lemma}  \label{lemma:arrange}
It is possible to arrange the transverse contact vector field $X$ for $\bdry M$ to be ${\bdry \over \bdry x}$
and the transverse contact vector field $Y$ for $\Sigma$ to be ${\bdry \over \bdry y}$.
\end{lemma}

\proof
Follows from Proposition \ref{prop:makeconvex}.
\qed

\medskip
Now cut $M$ along $\Sigma$ to obtain $M\backslash \Sigma$ (which we really mean to be $M\backslash
int(\Sigma\times I)$).  Then round the edges using the following
edge-rounding lemma:

\begin{lemma}[Edge-rounding] \label{edge}  Let $\Sigma_1$ and $\Sigma_2$
be convex surfaces with collared Legendrian boundary which intersect transversely inside the
ambient contact manifold along a common boundary Legendrian curve.
Assume the neighborhood of
the common boundary Legendrian is locally isomorphic to the neighborhood $N_\varepsilon=\{x^2
+ y^2  \le \varepsilon\}$ of $M = {\bf R}^2 \times ({\bf R}/{\bf Z})$ with
coordinates $(x, y, z)$ and contact 1--form $\alpha = \sin(2\pi n z) dx +
\cos(2\pi n z) dy$, for some $n \in {\bf Z}^+$, and that $\Sigma_1\cap N_\varepsilon = \{x = 0, 0 \le
y \le \varepsilon\}$
and $\Sigma_2\cap N_\varepsilon = \{y = 0, 0 \le x \le \varepsilon\}$.  If we join $\Sigma_1$
and $\Sigma_2$ along
$x = y = 0$ and round the common edge (take $((\Sigma_1\cup \Sigma_2)\backslash N_\delta)\cup
(\{(x-\delta)^2+(y-\delta)^2=\delta^2\}\cap N_\delta)$, where $\delta<\varepsilon$),
the resulting surface is convex, and the
dividing curve
$z = {k\over 2n}$ on $\Sigma_1$ will connect to the dividing curve $z = {k\over 2n} - {1\over
4n}$ on $\Sigma_2$, where $k = 0, \cdots, 2n - 1$.  Here we assume that the orientations of
$\Sigma_1$ and $\Sigma_2$ are compatible and induce the same orientation after rounding.
\end{lemma}

\begin{figure}[ht!]
\centerline{\small
\psfrag {S1}{$\Sigma_1$}
\psfrag {S2}{$\Sigma_2$}
\epsfysize=2in\epsfbox{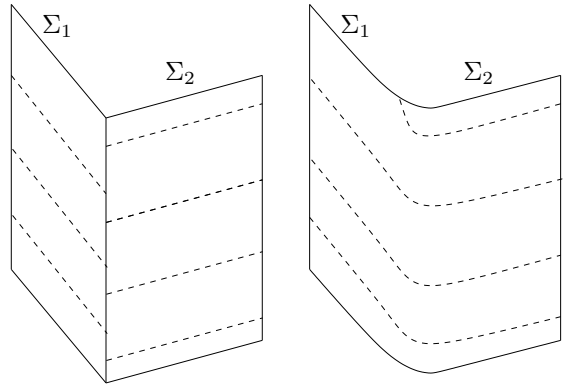}} 
\caption{Edge rounding:  Dotted lines are dividing curves.}
\label{rounding}
\end{figure}

Refer to Figure \ref{rounding}.

\proof This follows from Lemma \ref{lemma:arrange}, and taking the transverse vector field for
$\Sigma_1$ to be ${\partial \over \partial x}$ and taking the transverse vector field
for $\Sigma_2$ to be ${\partial \over \partial y}$.  The transverse vector field for
$\{(x-\delta)^2+(y-\delta)^2=\delta^2\}\cap N_\delta$ is the inward-pointing radial vector
$-{\bdry \over \bdry r}$ for
the circle    $\{(x-\delta)^2+(y-\delta)^2=\delta^2\}$.
\qed

\subsection{Bypasses}  \label{by}

\begin{figure}[ht!]
\centerline{\small
\psfrag {g1}{$\gamma_1$}
\psfrag {g2}{$\gamma_2$}
\epsfysize=2in\epsfbox{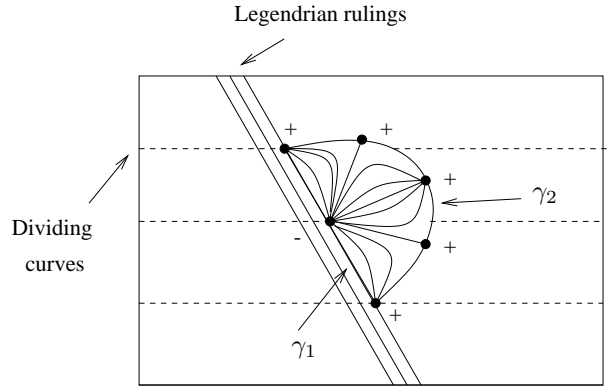}} 
\caption{A bypass}
\label{bypass}
\end{figure}

Let $\Sigma\subset M$ be convex surface (closed or compact with Legendrian boundary).
A {\it bypass} for $\Sigma$ is an oriented embedded half-disk $D$  with Legendrian boundary,
satisfying the following:
\be
\item $\bdry D$ is the union of two arcs $\gamma_1$, $\gamma_2$ which intersect at their
endpoints.
\item $D$ intersects $\Sigma$ transversely along $\gamma_1$.
\item $D$ (or $D$ with opposite orientation) has the
following tangencies along
$\bdry D$:
\be
\item positive elliptic tangencies at
the endpoints of $\gamma_1$ (= endpoints of $\gamma_2$),
\item one negative elliptic tangency on the interior of $\gamma_1$, and
\item only positive tangencies along $\gamma_2$, alternating between elliptic and hyperbolic.
\ee
\item $\gamma_1$ intersects $\Gamma_\Sigma$ exactly at three points, and these three points are
the elliptic points of $\gamma_1$.
\ee
Refer to Figure \ref{bypass} for an illustration.
We will often also call the arc $\gamma_2$ a {\it bypass for} $\Sigma$ or a {\it bypass for} $\gamma_1$.
We define the {\it sign} of a bypass to be the sign of the half-elliptic point at the center of
the half-disk.

\subsubsection{Bypass attachment lemma}

\begin{figure}[ht!]
{\epsfysize=1.5in\centerline{\epsfbox{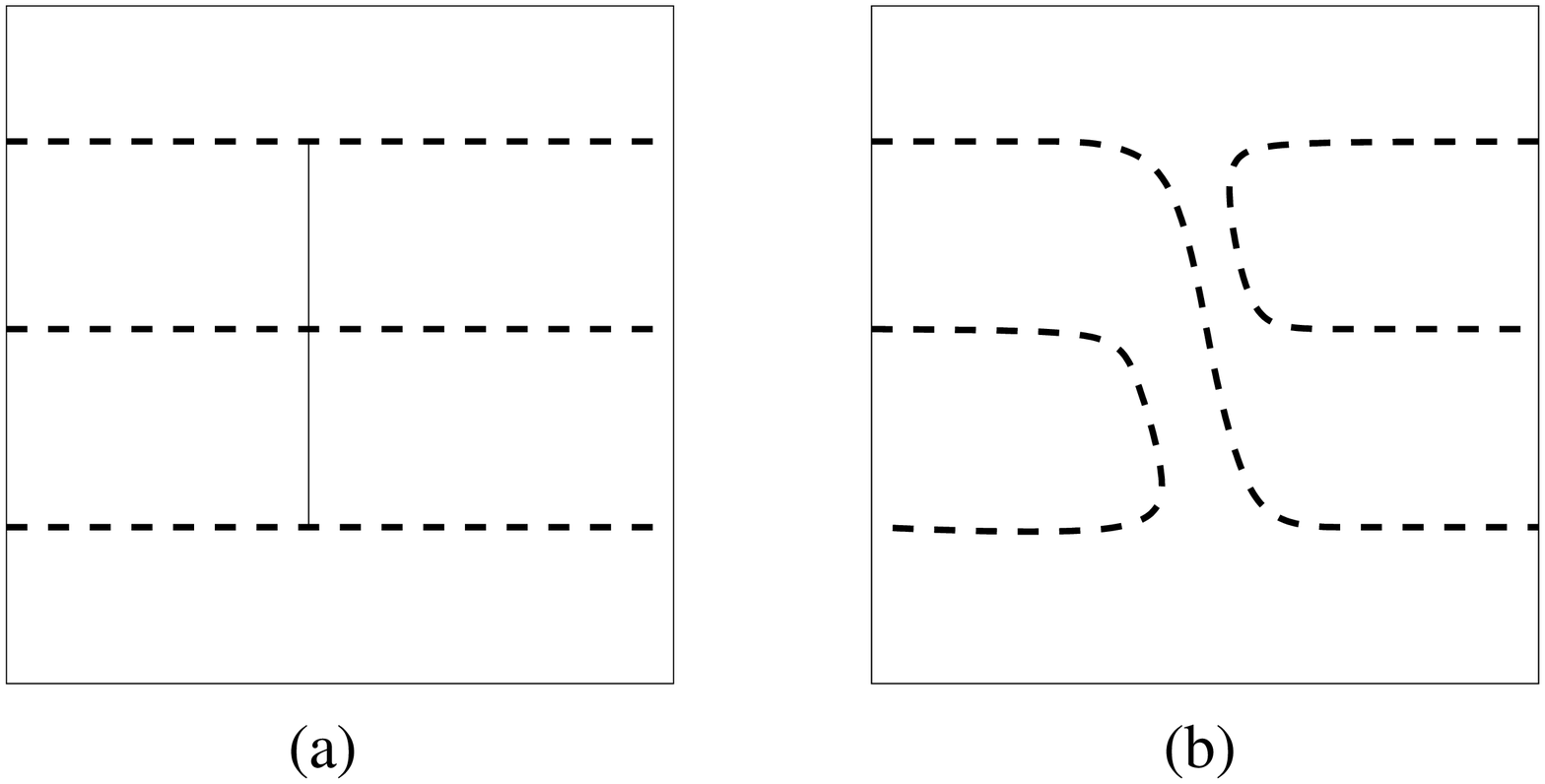}}}
\caption{Bypass attachment:\qua (a) Dividing curves on $\Sigma_0$.
(b) Dividing curves on $\Sigma_1$.  The dividing curves are dotted lines, and the Legendrian arc of
attachment $\gamma_1$ is a solid line.  We are only looking at the portion of $\Sigma_i$ where the
attachment is taking place.}
\label{bypassmove}
\end{figure}

\begin{lemma}[Bypass Attachment] Assume $D$ is a bypass for a convex $\Sigma$.  Then there exists a
neighborhood of $\Sigma\cup D\subset M$ diffeomorphic to $\Sigma\times [0,1]$, such that
$\Sigma_i=\Sigma\times \{i\}$, $i=0,1$, are convex, $\Sigma\times[0,\varepsilon]$ is $I$--invariant,
$\Sigma=\Sigma\times \{\varepsilon\}$, and $\Gamma_{\Sigma_1}$
is obtained from $\Gamma_{\Sigma_0}$ by performing the Bypass Attachment operation depicted
in Figure \ref{bypassmove} in a neighborhood of the attaching Legendrian arc $\gamma_1$.
\end{lemma}

\proof Extend $\gamma_1$ to a closed Legendrian curve $\gamma$ on $\Sigma$ using the Legendrian
Realization Principle.   We may also assume that $\gamma$ has an annular neighborhood of $\Sigma$ which
is in standard form, and that $D$ is a convex half-disk transverse to $\Sigma$.  Take an $I$--invariant
one-sided neighborhood
$\Sigma\times [0,\varepsilon]$ of $\Sigma$, where
$\Sigma=\Sigma\times\{\varepsilon\}$.  Now, $A'=\gamma\times [0,\varepsilon]\subset
\Sigma\times [0,\varepsilon]$ is an annulus in standard form transverse to $\Sigma\times \{0\}$.  Form
$A=A'\cup D$.  $A$ is convex, and we can take an $I$--invariant neighborhood $N(A)$ of $A$.  If $\bdry A$ was
smooth, then we take $(\Sigma\times\{0\})\cup N(A)$, and smooth out the four edges
using the Edge-Rounding Lemma.

To smooth out $\bdry A$, we use the Pivot Lemma, first observed by Fraser \cite{F}.  The proof is similar to
the Flexibility Theorem.

\begin{lemma}[Pivot] Let $S$ be an embedded disk in a contact manifold $(M,\xi)$
with a characteristic foliation $\xi|_{S}$ which consists only of one positive elliptic singularity $p$ and unstable orbits
from $p$ which exit transversely from $\bdry S$.
If $\delta_1, \delta_2$ are two unstable orbits meeting at  $p$, and $\delta_i\cap \bdry S=p_i$, then,
after a $C^\infty$--small
perturbation of $S$ fixing $\bdry S$, we obtain $S'$ whose characteristic foliation
has exactly one positive elliptic singularity $p'$ and
unstable orbits from $p'$ exiting transversely from $\bdry S$, and for which the orbits passing through $p_1$, $p_2$
meet tangentially at $p'$.
\end{lemma}

\begin{figure}[ht!]
{\epsfysize=1in\centerline{\epsfbox{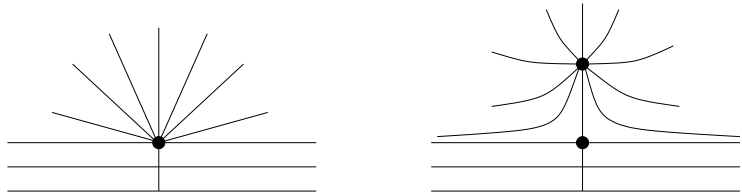}}}
\caption{Replacing a half-elliptic point by a half-hyperbolic point and a full elliptic point}
\label{fig21}
\end{figure}

\begin{figure}[ht!]
{\epsfysize=0.8in\centerline{\epsfbox{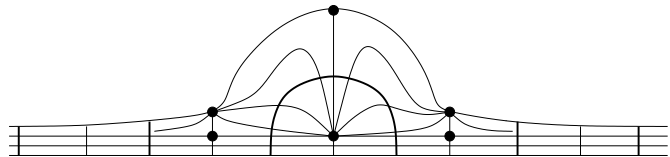}}}
\caption{Smoothing the corners of $A$ using the Pivot Lemma}
\label{fig22}
\end{figure}

\begin{figure}[ht!]
{\epsfysize=1.5in\centerline{\epsfbox{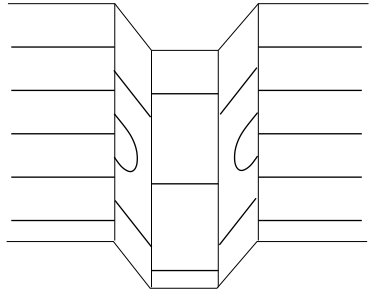}}}
\caption{Rounding the edges will give the desired dividing set}
\label{fig23}
\end{figure}

Now consider the half-elliptic singular points $q_1,q_2$ on $D$ which are also the endpoints of $\gamma_1$.
Modify $D$ near $q_i$ to replace $q_i$ by a pair $q^e_i$, $q^h_i$, where $q^e_i$ is a (full) elliptic point and
$q^h_i$ is a half-hyperbolic point as pictured in Figure \ref{fig21}.
Use the Pivot Lemma to smooth the corners of $A$ as in Figure \ref{fig22}.
$A$ is now convex with Legendrian boundary.  The dividing curves on $A$ are the thicker straight lines in
Figure \ref{fig22}.    Finally, we round the edges (see Figure \ref{fig23}) using the Edge-Rounding Lemma.
\qed

\begin{figure}[ht!]
{\epsfysize=1in\centerline{\epsfbox{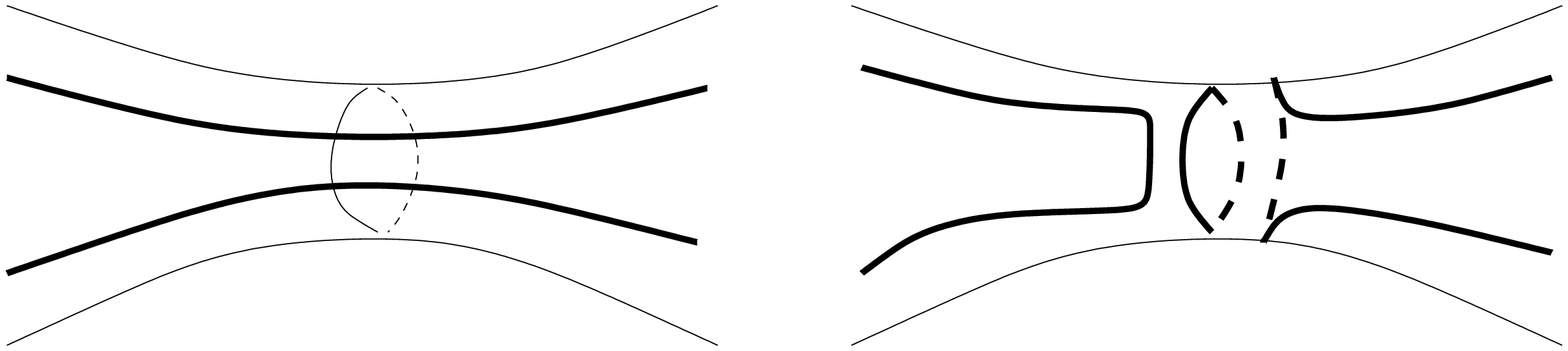}}}
\caption{Edge-Rounding for a singular bypass}
\label{fig24}
\end{figure}

\medskip
We can also a define a {\it singular bypass} to be an immersion $D\rightarrow M$ which satisfies all the
conditions of a bypass except one:  $D$ is an embedding away from $\gamma_1\cap \gamma_2$, and these
two points get mapped to one point on $\Sigma$.  In this case, the Bypass Attachment Lemma would be as
in Figure \ref{fig24}.

\subsubsection{Tori}
Let $\Sigma\subset M$ be a convex torus in standard form, identified with
$\R^2/\Z^2$.  With this identification we will assume that the Legendrian
divides and rulings are already linear,
and will  refer to {\it slopes} of
Legendrian divides and Legendrian rulings. The slope of the Legendrian
divides of $\Sigma$
will be called the {\it boundary slope} $s$ of $\Sigma$, and the slope of
the Legendrian
rulings will be the {\it ruling slope} $r$.
Now assume, after acting via $SL(2,\Z)$, that $\Sigma$ has $s=0$ and $r\not=0$
rational.
Note that we can normalize the Legendrian
rulings via an element $\left (\begin{array}{cc}
1 & m\\ 0 & 1 \end{array}\right ) \in SL(2,\Z)$, $m\in\Z$, so that
$-\infty<r\leq -1$.

In our later analysis on $T^2\times I$ we will find an abundance of bypasses, and
use them to stratify a given $T^2\times I$ with a tight contact structure
and convex boundary into thinner, more basic slices of $T^2\times I$.

\begin{lemma}[Layering] \label{layer} Assume a bypass $D$ is attached to
$\Sigma=T^2$ with slope $s(T^2)=0$, along a Legendrian ruling curve of slope $r$
with $-\infty < r \le -1$.
Then there exists a neighborhood $T^2\times I$ of
$\Sigma\cup D
\subset M$, with $\bdry (T^2\times
I)=T_1- T_0$,
such that $\Gamma_{T_0}=\Gamma_\Sigma$, and $\Gamma_{T_1}$ will be as follows, depending
on whether $\#\Gamma_{T_0}=2$ or $\#\Gamma_{T_0}>2$:
\be
\item If $\#\Gamma_{T_0}>2$, then $s_1=s_0=0$, but $\#\Gamma_{T_1}=\#\Gamma_{T_0}-2$.

\item If $\#\Gamma_{T_0}=2$, then $s_1=-1$, and $\#\Gamma_{T_1}=2$.
\ee
Here $s_i$ is the boundary slope of $T_i$.
\end{lemma}

\begin{figure}[ht!]
{\epsfysize=3.5in\centerline{\epsfbox{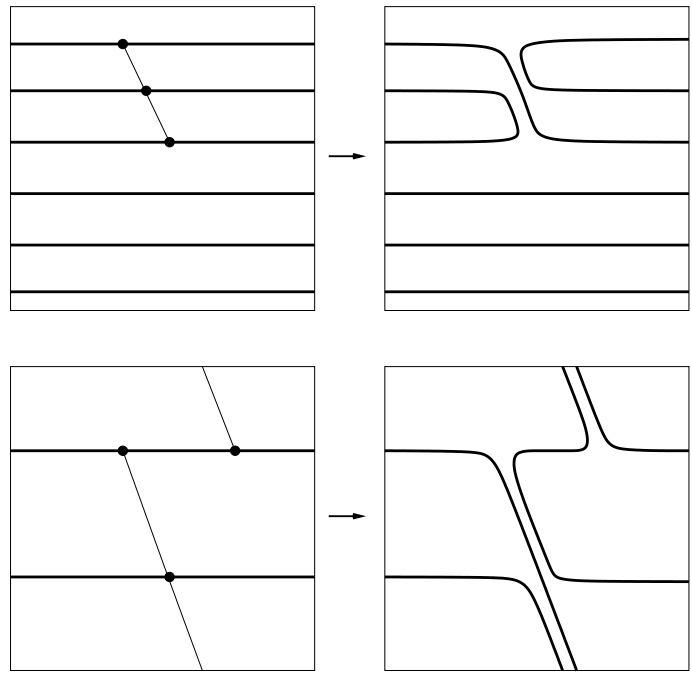}}}
\caption{Bypass attachments along $T^2$}
\label{fig25}
\end{figure}

\proof Follows from the Bypass Attachment Lemma.  Refer to Figure \ref{fig25} for the two possibilities.
\qed

\medskip
Notice that in the case $\#\Gamma_{T^2}=2$, a bypass attachment effectively performs a {\it positive Dehn twist}.

\subsubsection{Tessellation picture}
\begin{figure}[ht!]
	{\epsfysize=2.5in\centerline{\epsfbox{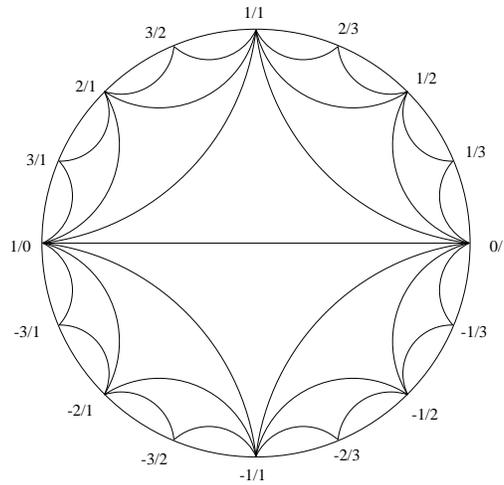}}}
	\caption{The standard tessellation of the hyperbolic unit disk}
	\label{tessellation}
\end{figure}

In this section we interpret the Bypass Attachment Lemma in terms of the standard (Farey)
tessellation of the hyperbolic unit disk
$\H^2=
\{(x,y)|x^2+y^2\leq 1\}$.
Recall we start by labeling $(1,0)$ as $0={0\over 1}$, and $(-1,0)$ as $\infty={1\over 0}$.
We inductively label points on $S^1=\bdry \H^2$ as follows (for $y>0$):  Suppose we have already
labeled $\infty \geq {p\over q}\geq 0$ ($p,q$ relatively prime) and $\infty \geq {p'\over q'}\geq 0$ ($p',q'$ relatively prime)
such that $(p,q)$, $(p',q')$ form a $\Z$--basis
of $\Z^2$.  Then, halfway between ${p\over q}$ and ${p'\over q'}$ along $S^1$ on the shorter arc (one for which
$y>0$ always), we label ${p+p'\over q+q'}$.
We then connect two points ${p\over q}$ and ${p'\over q'}$ on the boundary, if the corresponding shortest integral
vectors form an integral basis of $\Z^2$.  See Figure \ref{tessellation}.

By transforming the situation in Lemma \ref{layer} via $SL(2,\Z)$, we obtain the following rephrasing in more
invariant language.

\begin{lemma} \label{slopes}
Let $\Sigma=T^2$ be a convex surface with $\#\Gamma_{T^2}=2$ and  slope $s=s(T^2)$.
If a bypass $D$ is attached to $\Sigma$ along
a Legendrian ruling curve of slope $r\not= s$, then the resulting  convex surface $\Sigma'$ will have
$\#\Gamma_{T^2}=2$ and slope $s'$ which is obtained as follows:  Take the arc $[r,s]\subset \bdry \H^2$ obtained
by starting
from $r$ and moving counterclockwise until we hit $s$.  On this arc, let $s'$ be the point which is closest to
$r$ and has an edge from $s'$ to $s$.
\end{lemma}

\subsubsection{Abundance of bypasses}

In this section we will demonstrate that bypasses are usually quite abundant.  Suppose
$M$ is a 3--manifold with convex boundary, and we cut $M$ along a convex surface with
Legendrian boundary.  The following are ways in which bypasses can occur.

\begin{lemma}    \label{bccdisk}
Let $\Sigma=D^2$ be a convex surface with Legendrian boundary inside a
tight contact manifold, and $t(\bdry \Sigma, Fr_\Sigma)=-n<0$.  Then every component
of $\Gamma_\Sigma$ is an arc which begins and ends on $\bdry \Sigma$. There exists a bypass
along $\bdry\Sigma$ if $t(\bdry\Sigma)<-1$.
\end{lemma}

\proof If there is a closed dividing curve $\gamma$, then $\gamma$ must bound a disk, contradicting
Giroux's criterion.  Therefore, every dividing curve must be an arc which begins and ends on the
boundary.  Now, if we have arranged $\Sigma$ to have a collared Legendrian boundary
and all half-elliptic points, then the endpoints of the dividing curves will lie between the
half-elliptic points. There will be $2|t(\bdry\Sigma)|$ endpoints for dividing curves,
and hence $|t(\bdry\Sigma)|$ curves.  Now assume $t<-1$.  Then there will exist an `outermost'
dividing curve $\gamma$ ---
one that begins and ends on consecutive endpoints and cuts off a half-disk $D_1$ which
does not contain any other dividing curve.  Take an arc $\delta\subset \Sigma\backslash D_1$ which
is parallel to $\gamma$ and does not intersect $\Gamma$.  Using the Legendrian realization principle
(and the fact that $t<-1$, so that there are at least two half-elliptic points on $\Sigma\backslash D_1$),
we can take $\delta$ to be a Legendrian arc after possible modification;
$\delta$ cuts off a half-disk $D_2\subset \Sigma$
(containing $D_1$) which is a bypass.\qed

\begin{figure}[ht!]
{\epsfysize=3in\centerline{\epsfbox{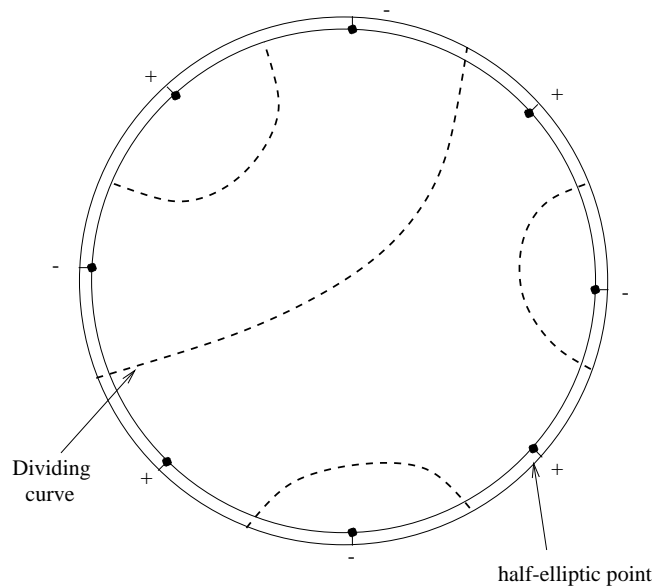}}}
\caption{Standardized convex disk with Legendrian boundary}
\label{disk}
\end{figure}

\medskip
Figure \ref{disk} illustrates a possible dividing set on $\Sigma=D^2$.

\begin{prop} [Imbalance Principle] \label{imba}  Let $\Sigma = S^1
\times [0,1]$ be convex with Legendrian boundary inside a tight contact manifold.
If $t(S^1 \times \{0\}) < t(S^1 \times\{1\}) \le 0$,
then there exists a bypass along $S^1
\times \{0\}$.
\end{prop}

\proof
Let $t_i=t(S^1\times\{i\})$, $i=0,1$.
There exist $2|t_0|$ endpoints of dividing curves on $S^1\times\{0\}$ and
$2|t_1|$ endpoints on $S^1\times \{1\}$.  If $t_0<t_1$, then there exist
two endpoints on $S^1\times\{0\}$
which are connected by the same dividing arc $\gamma$.  This $\gamma$ must bound a half-disk $D_1$,
and hence
there is a Legendrian arc $\delta$ which bounds a bypass half-disk $D_2\supset D_1$.
\qed

\medskip
Let $\Sigma$ be a convex surface with (nonempty) Legendrian boundary, and $\gamma$ be a dividing
curve which cuts off a half-disk $D\subset \Sigma$ which has no other intersections with $\Gamma_\Sigma$.
Such a dividing curve will be called a {\it boundary-parallel} dividing curve.
We can generalize the above discussion and state the following (the proof is immediate):
\begin{prop}     \label{bcc}
Let $\Sigma$ be a convex surface with Legendrian boundary, and $\gamma$ a boundary-parallel
dividing curve.  If $\Sigma$ is not a disk with $t(\bdry \Sigma)=-1$, then there exists a bypass
half-disk which contains the half-disk cut off by $\gamma$.
\end{prop}

\section{Layering of  $T^2 \times I$, $S^1\times D^2$, and $L(p,q)$}

\subsection{Basic building blocks}

In this section we will review the basic building blocks of tight contact manifolds.

\subsubsection{3--ball} Recall the following fundamental theorem of Eliashberg \cite{E}:

\begin{thm}   \label{thm:3-ball} Assume there exists a contact structure $\xi$ on a neighborhood of
$\bdry B^3$ which makes $\bdry B^3$ convex with $\#\Gamma_{\bdry B^3}=1$.
Then there exists a unique extension of $\xi$ to a tight contact structure on $B^3$,
up to an isotopy which fixes the boundary.
\end{thm}

The basic building blocks of tight contact manifolds are $B^3$, equipped with a unique tight
contact structure  if we prescribe the  boundary.

\subsubsection{Flexibility of characteristic foliation on boundary}

Let $M$ have nonempty boundary, and $\mathcal{F}$ be a characteristic foliation which is {\it
adapted} to a dividing set $\Gamma_{\bdry M}$.  Denote by $\mbox{Tight}(M,\mathcal{F})$ the
set of smooth contact 2--plane fields $\xi$ on $M$ which induce a characteristic foliation $\mathcal{F}$
on $\bdry M$.    Then $\pi_0(\mbox{Tight}(M,\mathcal{F}))$ consists of the isotopy classes of
tight contact structures on $M$ with fixed boundary characteristic foliation $\mathcal{F}$.
The Flexibility Theorem allows us to prove the following:

\prop \label{prop:flex}
Let $M$ be a compact, oriented 3--manifold with nonempty boundary.
Let $\mathcal{F}_1$ and $\mathcal{F}_2$ be two characteristic foliations on $\bdry M$ which
are adapted to $\Gamma_{\bdry M}$.
There exists a bijection $$\phi_{12}\co \pi_0(\mbox{Tight}(M,\mathcal{F}_1)) \rightarrow
\pi_0(\mbox{Tight}(M,\mathcal{F}_2)).$$\rm

\proof
The map $\phi_{12}$ is defined as follows:  Given any tight contact structure $\xi$ in
$\mbox{Tight}(M,\mathcal{F}_1)$, take an invariant neighborhood $\Sigma\times [0,\infty)\subset M$ for $\xi$,
where $\Sigma\times\{0\}=\bdry M$.  Take a parallel copy $\Sigma_k=\Sigma\times \{k\}$, for some large $k$.
Apply Giroux's Flexibility Lemma, with contact vector field ${\partial \over \partial t}$, where $t$ is the coordinate
for $[0,\infty)$. Starting with $\Sigma_k$ we obtain $\Sigma'\subset \Sigma\times (0,\infty)$
with characteristic foliation $\mathcal{F}_2$,
after a ${\bdry \over \bdry t}$--admissible isotopy (provided $k>>0$).  $\Sigma'$ divides
$M=M_1\cup M_2$, where $M_1\subset \Sigma\times[0,\infty)$.  We simply set $\phi_{12}(\xi)=\xi|_{M_2}$, where
$M_2$ is identified with $M$ and $\Sigma'$ is identified with $\Sigma$ via the flow of ${\bdry \over \bdry t}$.
$\phi_{12}$ does not depend on $k$, since we are considering contact 2--plane fields up to isotopy.
We now prove that $\phi_{12}$ is independent of the choice of contact vector field $X$.  Take a 1--parameter
family of contact vector fields $X_s$, $s\in[0,1]$, which are transverse to $\Sigma$.
Altering our perspective, this is equivalent to a 1--parameter family of ${\partial \over \partial t}$--invariant
contact 1--forms $\alpha_s$, $s\in[0,1]$, on $\Sigma\times[0,\infty)\subset M$.  The independence of the choice
of vector field then follows from observing that the proof of the Flexibility Lemma also applies to a family of
$\Sigma\times \R$'s.
We now show that $\phi_{21}$ is the inverse of $\phi_{12}$.  Refer again to $\Sigma\times[0,\infty)$
for $\xi\in \mbox{Tight}(M,\mathcal{F}_1)$.  Since ${\partial \over \partial t}$ is also a transverse contact
vector field for $\Sigma'$,  for sufficiently large $k'>>0$, $\Sigma_{k'}\subset M_2$.  Finally observe that
$\xi|_{M\backslash (\Sigma\times[k,k'])}$ is isotopic to $\xi$ itself.
\qed

\medskip
In view of the proposition, we will often write  $\mbox{Tight}(M,\Gamma)$ to stand for any of the
$\mbox{Tight}(M,\mathcal{F})$, where $\mathcal{F}$ is adapted to $\Gamma$.

\subsubsection{Standard neighborhoods of Legendrian curves} Let $\gamma\subset M$ be
a Legendrian curve with a negative twisting number $t(\gamma)=n$ with respect to a fixed framing.
The {\it standard tubular neighborhood} $N(\gamma)$ of a Legendrian curve $\gamma$
with $t(\gamma)$ negative is defined to
be $S^1\times D^2$ with coordinates $(z,(x,y))$ and contact 1--form
$\alpha=\sin (2\pi nz)dx+\cos(2\pi nz)dy$.  Here $\gamma=\{(z,(x,y))|x=y=0\}$.  With respect to this fixed
framing, we may identify $\bdry (N(\gamma))=\R^2/\Z^2$  by letting the meridian correspond
to $\pm (1,0)^T$ and the longitude (from the framing) correspond to $\pm (0,1)^T$.
With this identification, $s(\bdry (N(\gamma)))=-{1\over n}$.    On the other hand, we have the following
proposition, which is used by Kanda in \cite{K}, and is essentially proved in Makar--Limanov \cite{ML},
although phrased a bit differently.

\begin{prop} \label{prop:legnbhd} There exists a unique tight contact structure on $S^1\times D^2$
with a fixed convex boundary with $\#\Gamma_{\bdry (S^1\times D^2)}=2$  and slope
$s(\bdry (S^1\times D^2))=-{1\over n}$, where $n$ is a negative integer.  Modulo modifying the characteristic
foliation on the boundary using the Flexibility Lemma,  the tight contact structure is
isotopic to the standard neighborhood of a Legendrian curve with twisting number $n$.
\end{prop}

\proof
Using Proposition \ref{prop:flex}, we may assume that $T^2=\bdry(S^1\times D^2)$ has Legendrian ruling curves of
slope $0$.  Take a meridional disk $D$ with one Legendrian ruling curve $L$ on the boundary.  There exists
a collar annulus $A=L\times [0,1]$ with $L=L\times \{0\}$ transverse to $T^2$ along $L$.   Using
Proposition \ref{prop:makeconvex}, we may perturb $D$ to be convex with collared Legendrian boundary.
Since $t(L,Fr_{D})=-1$, there exists a unique dividing set $\Gamma_D$ consisting of one
arc from $\bdry D$ to $\bdry D$.  Using the Flexibility Lemma, we find that any $D$ can be normalized to
have a particular chosen characteristic foliation with this dividing set.  Given any two tight contact structures
$\xi_1$ and $\xi_2$ on $S^1\times D^2$ with given boundary condition, we may match them up along
$T^2\cup D$, after an isotopy (not necessarily contact).  The rest is a 3--ball $B^3=(S^1\times D^2)
\backslash (T^2\cup D)$ (after edge-rounding),
and we find an isotopy which matches $\xi_1$ and $\xi_2$ on $B^3$ fixing $\bdry B^3$, using Eliashberg's
Theorem  (Theorem \ref{thm:3-ball}).
\qed

\medskip
The following is a useful lemma:

\begin{lemma} [Twist Number Lemma] Let $(M, \xi)$ be a tight manifold with
a fixed framing $\mathcal{F}$.  Consider a Legendrian curve
$\gamma$ with $t(\gamma, Fr) = n, n \in {\bf Z}$,  and a standard tubular neighborhood $V$ of
$\gamma$ with boundary slope ${1\over n}$.  If there exists a bypass $D$ which is attached along a
Legendrian ruling curve of slope $r$, and ${1\over r}\geq n+1$, then
there exists a Legendrian curve with larger twisting number  isotopic (but not
Legendrian isotopic) to $\gamma$.
\end{lemma}

\proof Follows immediately from Lemma \ref{slopes}.\qed

\medskip
Notice that from this perspective the notion of {\it destabilization} due to Etnyre \cite{Et2} is basically identical
to our notion of a bypass.

\subsection{Relative Euler class}  Consider a tight contact structure $\xi$ on a manifold
$M$ with convex boundary $\bdry M$.   Assume $\xi|_{\bdry M}$ is trivializable, and choose
a nowhere zero section $s$ of $\xi$ on $\bdry M$.  Then we may form the {\it relative Euler class}
$e(\xi,s)\in H^2(M,\bdry M;\Z)$.  Consider the following exact sequence:
$$\begin{array}{ccccccc}
H^1(\bdry M) & \rightarrow & H^2(M,\bdry M) & \rightarrow & H^2(M) & \rightarrow & H^2(\bdry M)\\
& & e(\xi, s) & \mapsto & e(\xi) &\mapsto & 0
\end{array}   $$
This implies that a nonzero section $s$ of $\bdry M$ allows for a lift of $e(\xi)$ to $e(\xi,s)$.     Given
two nonzero sections $s_1$ and $s_2$, $e(\xi,s_1)$ and $e(\xi,s_2)$ will differ by an element which
is represented in $H^1(\bdry M)=Map (\bdry M,S^1)$.  The relative Euler class can be evaluated as follows:

\begin{prop}
Let $(M,\xi)$ be a contact manifold with convex boundary.  Fix a nonzero section $s$ of $\xi|_{\bdry M}$.
\be
\item If $\Sigma\subset M$ is a closed convex surface with positive (resp. negative) region $R_+$ (resp. $R_-$)
divided by
$\Gamma_\Sigma$, then $\langle e(\xi),\Sigma\rangle =\langle e(\xi,s),\Sigma\rangle=\chi(R_+)-\chi(R_-)$.
\item If $\Sigma\subset M$ is a compact convex surface with Legendrian boundary on
$\bdry M$ and regions $R_+$ and $R_-$, and $s$ is homotopic to $s'$ which coincides with
$\dot \gamma$ for every oriented connected component $\gamma$ of $\bdry \Sigma$,
then $\langle e(\xi, s),\Sigma\rangle=\chi(R_+)-\chi(R_-)$.
\ee
\end{prop}

\proof (1) follows from perturbing $\Sigma$ while fixing $\Gamma$ so that $\Sigma$ is singular Morse--Smale.
Then use a standard computation which says that $\langle e(\xi,s),\Sigma\rangle =d_+-d_-$,
where $d_\pm=e_\pm-h_\pm$, $e_+$ (resp. $e_-$) is the number of positive (resp. negative) elliptic points,
and $h_+$ (resp. $h_-$) is the number of positive (resp. negative) hyperbolic points.
(2) is almost identical.  The only difference is that the half-elliptic and half-hyperbolic points must be counted
properly.  This is done in Kanda \cite{K98}. \qed

\medskip
Let $T=\R^2/\Z^2$ be a component of a convex  $\bdry M$, where $\xi$ is tight.
Then, by the Flexibility Theorem, we may assume $T$ is in standard form with slope $s(T)$ and
Legendrian rulings with
slope $r$.   Take a nonzero section $s$ of $\xi|_T$ given by the tangent field of the rulings.  Let $\Sigma$ be
a compact surface with boundary along $T$.
Starting with $T_0=T$, there exists a 1--parameter family $T_t$, $t\in[0,1]$, of convex surfaces as in the
Flexibility Theorem, so that $T_1$ is in standard form and has Legendrian rulings of slope $r'$.
By excising and viewing $T_t$ as the new boundary of $M$, we obtain a 1--parameter family of
contact structures $\xi_t$ with $\xi_0=\xi$.  If we take $s'$ given by the tangent field of Legendrian rulings
of slope $r'$ on $T_t$,  then
$$\langle e(\xi,s),\Sigma\rangle=\langle e(\xi_1,s'),\Sigma\rangle,$$
since the relative Euler class remains invariant under homotopy.   This proves:

\begin{lemma}
Let $(M,\xi)$ be a tight contact manifold with convex boundary consisting of tori.
Then the relative Euler class $\langle e(\xi,s),\Sigma\rangle$ is independent of
the slope of the Legendrian rulings, if $s$ is a nonzero section of $\xi$ on a perturbation of $\bdry M$, given
by the tangent field of the rulings.
\end{lemma}

We now explain how to compute relative Euler classes for spaces of interest.  Assume $\xi$ is tight.
For $S^1\times D^2$ with convex boundary, use the Flexibility Theorem to make the Legendrian
rulings horizontal, take $s$ to be tangent to $\bdry (S^1\times D^2)$, and compute the relative Euler class
of a meridional convex disk $\Sigma$ with Legendrian boundary by taking $\langle e(\xi,s),\Sigma\rangle
=\chi(R_+)-\chi(R_-)$.    $e(\xi,s)\in H^2(M,\bdry M;\Z)=H_1(M;\Z)\simeq \Z$, so evaluation on a single meridional
disk completely determines the relative Euler class.

Similarly, for $T^2\times I$, we modify the boundary so the Legendrian rulings have the same slope $r$ for
both $T^2\times\{0\}$ and $T^2\times\{1\}$.   Take a convex annulus $A=\gamma\times I$ with Legendrian
boundary, where $\gamma$ is a closed
curve with slope $r$.   If we compute   $\langle e(\xi,s),A\rangle
=\chi(R_+)-\chi(R_-)$ for two annuli of two different slopes, then this determines the element
$e(\xi,s)\in H_1(T^2\times I;\Z)\simeq H_1(T^2;\Z)\simeq \Z^2$.

\subsubsection{Computation when $\bdry (T^2\times I)$ is nonsingular Morse--Smale}
We explain how to relate the relative Euler class computations in the two settings:
when $T_i=T^2\times \{i\}$, $i=0,1$, have nonsingular Morse--Smale characteristic foliations
versus when $T_i$ are in standard form.  We assume the dividing sets are unchanged when switching between
cases.  In the Morse--Smale case, take the nonzero section $s'_0$ given by the nonsingular
flow on $T_i$, or, equivalently, a nonzero section $s'$ of $\xi$ which is everywhere transverse to $T_i$.
In the standard form situation, take the nonzero section $s_0$ given by the Legendrian rulings,
or, equivalently, a nonzero section $s$ of $\xi$ which is transverse to the rulings and twists along the
ruling curves.    By comparing $s$ and $s'$, we see that `$s-s'$' is given by $\pm n\cdot PD(v_i)\in
H^1(T_i;\Z)$,
where $v_i$ is the shortest integral vector with slope $s(T_i)$ and $n$ is the torus division number.

\subsection{Basic slices}

In what follows, we will fix an identification $T^2=\R^2/ \Z^2$.
Let $T^2 \times I=\R^2/\Z^2
\times [0, 1]$ with coordinates  $(x, y, z)$, and $T_s=T^2\times\{s\}$, $s\in[0,1]$.
Recall the {\it boundary slope} $s_i=s(T_i)$
is  the slope of the dividing curves on $T_i$ (defined only when $T_i$ is convex).
We will call $(T^2\times I,\xi)$ a {\it basic slice} if
\be
\item $\xi$ is tight.
\item $T_i$ are convex and $\#\Gamma_{T_i}=2$, for $i=0,1$.
\item The minimal integral representatives of $\Z^2$ corresponding to $s_i$ form a
$\Z$--basis of $\Z^2$.
\item $\xi$ is {\it minimally twisting}, as defined in Section \ref{section:twisting}.
\ee

After a diffeomorphism of $T^2$, we may assume that a basic slice has boundary slopes
$s_1=-1$ and $s_0=0$.       Denote the subset of minimally twisting tight contact structures
in $\mbox{Tight}(T^2\times I,\mathcal{F})$ by  $\mbox{Tight}^{min}(T^2\times I, \mathcal{F})$.

\begin{prop} \label{basic}   Let $\Gamma_{T_i}$, $i=0,1$, satisfy $\#\Gamma_{T_i}=2$ and $s_1=-1$, $s_0=0$.
Then $|\pi_0(\mbox{Tight}^{min}(T^2\times I, \Gamma_{T_1}\cup \Gamma_{T_2}))|=2$.    (Here $|\cdot|$
denotes cardinality.)
The two tight contact structures are universally tight, and the Poincar\'e duals to the
relative Euler classes are given by
$\pm(0,1)\in H_1(T^2;\Z)$.
\end{prop}

\proof We will prove this proposition in steps.

\medskip

{\bf Step 1}\qua We will show that $|\pi_0(\mbox{Tight}^{min}(T^2\times I,
\Gamma_{T_1}\cup \Gamma_{T_2}))|\leq 2$. Assume the contact structure $\xi$ is tight.
Take $\Gamma_{T_i}$ to have $\#\Gamma_{T_i}=2$ and
$s_1=-1$, $s_0=0$, and choose $\mathcal{F}_i$ adapted to $\Gamma_{T_i}$, $i=0,1$, so that the
Legendrian rulings for both $T_i$ are vertical.  Take a vertical annulus
$A=\{0\}\times (\R/\Z)\times[0,1]$, whose boundary consists of two Legendrian
ruling curves,
each with twisting number $-1$  relative to $T_i$. After perturbation, $A$ is convex with collared Legendrian
boundary.   Assume that the endpoints of $\Gamma_A$ are $\{0\}\times\{0,{1\over 2}\}\times \{0,1\}$.

\begin{claim} All the dividing curves on $A$ must
connect from $T_0$ to $T_1$, ie, there are no boundary-parallel dividing curves.
\end{claim}

{\bf Proof of Claim}\qua
Otherwise, we obtain a singular bypass for
$T_0$ attached along a vertical Legendrian ruling curve by using the Imbalance Principle.
Using the Pivot Lemma, we smooth this bypass curve into a Legendrian curve $\gamma$  which has slope
$\infty$ when linearized.  There exist $T_{1\over 2}\supset \gamma$ for which the twist number $t(\gamma)$
is zero with respect to $T_{1\over 2}$.  Perturbing $T_{1\over 2}$ into a convex surface, we find that
$s(T_{1\over 2})=\infty$.  Therefore, this contradicts the assumption that $\xi$ is minimally twisting.
\qed

\medskip
Although the dividing curves connect from $T_0$ to $T_1$ and are parallel,
there are still infinitely many possible configurations for $\Gamma_A$, distinguished by
the {\it holonomy}.    We can define the {\it holonomy} $k(A)$
as follows:  pass to the cover $\{0\}\times \R\times I\subset S^1 \times \R \times I$ and let $k(A)$
be the integer such that there is a dividing curve which connects from $(0, 0, 0)$ to $(0, k(A),
1)$.

\begin{claim}
The holonomy function $k\co \mathcal{A}\rightarrow \Z$ is surjective, where $\mathcal{A}$ is the set of
convex annuli which have the same boundary as $A$ and are isotopic to $A$.
\end{claim}

{\bf Proof of Claim}\qua We explain how to apply {\it sliding} to modify $A$ to $A'$ (with the same
boundary) so that $k(A')=k(A)\pm 1$.  This would then imply the surjectivity.
Assume $A=\{0\}\times S^1\times [0,1]$ is convex with Legendrian boundary and holonomy $k(A)$.
Let $N(A)=[-\varepsilon,\varepsilon]\times S^1\times [0,1]$ be an $I$--invariant neighborhood of $A$.
Take $A'=(\{-\varepsilon\}\times S^1\times [\varepsilon,1-\varepsilon] ) \cup
((S^1\backslash (-\varepsilon,0))\times S^1\times
\{\varepsilon,1-\varepsilon\}) \cup (\{0\}\times S^1\times ([0,\varepsilon]\cup[1-\varepsilon,1]))$
and round the edges using the Edge-Rounding Lemma.
Informally we are adjoining copies of $T_0$ and $T_1$ which are cut open along
$\bdry A$, and rounding. $k(A')=k(A)+1.$  We can obtain $A'$ with $k(A')=k(A)-1$ similarly. \qed

\medskip
Therefore, after an isotopy fixing its boundary, $A$
can be put into standard form, with $k(A)=0$ and vertical
Legendrian rulings.  Now cut along $A$ to obtain
$S^1\times D^2$ with boundary slope $s(\bdry(S^1\times D^2))=-2$
and vertical Legendrian rulings, after rounding the edges.

\begin{figure}[ht!]
{\epsfysize=2in\centerline{\epsfbox{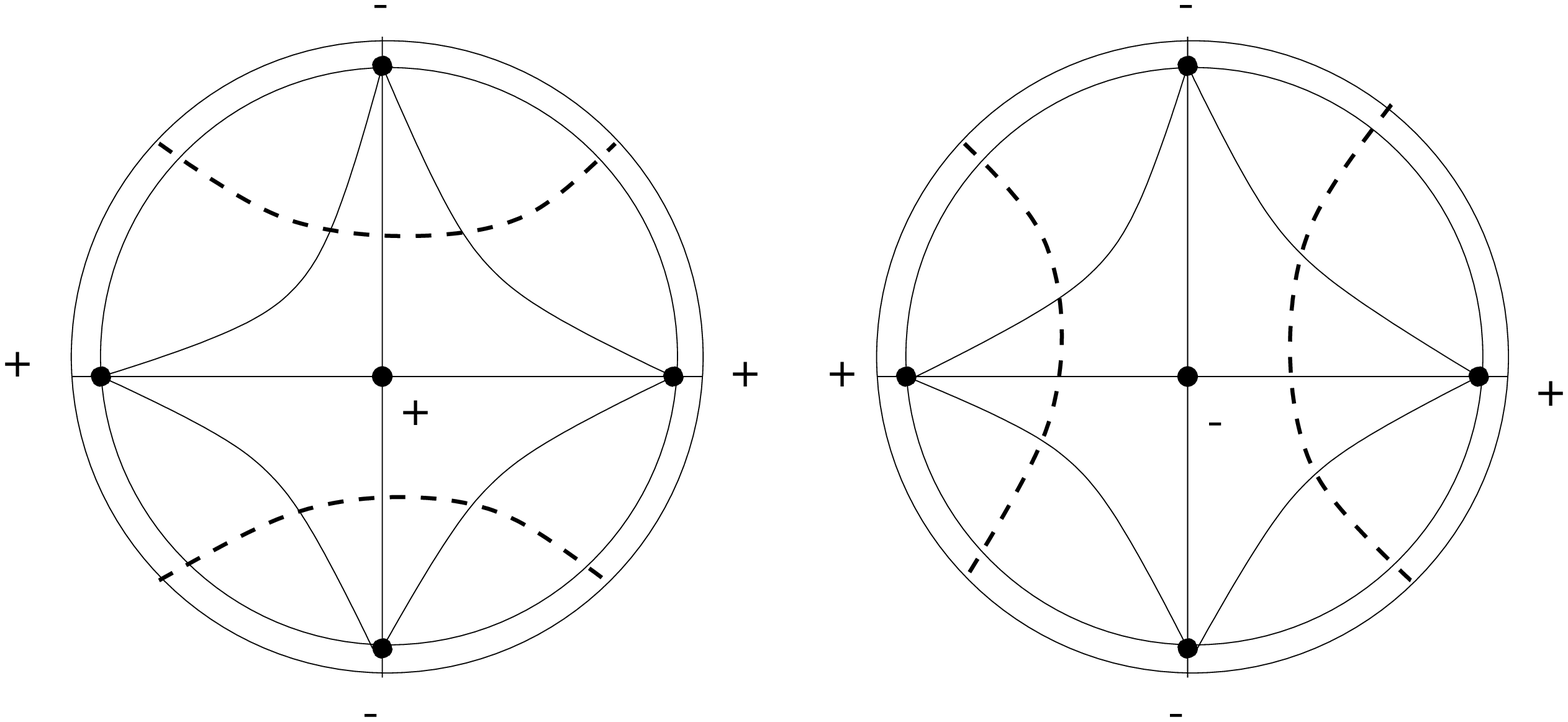}}}
\caption{Two possibilities on $D^2$ with $t(\bdry D)=-2$:  The dotted lines are dividing curves.}
\label{disk2}
\end{figure}

Next, using the Flexibility Lemma, we make the Legendrian rulings horizontal, and
take a meridional disk $D^2$ of the solid torus, which we assume is convex with
collared Legendrian boundary.  There are two possible
configurations of dividing curves, pictured in Figure \ref{disk2}.
Now, given two tight contact structures $\xi_1$ and $\xi_2$ on $T^2\times I$ with the given boundary
conditions, $\xi_1$ and $\xi_2$ can be isotoped so that they agree on $T_0\cup T_1\cup A$.  If
the $\Gamma_D$ are isotopic, then $\xi_1$ and $\xi_2$ can be matched up on $D$ in addition, and
Eliashberg's theorem (Theorem \ref{thm:3-ball}) implies that $\xi_1$ and $\xi_2$ are contact isotopic rel the
boundary. Therefore
we have at most two tight structures on a basic slice up to an
isotopy which fixes the boundary.

\medskip

{\bf Step 2}\qua  Let us compute the relative Euler class.  We already found that if the Legendrian rulings
were made to have slope $r=\infty$, then the annulus $A=\gamma\times I$ with $\gamma$ a closed
curve of slope $\infty$ satisfies $\langle e(\xi,s),A\rangle =0$.
We now compute $\langle e(\xi,s),B\rangle$, for the annulus
$B=\gamma\times [0,1]$, where $\gamma$ is a closed curve with slope $1$.  Here the Legendrian rulings
for $T_0$, $T_1$ have slope 1, and $B$ is a convex surface with Legendrian boundary (we have fixed
an orientation for $B$).   Write $\gamma_i
=\gamma\times\{i\}$, $i=0,1$.
$t(\gamma_0)=-1$ and $t(\gamma_1)=-2$ with respect to $B$, so there
exists a boundary-parallel dividing curve on $B$ along $\gamma_1$ by the Imbalance Principle.
We argue as in Step 1 to show that (1) two of the dividing curves on $B$ must go across from $\gamma_0$
to $\gamma_1$; otherwise minimal twisting is violated, (2) we may normalize the holonomy $k(B)$ of the
two dividing curves which go across, and (3) once $B$ is normalized, the cut-open solid torus has boundary
slope $-1$, hence is unique.  Therefore, we find that  $\langle e(\xi,s),B\rangle=\pm 1$, and
$PD(e(\xi,s))=\pm (0,1)$.

\medskip

{\bf Step 3}\qua
The two possible candidates for tight structures on the basic slice
are tight (and even universally tight).  We find an explicit model which
can be embedded in $(T^3,\xi_1)$, where $T^3=\R^3/\Z^3$ has coordinates
$(x,y,z)$ and  $\xi_1$ is
given by the 1--form $\alpha_1=\sin(2\pi z)dx+\cos(2\pi z)dy$.
We can choose $T^2\times [0, {1\over 8}]\subset T^3$, and perturb the boundary
so that $\#\Gamma=2$ and in standard form for both boundary components, with boundary slopes
$s_{{1\over 8}} = -1$ and
$s_0 = 0$.  If we rotate this tight structure by $\pi$, then
we obtain the other candidate.
Although not isotopic (distinguished by the relative Euler class), the two tight structures are diffeomorphic via a
diffeomorphism
isotopic to $-id$, where $id$ is the identity map on $T^2=\R^2/\Z^2$.   The relative Euler class
is computed by perturbing the boundary of $T^2\times [0,{1\over 8}]$ so the characteristic foliation is
Morse--Smale. The annulus $B=\gamma\times [0,{1\over 8}]$ with transverse boundary,
where $\gamma$ has slope $1$, will give
$\langle e(\xi,s'),B\rangle=0$ if $s'$ is tangent to the boundary.  Converting this to $s$ which is tangent to
the Legendrian rulings for the characteristic foliation in standard form, we obtain
$\langle e(\xi,s),B\rangle=\pm 1$.

\medskip

{\bf Step 4}\qua
It remains to show that the tight structure on $N=T^2\times [0, {1\over 8}]
\subset T^3$
is minimally twisting. Assume the existence of a torus $T'\subset N$
parallel to
$T_{{1\over 8}}$ and $T_0$, for which the boundary slope $s'$ is not
between $-1$ and $0$.
This is equivalent to the existence of a linear Legendrian curve
$\gamma_0\subset N$ with slope
$s'$  and $t(\gamma_0,Fr_{T^2})=0$.
We will pass to the universal cover $(\widetilde{N}=\R^2\times [0, {1\over
8}],
\widetilde{\xi}_1)$ to find an
overtwisted disk.

Assume $s'>0$.  Pick a point $p=(x_0,y_0,z_0)$ on
$\gamma_0$ with the smallest $z$--coordinate, and view $\gamma_0$
as starting and ending at $p$.  A lift $\widetilde{\gamma}_0$ will have
endpoints $\widetilde{p}_1
=(x_1,y_1,z_0)$,
$\widetilde{p}_2=(x_2,y_2,z_0)$ which are lifts of $p$. Let
$\widetilde{\gamma}_1$ be
the linear Legendrian curve from $(x_1,y_1,0)$ to $(x_1,y_1,z_0)$,
$\widetilde{\gamma}_2$
be the linear Legendrian curve from $(x_2,y_2,z_0)$ to $(x_2,y_2,0)$, and
$\widetilde{\gamma}_3$
be the linear Legendrian curve from $(x_2,y_2,0)$ to $(x_1,y_2,0)$.  Then
the composite
$\widetilde{\gamma}=\widetilde{\gamma}_1 + \widetilde{\gamma}_0+
\widetilde{\gamma}_2+\widetilde{\gamma}_3$ is
a Legendrian curve which projects to a closed curve onto the $xz$--plane and has
positive holonomy.  It is easy to decrease its holonomy by adding
a curve $\widetilde{\gamma}'$ which projects to  $\gamma'$ in the
$xz$--plane and satisfies
$\gamma'=\bdry\Omega$,
where $\Omega$ is a region in the $xz$--plane.  Therefore, we obtain an
overtwisted disk bounded by
$\widetilde\gamma+\widetilde\gamma'$.  Notice that $t(\gamma_0)=0$
translates to
$tb(\widetilde\gamma+\widetilde\gamma')=0$.
We argue similarly for $s'<-1$, and find that $N$ is minimally
twisting.
\qed

\medskip
We also have the following corollary.  A {\it pre-Lagrangian torus} is a torus with linear characteristic foliation.

\begin{cor}   \label{allslopes}
Let $(T^2\times I,\xi)$ be a basic slice, with boundary slopes $s_0$ and $s_1$.
Then for any slope $s$ {\it between}  $s_1$ and $s_0$
(see Section \ref{section:twisting} for the definition), there exists a convex torus $T$
parallel to $T^2\times\{pt\}$ with slope $s(T)=s$.    For any slope $s$ between $s_1$ and $s_0$
(but $\not=s_0,s_1$), there exists a pre-Lagrangian torus  $T$ parallel to $T^2\times \{pt\}$ with slope
$s$.
\end{cor}

\proof  This follows from the explicit model in the proof of Proposition \ref{basic}.  A basic slice will
have pre-Lagrangian tori of all slopes between $s_1$ and
$s_0$, and any pre-Lagrangian torus can be perturbed into a convex torus with the same slope.\qed

\subsection{Decomposition of $T^2\times I$ into layers}

Assume that $\xi$ on $T^2\times I$ is tight.  In this section we will also
assume the following: (1) $\#\Gamma_{T_i}=2$, $i=0,1$, (2) $\xi$
has minimal twisting.  It is most convenient to arrange
the boundary slopes, via an action of $SL(2,\Z)$, as follows:
$-\infty<s_1\leq -1$ and $s_0=-1$.
Write $s_1=-{p\over q}$, where
$p\geq q>0$ are integers and $(p,q)=1$.

\subsubsection{Nonrotative case}

\begin{prop}  [Nonrotative case]           \label{nonrotative}
Let $\Gamma_{T_i}$, $i=0,1$, satisfy $\#\Gamma_{T_i}=2$ and $s_0=s_1=-1$.  Then
there exists a holonomy map $k:
\pi_0(\mbox{Tight}^{min}(T^2\times I, \Gamma_{T_1}\cup \Gamma_{T_2})) \rightarrow \Z$ which is
bijective.
\end{prop}

\proof Use the Flexibility Lemma to obtain rulings of slope $r_0=r_1=0$, take a horizontal annulus
$S^1\times \{0\}\times I$ with Legendrian boundary, and perturb it into a convex surface with collared
Legendrian boundary. If the dividing curves of the annulus $A$ do not cross from $T_0$ to $T_1$,
then, by Proposition \ref{bcc},  there exists a boundary-parallel dividing curve on $A$ along $T_1$,
and the corresponding singular bypass gives rise to a factoring $T^2\times[0,1]=T^2\times[0,{1\over 2}]
\cup T^2\times [{1\over 2},1]$, where the intermediate layer
$T_{1\over 2}$ is convex with slope $s_{1\over 2}=0$.  This contradicts our minimal twisting assumption.
Therefore, both dividing curves on $A$ cross from $T_0$ to $T_1$.     Put $A$ in standard form,
cut along $A$, and perform Edge-Rounding to obtain a solid torus with boundary slope $-1$.    There
exists a unique tight contact structure on this solid torus by  Proposition \ref{prop:legnbhd}.  This implies that,
for every choice of $\Gamma_A$, there exists at most one tight contact structure.

Now define the {\it holonomy} $k(A)$  by passing to the cover $\R\times\{0\} \times I \subset \R\times S^1\times I$
and letting $k(A)$ be the integer such that there exists a dividing curve connecting from
$(0,0,0)$ to $(k(A),0,1)$ (assume that the endpoints of all the possible dividing curve configurations are fixed).
To write down a tight contact structure $\xi_0$ with $k(A)=0$, simply take the $I$--invariant neighborhood of a
convex $T^2$ with $\#\Gamma=2$, $s(T^2)=-1$, and horizontal Legendrian rulings.  If we take $\xi_0$
and isotoped $T_1$ via $(x,y)\mapsto (x, y+k)$ ($k\in\Z$), while fixing $T_0$, then we obtain $\xi_k$ with $k(A)=k$.
The $I$--invariant contact structure is embeddable into a basic slice, hence it is universally tight.  Moreover, since the
basic slice is minimally twisting, so is the $I$--invariant tight structure.

We claim that $k(A)$ takes constant values in $\mathcal{A}$, the set of convex annuli which have the
same boundary as $A$ and are isotopic to $A$, provided $\xi$ is fixed. Assume $A'\in
\mathcal{A}$ with $k(A')\not=0$ (assume $k(A)=0$).   The proof follows a
strategy due to Kanda \cite{K}.   The strategy is to pass to
$\widetilde M=S^1\times [-n,n]\times I$  ($n$ large) and pick nonintersecting lifts
$\widetilde{A'}$ and $\widetilde{A}$ of $A'$ and $A$.  If
$k(A')>0$, then take $\widetilde {A'}$ to lie above $\widetilde A$.    Pick $N\subset \widetilde M$ bounded
above by $\widetilde{A'}$ and below by $\widetilde{A}$, and round the edges.  We find that
the boundary slope of the rounded $N$ is $0$ or a positive integer. If the slope is zero,
we have an overtwisted disk. Assume the slope is a positive integer.   Make $N$ have horizontal Legendrian
rulings, and take a convex meridional disk $D$ with a Legendrian collar boundary.  There exists a bypass by
Lemma \ref{bccdisk}, and, after bypass attachment, the slope is $\infty$.  Therefore, $N$ is the
standard neighborhood of a Legendrian curve $\gamma$ isotopic to $S^1\times \{0\}$ with twist number $0$, and
$M\supset N$ is the standard neighborhood of a Legendrian curve with twist number $-1$.
This is a contradiction.
\qed

\medskip

Although the holonomy gives infinitely many tight contact structures up to isotopy (fixing the boundary),
this turns out to be a special feature of the {\it nonrotative} case.
In the {\it rotative} case, Proposition \ref{basic} allows us to reduce the infinitely many possible dividing sets
to a finite collection.

\subsubsection{Rotative case}

\begin{prop} [Minimal twisting, rotative case]    \label{prop:rotative}
Let $\Gamma_{T_i}$, $i=0,1$, satisfy $\#\Gamma_{T_i}=2$ and $s_0=-1$, $s_1=-{p\over q}$, where
$p>q>0$.  Then
\begin{equation}\label{eqn6}
|\pi_0(\mbox{Tight}^{min}(T^2\times I, \Gamma_{T_1}\cup
\Gamma_{T_2}))| \leq  |(r_0+1)(r_1+1)\cdots (r_{k-1}+1)(r_k)|,
\end{equation}
where $-{p\over q}$ has a
continued fraction expansion $$-{p\over q}=r_0-{1\over r_1-{1\over r_2-\cdots {1\over r_k}}},$$
with all $r_i<-1$ integers.
\end{prop}

The proof will consist of a factorization $T^2\times I =\bigcup _{i=0}^k (T^2\times [{i\over k},{i+1\over k}])$
where $T_{i\over k}$, $i=0,\cdots k$, are convex with $\#\Gamma_{T_{i/k}}=2$ and slopes $s_{i\over k}$
arranged as
$s_0>s_{1\over k} > s_{2\over k}> \cdots > s_{k\over k}=s_1$;  this is followed by a shuffling argument
which reorders the layers. This is sufficient to prove the upper bound in the proposition. The proof
will occupy the next three sections. To prove
that the upper bound is exact requires embeddings into lens spaces.
This will be done in Section \ref{section:surg}.

\subsubsection{Factoring}
Take $r_1=r_0=0$ as before, and consider the
horizontal annulus $A$. Since $t(S^1\times \{0\}\times \{1\})=-p
<t(S^1\times \{0\}\times \{0\})=-1$, there must exist
a bypass along $T_1$.  Therefore, we can factor $T^2\times I$ into
$T^2\times [0,{1\over 2}]$ and $T^2\times [{1\over 2},1]$, where the latter
is a basic slice.  This follows from the following lemma:

\begin{lemma} \label{peeloff}
$T_{1\over 2}$ will have boundary slope $-{p'\over q'}$, where
$pq'-qp'=1$, $p>p'>0$, and $q\geq q'>0$.
\end{lemma}

\proof In order to use Lemma \ref{slopes}, we need to reflect $T_1$ and transform
via $SL(2,\Z)$ so that the boundary slope is 0.  Reflection gives us $-T_1$
with boundary slope ${p\over q}$ and rulings of slope $0$.  Then
$A_0=\left(\begin{array}{cc}
p' & -q' \\ p & -q
\end{array}
\right)$
sends $(q, p)^T\mapsto (-1, 0)^T$, $(1, 0)^T\mapsto (p', p)^T$.
Since $p>p'>0$, ${p\over p'}>1$,  the
boundary slope must be $\infty$ by Lemma \ref{slopes}.  Now, $A_0^{-1}:
(0,1)^T\mapsto (q',p')^T$,
and we have the lemma.\qed

\medskip
Applying  Lemma  \ref{peeloff} inductively, we obtain basic slices whose
boundary slopes increase from $-{p\over q}$ to $-1$ in a finite number of
steps.

\medskip
{\bf Example}\qua  Assume $s_1=-{10\over 3}$ and $s_0=-1$.  Then the boundary
slopes are $-{10\over 3}, -{3\over 1}=-3, -2, -1$, so we have a factorization
into 3 layers.

\medskip

\subsubsection{Continued fractions}
There exists a natural interpretation of the layering process in terms of
continued fractions.

Let $-{p\over q}$ have the following continued fraction expansion:
$$-{p\over q}=r_0-{1\over r_1-{1\over r_2-\cdots {1\over r_k}}},$$
with all $r_i<-1$ integers.  We identify $-{p\over q}$ with
$(r_0,r_1,\cdots,r_k)$.  Then $-{p'\over q'}$ as given in Lemma
\ref{peeloff} will correspond to $(r_0,r_1,\cdots,r_k+1)$, where
we identify $(r_0,\cdots,r_{k-1}+1)\sim (r_0,\cdots,r_k+1)$ if
$r_k=-2$.  This follows inductively from observing that if ${a\over b}$,
${a'\over b'}$
satisfy $ab'-ba'=1$, then $r-{1\over a/b}={ra-b\over a}$
and $r-{1\over a'/b'}={ra'-b'\over a'}$ satisfy
$$(ra-b)a'-(ra'-b')a=1.$$
Therefore, the boundary slopes of the factorization can be obtained in order
by decreasing the last entry of the corresponding continued fraction
expansion.

Notice that this layering process corresponds to taking a sequence
$-{p\over q} = -{p_0 \over
q_0} < -{p_1\over q_1}< \dots < -1$ where the consecutive slopes correspond
to pairs of vectors
which form an integral basis of ${\bf Z}^2$.  Moreover, the slopes on each
basic slice represent
a positive Dehn twist from the front face to the back face.  Therefore, we
have the following Factoring Lemma:

\begin{lemma}     \label{lemma:hops}
Let $\xi$ be a minimally twisting tight contact structure on $T^2\times I$.  Then $T^2\times I$ admits
a decomposition  $T^2\times I =\bigcup _{i=0}^k (T^2\times [{i\over k},{i+1\over k}])$,
where $T_{i\over k}$, $i=0,\cdots k$, are convex with $\#\Gamma_{T_{i/k}}=2$ and slopes $s_{i\over k}$.
The sequence of slopes is obtained by taking the shortest sequence of positive Dehn twists from
$-{p\over q}$ to $-1$.    Alternatively, in the tessellation picture, $s_1,s_{k-1\over k},\cdots, s_0$,
is the shortest sequence of hops along edges from $s_1$ to $s_0$, subject to the constraint that
$s_{i\over k}$ sit on the arc $[s_1,s_0]\subset \bdry \H^2$ (counterclockwise starting from $s_1$).
\end{lemma}

\subsubsection{Sliding maneuver}

There exists a natural grouping of the layers into blocks via continued fractions.
The blocks are isomorphic to $T^2\times I$ with minimal twisting, $\#\Gamma_{T_i}=2$, $i=0,1$,
and boundary slopes $s_1=-m$, $s_0=-1$, where $m\in \Z^+$, $m>1$.   Such blocks will
be called {\it continued fraction blocks}, and are special because the basic layers that comprise
a continued fraction block can be `shuffled'.

\begin{prop} \label{m}
Let $\Gamma_i=\Gamma_{T_i}$, $i=0,1$, be dividing sets satisfying $\#\Gamma_i=2$,
$s_0=-1$, $s_1=-m$, $m\in \Z^+$, $m>1$.
Then $|\pi_0(\mbox{Tight}^{min}(T^2\times I, \Gamma_0\cup \Gamma_1))| \leq m$.
\end{prop}

\proof  Let $(T^2\times I, \xi)$ have minimal twisting, $\#\Gamma_i=2$, slopes $s_1=-m$ and
$s_0=-1$, and coordinates $((x,y),z)$.  Consider a convex annulus
$A=S^1\times \{0\}\times I$ with Legendrian boundary (after perturbation of the boundary) and
oriented normal ${\bdry \over \bdry y}$. The minimal twisting condition guarantees
the existence of two dividing curves on $A$ which go across from $T_0$ to $T_1$, and
$m-1$ dividing curves from $T_1$ to itself.
Since there must be at least one bypass, we can peel off a layer and obtain a
basic slice $T^2\times [{m-2\over m-1}, 1]$ with $s_1 = -m$ and $s_{{m-2\over m-1}} = -(m-1)$.
The horizontal annulus from $T_{m-2\over m-1}$ to $T_1$ will have $2(m-1)$ dividing
curves which go across from $T_{m-2\over m-1}$ to $T_1$,
and 1 dividing curve from $T_1$ to itself, which is
the boundary-parallel curve used to peel off $T^2\times [{m-2\over m-1}, 1]$.
The tight structure on the basic slice is determined
by whether the half-disk separated by the boundary-parallel curve is
positive or negative. (Recall that $PD(e(\xi,s))=\pm (0,1)\in
H_1(T^2;\Z)$ by Proposition \ref{basic}.)  In a similar manner,
we successively peel off $T\times [{i-1\over m-1}, {i\over m-1}]$, with boundary slopes $-i$ and $-(i + 1)$.
Let us say that the layer $T\times [{i-1\over m-1}, {i\over m-1}]$ is positive (resp. negative) if the
sign of half-disk separated by the boundary-parallel dividing curve is positive (resp. negative).
The proof then follows from repeated applications of the following lemma.
\qed

\begin{lemma}[Shuffling]
Consider a minimally twisting tight $(T^2\times I,\xi)$ with $\#\Gamma_i=2$, $i=0,1$,
boundary slopes $s_1=-k$ and $s_0=-k+2$, $m\geq k,k-2\geq 1$.  Given a factorization
$T^2\times I=N_1\cup N_2$, $N_1=T^2\times [0,{1\over 2}]$, $N_2=T^2\times [{1\over 2},1]$,  into basic layers
($s_{1\over 2}=-k+1$), where $N_1$ is positive and $N_2$ is negative, there exists another factorization
$T^2\times I=N_1'\cup N_2'$ so that $N_1$ is negative and $N_2$ is positive.
\end{lemma}

\proof The lemma follows from applying the {\it sliding maneuver}, which we have already
used once to prove Proposition \ref{basic}.  Assume the Legendrian rulings for $T_0, T_{1\over 2}, T_1$
are all horizontal. Let $A_1=S^1\times \{0\}\times [0,{1\over 2}]$ and $A_2=S^1\times \{0\}\times [{1\over 2},1]$.
$A_1$ (resp. $A_2$) has $2(k-2)$ (resp. $2(k-1)$) dividing curves which go across and $1$
boundary parallel dividing curve along $T_{1\over 2}$ (resp. $T_1$).
Let $L=A_1\cap T_{1\over 2}=S^1\times \{0\}\times \{{1\over 2}\}$, and $\Gamma_{A_1}\cap L=
\{0,{1\over 2(k-1)}, {2\over 2(k-1)},\cdots, {2(k-1)\over 2(k-1)}\}\times \{0\}\times \{{1\over 2}\}$.
We modify $A_1$ to $A_1'$ by an isotopy which fixes $\bdry N_1$ so that the
boundary-parallel curve of $\Gamma_{A_1'}$ has endpoints which have been shifted
by $\pm {2\over 2(k-1)}$ along $L$.

\begin{figure}[ht!]
\centerline{\small
\psfrag {T1}{$T_{1\over2}$}
\psfrag {T0}{\smash{\rlap{\kern -4pt \raise -2pt\hbox{$T_0$}}}}
\psfrag {A}{$A_1$}
\psfrag {Intersection with}{\smash{\rlap{\kern -6pt \raise 
0pt\hbox{\scriptsize Intersection with}}}}
\epsfysize=2in\epsfbox{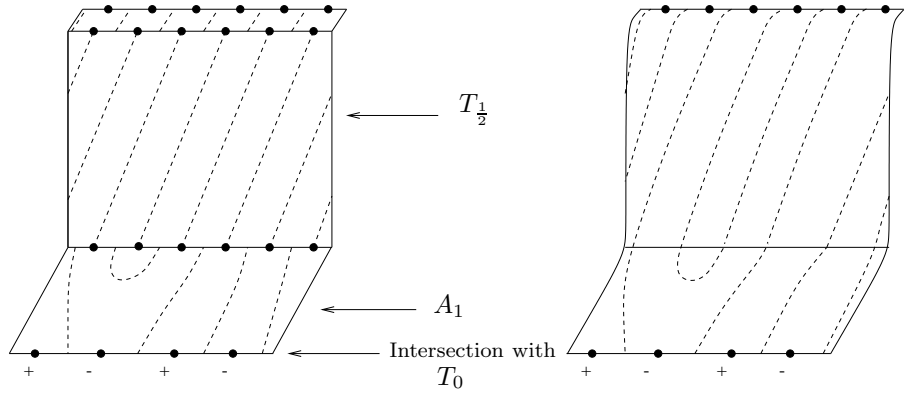}}
\caption{Sliding maneuver}
\label{swindle}
\end{figure}

Informally, we attach copies of $T_0$ and $T_{1\over 2}$ and round the edges.
Let $N(A_1)=S^1\times [-\varepsilon,\varepsilon]\times[0,{1\over 2}]$ be an $I$--invariant neighborhood
of $A_1$.  Take $A_1'=(S^1\times \{\varepsilon\}\times[\varepsilon,{1\over 2}-\varepsilon])\cup
(S^1\times (S^1\backslash (0,\varepsilon))\times\{\varepsilon,{1\over 2}-\varepsilon\})\cup
(S^1\times \{0\}\times ([0,\varepsilon]\cup[{1\over 2}-\varepsilon,{1\over 2}]))$, and round the edges
using the Edge-Rounding Lemma.
This moves the endpoints of the boundary-parallel curve by $-{2\over 2(k-1)}$ along $L$.
See Figure \ref{swindle} for an illustration.  Note that the copy of $T_0$ is not attached in this picture,
but we can still see that the  bypass has been slid along $L$.
Using the sliding maneuver, we may arrange $A_1\cup A_2$ so that the two dividing curves
with endpoints on $T_1$ are not nested, ie, they are both boundary-compressible dividing curves for
$A_1\cup A_2$.  We then have the freedom to choose which bypass to peel off first.
\qed

\medskip

{\bf Proof of Proposition \ref{prop:rotative}}\qua  We now group the layers of $T^2\times I$ with minimal boundary, minimal
twisting, and boundary slopes $-{p\over q}$ and $-1$ as follows: Act via
$A_0=\left(
\begin{array}{cc}
-r_0 & 1 \\ -1 & 0
\end{array}
\right)$.  Then $(1,-1)^T\mapsto (-r_0-1,-1)^T$  and $(q,-p)^T\mapsto (-r_0q-p, -q)^T$.
The boundary slopes are now $s_1={q\over r_0q+p}=r_1-{1\over {r_2\cdots}}$ and $s_0={1\over
r_0+1}$. Peel off a block with slopes $s_{1\over 2}=-1$ and $s_0={1\over r_0+1}$, which
is diffeomorphic to the form treated in Proposition \ref{m}, then continue. We will
then obtain $k$ blocks, each with minimal twisting, minimal boundary, and
boundary slopes $-1$, ${1\over r_i+1}$ (or, equivalently, $r_i+1$ and $-1$),
and one last block (at the very front) with boundary slopes $r_k$ and $-1$.
This completes the proof of Proposition \ref{prop:rotative}.\qed

\subsection{Factoring solid tori}

Let $(S^1\times D^2,\xi)$ be a solid torus with convex boundary $T$ and $\#\Gamma_T=2$. Fix
a framing $\mathcal{F}$ so that the boundary slope
$-{p\over q}$ satisfies $-\infty<-{p\over q}\leq -1$.    This is possible by normalizing via a suitable element of
$SL(2,\Z)$.
Here we view $T=\R^2/\Z^2$, where $(1,0)^T$ is the meridional circle and $(0,1)^T$ is the longitude with
respect to $\mathcal{F}$.

\begin{prop}   \label{torifactor}
Let $\Gamma_0$, $\Gamma_1$ be dividing sets on $T^2\simeq \R^2/\Z^2$ with $\#\Gamma_i=2$, $i=1,2$, and
slopes $s_0=-1$, $s_1=-{p\over q}$ ($-\infty<-{p\over q}\leq -1$).
Assume we have identified $\bdry (S^1\times D^2)\simeq \R^2/\Z^2$.
Let $\xi$ be a tight contact structure on $M\simeq S^1\times D^2$ with convex boundary
condition $\Gamma_1$.
Then there exists a factorization $M=N\cup(M\backslash N)$, where
$N$ is the standard neighborhood of a core Legendrian curve with twist number $-1$,
$M\backslash N\simeq T^2\times I$,
and $\xi|_{T^2\times I}$ is minimally twisting with boundary dividing sets $\Gamma_0$, $\Gamma_1$.
Hence we have
$$|\pi_0(\mbox{Tight}(S^1\times D^2, \Gamma_1))|   \leq
|\pi_0(\mbox{Tight}^{min}(T^2\times I, \Gamma_0\cup \Gamma_1))|$$
\end{prop}

\proof Let $\gamma$ be a Legendrian curve isotopic to the core $S^1$, satisfying
$t(\gamma)=-m$, $m\in \Z^+$. Such a Legendrian curve exists because any closed
curve $C'$ is $C^0$--close approximated by a Legendrian curve $C$ isotopic to $C'$.
Take a standard neighborhood $N'$ of $\gamma$ so that $\bdry N'$ is convex
with $s(\bdry N')=-{1\over m}$ and $\#\Gamma_{\bdry N'}=2$.  Now consider $M\backslash N'$ with
boundary slopes $-{p\over q}$ and  $-{1\over m}$.

We claim that the tight contact structure on $M\backslash N'$ is minimally twisting.  Assume otherwise.
Then there exists a factorization of $M\backslash N'=T^2\times I$ as $(T^2\times[0,{1\over 2}])\cup
(T^2\times [{1\over 2},1])$, where $s_0=-{1\over m}$, $s_1=-{p\over q}$, and $s_{1\over 2}$ is not between
$s_1$ and $s_0$.  Proposition \ref{prop:allslopes} below implies that there exists a convex torus
with any slope between $s_1$ and $s_{1\over 2}$ and any slope between $s_{1\over 2}$ and $s_0$.
In particular, $s_{1\over 2}=0$ is realized.  Now, a Legendrian divide on $T_{1\over 2}$ has twisting number
zero with respect to a meridional disk it bounds.  Hence $M\backslash N'$ is minimally twisting.

Since $-{p\over q}<-1$, the layering procedure for $M\backslash N'$ will give us a convex
torus $T'$ with boundary slope $-1$, parallel to $T$.  Factor $M=N\cup (M\backslash N)$ along $T$.
By Proposition \ref{prop:legnbhd} $N$ is a standard neighborhood of a Legendrian curve with twisting number
$-1$. Hence, the number of potential tight structures on $S^1\times D^2$ with
$\#\Gamma_{\bdry (S^1\times D^2)}=2$  and boundary slope $-{p\over q}$ is bounded above by the
number of minimally twisting tight contact structures on $T^2\times I$
with $\#\Gamma_{T_i}=2$ and  boundary slopes $s_1 = -{p\over q}$ and $s_0= -1$.    \qed

\begin{prop}      \label{prop:allslopes}
Let $(T^2\times I,\xi)$ be tight with convex boundary, and let $s_0$, $s_1$ be the boundary slopes.
Given any $s$ between $s_1$ and $s_0$, there exists a convex torus parallel to $T^2\times \{pt\}$ with
slope $s$.
\end{prop}

\proof   Let $s_0=-1$, $s_1=-{p\over q}$, with $p>q$ positive integers.  Let $T_0$, $T_1$ have
Legendrian rulings of slope $0$, and take a convex annulus $S^1\times \{0\}\times I$
with Legendrian boundary which are ruling curves of $T_i$.  There will exist a boundary-parallel dividing
curve, and if we attach the corresponding bypass we obtain a slope $-{p'\over q'}$ as in Lemma \ref{peeloff}.
After enough steps we arrive at a slope of $-1$.  Now, by Corollary \ref{allslopes}, any $s$ between $s_1$ and $s_0$
is represented by a convex torus.  \qed

\subsection{Lens spaces}

\subsubsection{Decomposition of lens spaces}
Consider the lens space $M=L(p,q)$, with $p>q>0$.  $L(p,q)$ is obtained
by gluing two solid tori $V_0$ and $V_1$ via $A_0\co  \bdry V_0\rightarrow \bdry V_1$ given by
$\left(
\begin{array}{cc}
-q & q' \\ p &-p'
\end{array}
\right)\in -1\cdot SL(2,\Z)$.
Here, $(1,0)^T$ is the meridional direction of $V_i$, and $(0,1)^T$ is the direction of the
core curve $C_i$ of $V_i$.  Note that $A_0$ is not unique --- we can compose $A_0$ with
Dehn twists to the left and the right.  However, we will fix a framing for $V_i$,
and assume $pq'-qp'=1$, $p>p'>0$ and $q\geq q'>0$.

\begin{prop}  \label{lensfactor}
Let $\Gamma_0,\Gamma_1$ be dividing sets on $\bdry (S^1\times D^2)\simeq \R^2/\Z^2$ with $\#\Gamma_i=2$,
$i=0,1$, and
slopes $s_0=-1$, $s_1=-{p'\over q'}$ ($-\infty<-{p'\over q'}\leq -1$).  Assume $-{p\over q}$ has continued fraction
representation $(r_0,\cdots, r_k)$ and $-{p'\over q'}$ has continued fraction representation
$(r_0,\cdots, r_k+1)$.
Then
\begin{eqnarray}|\pi_0(\mbox{Tight}(L(p,q)))|&\leq&|\pi_0(\mbox{Tight}(S^1\times D^2, \Gamma_1))| \\
&\leq & |\pi_0(\mbox{Tight}(T^2\times I, \Gamma_0\cup\Gamma_1))|\\
&\leq &|(r_0+1)(r_1+1)\cdots (r_{k-1}+1)(r_k+1)|.
\end{eqnarray}
\end{prop}

\proof The proof is very similar to Proposition \ref{torifactor}.  The goal is to thicken the
core Legendrian curve isotopic to $C_0$.  Note that
the meridional slope of $V_0$, when mapped to $\bdry V_1$, will have slope $-{p\over q}$ on
$\bdry V_1$. Let $\gamma$ be a Legendrian curve in $M = L(p, q)$,
isotopic to $C_0$, and with twisting number $n\leq 0$.  Recall it is
always possible to reduce the twisting number if necessary. Let $V_0$ to be the standard neighborhood of
$\gamma$ and $V_1=M\backslash V_0$.  Then $A_0$ maps $(n,1)^T\mapsto (-qn+q', pn-p')^T $,
and the corresponding boundary slope on $\bdry V_1$ is ${pn-p'\over -qn+q'}$.
Note that $-{p'\over q'}$ is the point on $\bdry \H^2$ with an edge in $\H^2$ to $-{p\over q}$ which is closest
to $-1$ on the arc $(-{p\over q},-1)\subset \bdry \H^2$. There exists a convex torus $T\subset V_1$ with
boundary slope $-{p'\over q'}$, using the factorization in Proposition \ref{torifactor} and Corollary \ref{allslopes}.
Modify $V_i$ so that $M$ is split along $T$ into $V_0$, $V_1$.  Now $n=0$ by Proposition
\ref{prop:legnbhd} and the boundary slope of $V_1$ is $-{p'\over q'}$.
Now we count the number of (possible) tight structures on
$V_1$ with $\#\Gamma_{\bdry V_1}=2$ and boundary slope $-{p'\over q'}$.  According to
Proposition \ref{prop:rotative},
an upper bound is given by $|(r_0+1)(r_1+1)\cdots (r_{k-1}+1)(r_k+1)|$, where
$(r_0,\cdots, r_k)$ is the continued fraction representation of $-{p\over q}$,
and $(r_0,\cdots, r_{k-1}, r_k+1)$ is the continued fraction
representation of $-{p'\over q'}$.        \qed

\medskip
Hence we have embedded a (candidate) minimally twisting tight contact structure on $T^2\times I$
as follows: $$T^2\times I\hookrightarrow S^1\times D^2\hookrightarrow L(p,q).$$
It remains to prove:

\begin{prop} \label{lowerbound}
$|\pi_0(\mbox{Tight}(L(p,q)))|\geq |(r_0+1)(r_1+1)\cdots (r_{k-1}+1)(r_k+1)|$, where
$-{p\over q}$ has continued fraction representation $(r_0,\cdots, r_k)$.   All the tight contact
structures in the lower bound are given by Legendrian surgery.
\end{prop}

The proof will be presented in the next section, after a discussion of Legendrian surgeries.
Observe that Proposition  \ref{lowerbound} together with Proposition \ref{lensfactor} prove
Theorems \ref{1} and \ref{3} as well as Part 2(a) of Theorem \ref{2}.

\subsubsection{Legendrian surgeries of $S^3$} \label{section:surg}

In this section we will realize all of the possible tight structures from the
previous sections inside Legendrian surgeries of links of unknots in $S^3$.
Recall the following theorem due to Eliashberg \cite{E90}.

\begin{thm}  \label{Legendriansurgery}
Let $K_1,\cdots,K_n$ be mutually disjoint Legendrian knots in the standard tight contact structure $\xi$ on $S^3$.
Then $M$, obtained from $B^3$ by $(tb(K_i)-1)$--surgery (usually called {\em Legendrian surgery})
along all the $K_i$, $i=1,\cdots,n$, is
holomorphically fillable and therefore tight.
\end{thm}

\begin{figure}[ht!]
\centerline{
\psfrag {g0}{$\gamma_0$}
\psfrag {g1}{$\gamma_1$}
\psfrag {g2}{$\gamma_2$}
\psfrag {g3}{$\gamma_3$}
\psfrag {r0}{\smash{\rlap{\kern 0pt \raise 6pt\hbox{$r_0$}}}}
\psfrag {r1}{\smash{\rlap{\kern 0pt \raise 6pt\hbox{$r_1$}}}}
\psfrag {r2}{\smash{\rlap{\kern 0pt \raise 6pt\hbox{$r_2$}}}}
\psfrag {r3}{\smash{\rlap{\kern 0pt \raise 6pt\hbox{$r_3$}}}}
\epsfysize=1in\epsfbox{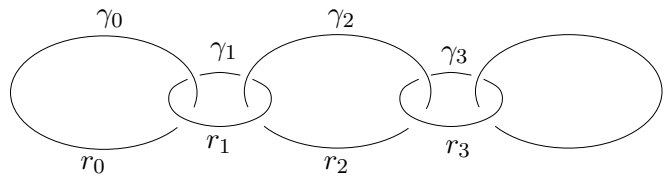}}
\caption{Surgery along link}
\label{link}
\end{figure}

Observe that for the lens space $L(p,q)$, $p>q>0$, and
the continued fraction expansion $(r_0,r_1,..,r_k)$ for $-{p\over q}$, we have
a linked chain of unknots in $S^3$ with framings $r_0$, $r_1$, ..., $r_k$ (in order along
the chain), along which we can do Legendrian surgery to obtain $L(p,q)$.
Denote the unknots by $\gamma_0,\cdots, \gamma_k$. See Figure \ref{link}.
To perform Legendrian surgery, $\gamma_i$ must have Thurston--Bennequin invariant
$tb(\gamma_i)=r_i+1$. There however are $|r_i+1|$ choices for the rotation number
$r(\gamma_i)$: $r_i+2, r_i+4,\cdots, r_i+2|r_i+1|$.

\begin{proof}
[Proof of Proposition \ref{lowerbound}.]  We will take an easy way out by using the
following theorem, due to Lisca and Mati\'c \cite{LM}:
\begin{thm} [Lisca--Mati\'c]
Let $X$ be a smooth 4--manifold with boundary.  Suppose $J_1$, $J_2$ are two Stein structures
with boundary on $X$.  If the induced contact structures $\xi_1$, $\xi_2$ on $\bdry X$ are
isotopic, then $c_1(J_1)=c_1(J_2)$.
\end{thm}

Let $X$ be the Stein surface obtained from $B^4$ by attaching 2--handles $H_1,\cdots,$ $H_k$
corresponding to Legendrian surgeries with coefficients $r_1,\cdots,r_k$ along the link in
Figure \ref{link}.  If $c_1(X)$ is the canonical class and $h_i$ is a 2--dimensional class
supported on $H_i$, then $\langle c_1(X),h_i\rangle= r(\gamma_i)$.  For the various
$r(\gamma_i)$, the $c_1(X)$ are distinct.
\end{proof}

{\bf Remark}\qua  Theorems \ref{2} and \ref{3} can be thought of as
a generalization of Eliashberg and Fraser's classification of Legendrian unknots \cite{F}.

\subsection{Homotopy classification}

In this section we will distinguish the minimally twisting tight structures on
$T^2\times I$ and tight structures on $S^1\times D^2$ using the relative Euler class. Observe
that the proof of Part 2(a) of Theorem \ref{2} implies the following lemma:

\begin{lemma}
Let $(T^2\times I,\xi)$ be a contact manifold which admits a factorization
$T^2\times I=\cup_{i=0}^{k-1} N_i$, where each $N_i=T^2\times
[{i\over k},{i+1\over k}]$ is a {\it basic slice}, and $s_0=-1>s_{1\over k}>s_{2\over k}>\cdots
> s_1=-{p\over q}$, $p>q>0$ integers, is obtained by taking the shortest counterclockwise
sequence from $s_1$ to $s_0$ on $\bdry \H^2$ as in Lemma \ref{lemma:hops}.
Then $\xi$ is tight and minimally twisting.  Moreover, such a factorization is unique up to a
shuffling within a continued fraction block.
\end{lemma}

\proof The fact that $\xi$ is tight and minimally twisting follows from observing that
Equation \ref{eqn6} is actually an equality.  This means that every gluing of basic layers is
tight, provided the slopes $s_{0 \over k}, s_{1\over k},\cdots, s_{k\over k}$ are obtained by
taking the shortest counterclockwise sequence from $s_1$ to $s_0$ on $\bdry \H^2$, since
the number of contact structures obtained this way is at most the right-hand side of Equation
\ref{eqn6}.    If the factorization was not unique up to a shuffling within a continued
fraction block, the number of potential tight contact structures will be less than the actual
number of tight contact structures, a contradiction. \qed

\subsubsection{Minimally twisting $T^2\times I$}

\begin{prop} \label{euler1} The minimally twisting tight contact structures on $T^2\times I$
with $\#\Gamma_{T_i}=2$ and fixed $s_0$, $s_1$ can be distinguished by the relative Euler
class. \end{prop}

\proof  For convenience, set $s_0=-1$, $s_1=-{p\over q}$, $p>q>0$ integers.  Consider the
factorization $T^2\times I=\cup_{i=0}^{k-1} N_i$, where each
$N_i=T^2\times [{i\over k},{i+1\over k}]$ is a {\it basic slice}, and $s_0>s_{1\over
k}>s_{2\over k}>\cdots > s_1$, is obtained by taking the shortest counterclockwise sequence
from $s_1$ to $s_0$ on $\bdry \H^2$.   Let $v_i$ the shortest integral vector with slope
$s_{i\over k}$ and negative $x$--coordinate.  Then $PD(e(\xi_{N_i},s))=\pm(v_{i+1}-v_i)$,
and
\begin{equation}  \label{eqn:euler}
PD(e(\xi,s))=\sum_{i=0}^{k-1} \pm (v_{i+1}-v_i),
\end{equation}
for $(T^2\times I,\xi)$.

Let $A$ be a horizontal convex annulus with Legendrian boundary, after a perturbation of
$T_i$.  We claim that $\langle e(\xi,s),A\rangle$ are distinct for the different $\xi$.
Let $(r_0,r_1,\cdots, r_k)$ be the continued fraction representation of
$-{p\over q}$.  We will track the change in $\langle e(\xi|_{T^2\times[0,i+1]},s),A_i\rangle$,
where $A_i$ is the horizontal convex annulus for $N_0\cup\cdots \cup N_i$, starting from the
innermost layer with boundary slope $-1$, and moving out to $-{p\over q}$.  Consider the
boundary slope $s_i=-{p_i\over q_i}={-ar_j+b\over cr_j-d}$, corresponding to the continued
fraction representation $(r_0,\cdots,r_j)$, where $(-c,a)$ and $(-d,b)$ form an oriented basis
and $a >c\geq 0$, $b\geq d>0$. Inductively we have $|\langle
e(\xi|_{T^2\times[0,i+1]},s),A_i\rangle|< p_i$. Then $(r_0,\cdots,r_j-1,-2)$ corresponds to
$s_{i+1}={(-ar_j+b)+(-a(r_j-1)+b)\over (cr_j-d)+(c(r_j-1)-d)}$, and
\begin{eqnarray}
|\langle e(\xi|_{T^2\times[0,i+2]},s),A_{i+1}\rangle-
\langle e(\xi|_{T^2\times[0,i+1]},s),A_{i}\rangle|\kern-1in&&\nonumber\\
&=&(-a(r_j-1)+b)\\
& \geq&  -ar_j+b\\
& >&  |\langle e(\xi|_{T^2\times[0,i+1]},s),A_{i}\rangle|
\end{eqnarray}
We find that $\langle e(\xi,s),A\rangle$ determines the tight contact
structure. \qed

\subsubsection{Solid tori}
Let us now give a homotopy classification of the potential tight structures on $S^1\times D^2$
with $T=\bdry (S^1\times D^2)$, $\#\Gamma_{T}=2$, and boundary slope $-{p\over q}$.

\begin{prop}
The elements $[\xi]$ of $\pi_0(\mbox{Tight}(S^1\times D^2, \Gamma))$, $\#\Gamma=2$, $s=-{p\over
q}$ are distinguished by $r(\bdry D)=\langle e(\xi, s),D\rangle=\#(\mbox{Components of $R_+$})
-\#(\mbox{Components of $R_-$})$, where $D$ is a convex meridional disk with Legendrian
boundary.     Here $r$ denotes the rotation number.
\end{prop}

\proof Follows from Proposition \ref{euler1} and noting that every connected component of
$D\backslash \Gamma_D$ has Euler characterstic $1$.
\qed

\subsubsection{Lens spaces}
\begin{prop}   The homotopy classes of the tight contact structures on $L(p,q)$ are all
distinct.
\end{prop}

\proof
Let us use the same notation as before.  In particular,
$V_0$ is  the standard neighborhood of the Legendrian core curve $C_0$ with the largest
twisting number, and $V_1=L(p,q)\backslash V_0$.   Every tight contact structure is obtained by
Legendrian surgery along $\gamma_i$, $i=1,...,k$,  in Figure \ref{link}.
Let $V'_i$ be small standard neighborhoods of $\gamma_i\subset
S^3$, with boundary slopes ${1\over r_i+1}$ (use the standard framing on $S^3$).  Also let
$V''_i$ be standard neighborhoods of Legendrian curves with twisting number $-1$.
We remove $V'_i$ from $S^3$, and glue in $V''_i$ by mapping $\bdry
V''_i\rightarrow -\bdry (S^3\backslash V'_i)$ via
$\left(\begin{array}{cc}
-r_i & 1 \\ -1 & 0
\end{array}
\right)$.
For $\bdry V''_i$, $(1,0)^T$ is the meridian of $V''_i$ and $(0,1)^T$ the direction of the
Legendrian core curve with twist number $-1$. For $-\bdry(S^3\backslash V'_i)=\bdry V'_i$,
$(1,0)^T$ is the meridian of $V'_i$ and $(0,1)^T$ the longitude for
the framing for $V'_i$.  We now identify $V''_0\simeq V_0$ via a Dehn twist
to match up the framings (the Legendrian core curve of minimized twisting number 0 for $V_0$
must go to the Legendrian core curve of twisting number $-1$ for $V''_0$).
This then gives rise to a map $\bdry V_0\rightarrow \bdry V_1=
\bdry (S^3\backslash V'_0)$
$$\left(\begin{array}{cc} 0 & 1 \\ 1 &0 \end{array} \right)
\left(\begin{array}{cc} -r_0 & 1 \\ -1 & 0 \end{array}\right)
\left(\begin{array}{cc} 1 & -1 \\ 0 & 1\end{array}\right)
=\left(\begin{array}{cc} -1 & 1  \\ -r_0 & r_0+1 \end{array}\right).
$$
Here, $(1,0)^T$ and $(0,1)^T$ are the same as before for $\bdry V_0$,
and for $\bdry V_1$, $(0,1)^T$ is the meridional direction for $\bdry V'_0$
and $(1,0)^T$ is the meridional direction of $V_1$.

Now consider $T^2\times I=S^3\backslash(V'_0\cup V'_1)$ with boundary slopes $r_0+1$ and
${1\over r_1+1}$.
There exist $|r_0+1|$ possibilities for $\langle e(\xi,s),A\rangle$
on a horizontal annulus $A=S^1\times\{0\}\times I$ with Legendrian
boundary, depending on the rotation number of $\gamma_0$.
On the other hand, we have $|r_1+1|$ possibilities for a vertical
annulus $B$, depending on the rotation number of $\gamma_1$.  Therefore, we find
that all $|(r_0+1)(r_1+1)|$ possible tight structures on
$T^2\times I$ with the given boundary slopes are realized.

Next, we transform $T^2\times I$ via $\left( \begin{array}{cc}
0 & -1 \\ 1 & -r_1 \end{array} \right)$ to get boundary slopes
${1-r_1(r_0+1)\over
-(r_0+1)}$ and $-1$.  Notice that $\gamma_2$ is now vertical, with
boundary slope ${1\over r_2+1}$.  Consider a horizontal annulus $A$ with Legendrian boundary
for this (transformed) $T^2\times I$.  It will cut through $V'_2$, and $\langle
e(\xi,s),A\rangle$  will uniquely determine the homotopy class of the tight structure by
Proposition \ref{euler1}. Now take $N=S^3\backslash(V'_0\cup V'_1\cup V'_2) \cup V''_1$, ie,
we fill in $V''_1$ and remove $V'_2$.  Consider the new
horizontal annulus $A'$, obtained by removing the meridional disk of $V_2'$ and adding in the
meridional disk of $V_1''$.  Then $\langle e(\xi,s),A\rangle=\langle e(\xi,s),A\rangle$, where
the relative Euler class is taken in the respective manifolds.
Now, $\langle e(\xi,s), B\rangle$ for  the vertical annulus $B$ with Legendrian boundary
spanning from $\bdry V''_1$ to $\bdry V'_2$ corresponds to the rotation number of $\gamma_2$.
Therefore we see that all $|(r_0+1)(r_1+1)(r_2+1)|$ possible tight structures are represented
on $N$.

Take $V_1$ with convex boundary and horizontal Legendrian rulings, and perturb the
characteristic foliation into a nonsingular Morse--Smale characteristic foliation; also take a
meridional disk $D$ for $V_1$ with Legendrian boundary and perturb into $D'$ with transverse
boundary. Etnyre in \cite{Et} relates the number of positive elliptic points on $D$ (or the
self-linking number $sl(\bdry D')$) to the homotopy classes of 2--plane fields on $L(p,q)$, and
shows, in particular, that the homotopy classes of tight structures on $L(p,q)$ are distinct if
the self-linking numbers are distinct. Our proposition follows from observing that $r(\bdry D)$
are distinct for the contact structures with distinct $r(\gamma_i)$ and using
$sl(\bdry D')=tb(\bdry D) \pm r(\bdry D)$.\qed

\subsubsection{Gluing}
As a consequence of the classification of minimally twisting tight contact structures on
$T^2\times I$ we have the following gluing theorem:

\begin{thm}[Gluing $T^2\times I$]  \label{gluing} Let $\xi$ be a contact structure on
$T^2\times [0,n]$, where each $N_i=T^2\times [i,i+1]$ is a basic slice.  Assume all $s_i$ lie
on the counterclockwise arc $[s_n,s_0]\subset \bdry \H^2$, and $s_n<s_{n-1}<s_{n-2}<\cdots
<s_0$. Here we write $a<b$ if $b$ is closer to $s_0$ than $a$ is on the arc $[s_n,s_0]$.
Then $\xi$ is tight if and only if one of the following holds:
\be
\item $s_n,s_{n-1},\cdots,s_0$  is the shortest sequence from $s_n$ to $s_0$.
\item $s_n,\cdots,s_0$ is not the shortest sequence and there is a triple $s_{i+1},s_i,s_{i-1}$
where $s_i$ is removable from the sequence, ie, there exists an edge from $s_{i+1}$ to
$s_{i-1}$ along $\bdry \H^2$. $T^2\times[i-1,i+1]$ is a basic slice and we shorten the sequence
by omitting $s_i$. By repeating this procedure we get to Case (1). \ee
\end{thm}

In Theorem \ref{gluing}, we may need to determine when $T^2\times[i-1,i+1]$ is a basic slice,
given that $N_{i-1}=T^2\times [i-1,i]$ and $N_i=T^2\times [i,i+1]$ are basic slices.  The
relative Euler class is useful for this.     Let $v_i$ be a shortest integral vector for $s_i$,
chosen consistently so that there exists an element of $SL(2,\Z)$ which maps
$v_{i+1},v_i,v_{i-1}$ to $(1,0),(1,1),(0,1)$.  For the relative Euler classes of the two
component basic slices to add up to a relative Euler class for basic slice we need
$PD(e(\xi|_{N_{i-1}},s))=v_i-v_{i-1}$ and $PD(e(\xi|_{N_i},s))=v_{i+1}-v_i$, or both signs
reversed.

\section{Tight contact structures on $T^2\times I$}
\label{whatever}

\subsection{Universal tightness}

In this section we will precisely determine which minimally twisting
tight structures on $T^2\times I$, $S^2\times D^2$,
and $L(p,q)$ are universally tight. Let $\Sigma$ be an annulus with a collared Legendrian
boundary and negative
twisting number on both boundary components.  If $\Gamma$ is the dividing set,
then denote the connected components of $\Sigma\backslash \Gamma$ by $\Sigma_i$.  We
call $\Sigma_i$ a {\it one-sided component}
if it intersects only one boundary component of $\Sigma$.
$\Sigma_i$ is {\it boundary-parallel} if it is a half-disk which intersects a single
dividing curve $\gamma$ and $\gamma$ is boundary-parallel.

Recall
$\Sigma_i$ is positive if the oriented flow exits from $\bdry\Sigma$.

\begin{prop} \label{universal}
\be
\item There are exactly two tight constact structures on $M=T^2\times I$ with minimal twisting, $\#\Gamma_i=2$,
and boundary slopes $s_1=-{p\over q}, s_0=0$  ($p>q>0$ positive integers) which are
universally tight.
They satisfy $PD(e(\xi,s))=\pm( (-q,p)-(-1,0))$.
\item There are exactly two tight contact structures on $M=S^1\times D^2$ with $\#\Gamma_{\bdry M}=2$,
and boundary slope $s=-{p\over q}<-1$ which are universally tight.  (If $s=-1$ there is exactly one.)
\item There are exactly two tight contact structures on $M=L(p,q)$ with $q\not=p-1$  which are universally
tight.  (If $p=p-1$ there is exactly one.)
\ee
\end{prop}

The two universally tight structures on $T^2\times I$ are diffeomorphic via $-id$,
where $id$ is the identity map on $T^2$.

\medskip



\begin{proof}
(1) Let $A=S^1\times \{0\}\times I$. Consider its one-sided components $A_i$ (they are all along
$S^1\times \{0\}\times\{1\}$).
If $PD(e(\xi,s))\not=\pm( (-q,p)-(-1,0))$, not all the one-sided components have the same sign.
We have two possibilities:
(A) there exists a positive one-sided $A_1$ and negative one-sided $A_2,\cdots, A_k$ which lie further
toward $S^1\times\{0\}\times\{0\}$ as in the left-hand side of Figure \ref{one-sided} (or signs reversed), or
(B) there is a positive boundary-parallel $A_1$ as well as a negative  boundary-parallel $A_2$.
Let $\gamma_i$ be the dividing curve on $A_i$ which is `farthest' from $S^1\times\{0\}\times \{1\}$ (ie,
the half-disk cut off by $\gamma_i$ contains the other dividing curves which bound $A_i$).
\begin{figure}[ht!]
\centerline{
\psfrag {g1}{$\gamma_1$}
\psfrag {g2}{\smash{\rlap{\kern -2pt \raise 0pt\hbox{$\gamma_2$}}}}
\psfrag {g3}{$\gamma_3$}
\psfrag {A1}{$A_1$}
\psfrag {A2}{\smash{\rlap{\kern -2pt \raise 0pt\hbox{$A_2$}}}}
\psfrag {A3}{\smash{\rlap{\kern -2pt \raise 0pt\hbox{$A_3$}}}}
\epsfysize=1.8in\epsfbox{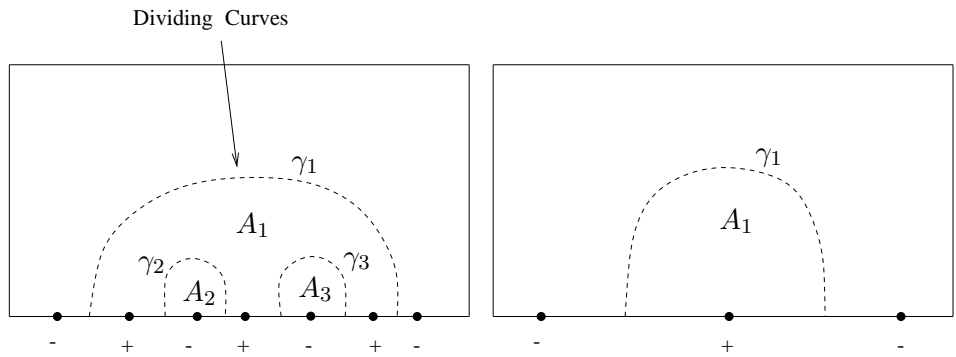}}
\caption{One-sided components}
\label{one-sided}
\end{figure}

For both cases, pass to the cover $\widetilde{M}=S^1\times \R\times I$.
Let us first consider Case (B). There exist lifts
$\widetilde A=S^1\times \{0\}\times I$ and $\widetilde A'=S^1\times
\{m\}\times I$, $m\in \Z^+$, for which $N_1=\widetilde A\cup \widetilde A'\cup
(S^1\times [0,m]\times \{1\})$, after rounding the edges, has a dividing curve $\gamma$ which bounds
a disk.  In fact, a lift of $\gamma_2$ on $\widetilde A$ will connect up to
a lift of $\gamma_1$ on $\widetilde A'$ for suitably chosen $m$.
The existence of a null-homotopic dividing curve then implies that $\widetilde{M}$ is overtwisted.

For Case (A), take $\widetilde A$, $\widetilde A'$ as above, as well as lifts
$\tilde{\gamma}_i$ on $\widetilde A$ and $\tilde\gamma_i'$ on $\widetilde A'$.
Pick $m\in \Z^+$ so that $\tilde{\gamma}_2$ connects the left endpoint of $\tilde\gamma_1'$
to the left endpoint of $\tilde\gamma_2'$, $\tilde{\gamma}_3$ connects the
right endpoint of $\tilde\gamma_2'$ to the left endpoint of
$\tilde\gamma_3'$, and so on.
What we still lack is a dividing curve connecting the right endpoint of
$\tilde\gamma_k'$
to the right endpoint of $\tilde\gamma_1'$.  Take $\widetilde A''=
S^1\times \{m'\}\times I$, $m'\in\Z^-$, as well as $N_2=\widetilde A'\cup
\widetilde A''\cup (S^1\times[m',m]\times I)$, after rounding the edges.
If we pick $m'$ appropriately, we can make the desired connection along
$N_2$.  Now, the dividing curve sits on the branched surface $N_1\cup N_2$, and there exists
an overtwisted disk on this branched surface.

If the tight contact structure on $M$ has a horizontal convex annulus $A$, all of whose
one-side components are boundary-parallel with the same sign, then $M$ can be embedded into,
and is universally tight because $(T^3,\xi_1)$ is.

(2) and (3) are left for the reader.
\end{proof}

{\bf Note}\qua Any tight contact structure $\xi$ on $M=T^2\times I$ with minimal
twisting, $\#\Gamma_i=2$, and boundary slopes $-{p\over q}$ and $-1$
factors into continued fraction blocks of the form $N=(\R^2/\Z^2) \times I$ with boundary slopes
$s_1=-m$, $s_0=-1$, $m\in \Z^+$.  Consider the block $N$.
According to Shuffling Lemma, we can arrange
the dividing curves on a horizontal annulus $A$ so we have the following:
(1) two dividing curves $\gamma_1$, $\gamma_2$ which go across,
(2) the rest are boundary-parallel curves.
If there exist both positive and negative half-disk cut off by the boundary-parallel curves,
then the double cover $\widetilde{N}=(\R/\Z)\times (\R/2\Z)\times I$ will
be overtwisted, using the methods of Proposition \ref{universal} and the special form of the
dividing curves on $A$.  Hence, if any of the blocks of $M$ have mixed
signs, then a double cover of $M$ is overtwisted.

\subsection{Non-minimal twisting for $T^2\times I$}

We will now finish the proof of  Parts 2(b) and 3 of Theorem \ref{2}.
Consider a basic slice $(N_0=T^2\times I,\overline\xi)$ with boundary slopes $s_1=0$,
$s_0=\infty$.   If we fix a boundary characteristic foliation compatible with $\Gamma_i$, there are 2 possible tight
structures on $N_0$.   Let $\overline\xi$ be the tight structure on $N_0$ for which
$PD(e(\xi,s))=(1,-1)\in H_1(T^2;\Z)$.

Let $N_{n\pi\over 2}$ be $N_0$ rotated counterclockwise by ${n\pi\over 2}$, $n\in \Z$.
Take $\xi^+_1= N_0 \cup N_{\pi\over 2}$, $\xi^+_2=N_0\cup N_{\pi\over 2}
\cup N_{\pi}\cup N_{3\pi\over 2}$, ... ,
$\xi^-_1= N_\pi\cup N_{3\pi\over 2}$, $\xi^-_2=N_\pi\cup N_{3\pi\over 2}\cup
N_{2\pi}\cup N_{5\pi\over 2}$, ... ,
where the $T_0$ of $N_{n\pi\over 2}$ is identified with the $T_1$ of
$N_{(n+1)\pi\over 2}$. These $\xi_n^\pm$ can be embedded inside some
$(T^3,\xi_m)$, $m\in\Z^+$, and
are therefore universally tight.

\begin{lemma} \label{parallelends}
A tight $(M=T^2\times I,\xi)$ with $\#\Gamma_{T_i}=2$, $i=0,1$, non-minimal twisting,
and $s_1=s_0=0$ is isotopic to one of the
$\xi_n^\pm$, $n\in \Z^+$.
\end{lemma}

\begin{figure}[ht!]
\centerline{
\psfrag {T0}{$T_0$}
\psfrag {T1}{$T_1$}
\epsfysize=2.5in\epsfbox{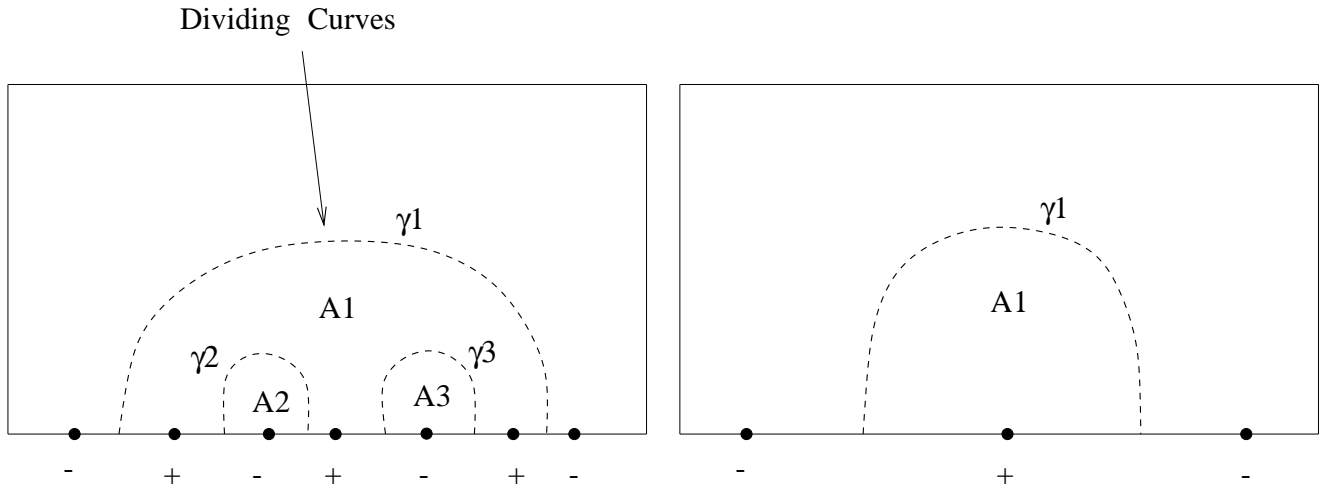}}
\caption{Configurations of dividing curves on $B$}
\label{B}
\end{figure}

\proof Let $\xi$ be a tight structure on $T^2\times I$ with minimal
boundary and $s_1=s_0=0$.  Assume $r_1=r_0=\infty$.  Let $B=\{0\}\times S^1\times I$
be a vertical convex annulus with Legendrian boundary and oriented normal
${\partial \over \partial x}$.  Also assume that $\#\Gamma_B$ is minimal among all
vertical convex annuli in its isotopy class rel boundary. See Figure \ref{B} for
possible configurations of dividing curves on $B$.
If $\Gamma_B$ does not have
any boundary-parallel dividing curves, then $\#\Gamma_B=2$ and the two dividing curves will
go across from $T_0$ to $T_1$; rounding the edges, we find that we are in the minimally twisting,
nonrotative case.  Therefore $\Gamma_B$ must have boundary-parallel dividing curves.
We then cut along $B$ and perform edge-rounding to obtain a solid
torus $S^1 \times D^2$ with $2 + 2i$ vertical dividing curves, where $i$ is the number of
closed dividing curves (parallel to the boundary) on $B$.

Next cut $S^1 \times D^2$ along a meridional disk $D$ after modifying the boundary to be
standard with horizontal rulings.  The configuration of dividing curves on $D$ is completely
determined by the condition that the number of dividing curves on $B$ be
minimal.  Let $\gamma_0$ and $\gamma_1$ be the dividing curves on $\bdry(S^1 \times D^2)$ which
intersect $T_0$ and $T_1$ (ie,  $\Gamma_{T_0}, \Gamma_{T_1}$ become part of $\gamma_0$, $\gamma_1$
after edge-rounding). Then all $\gamma\in\Gamma_D$ must separate $D \cap \gamma_1$ from $D
\cap \gamma_0$ (hence the dividing curves of $D$ are parallel segments, with only two
boundary-parallel components, each containing one $D\cap \gamma_i$ as the half-elliptic point
on the interior); otherwise there would exist a bypass which allows for a reduction in the number of
dividing curves on $B$.

Therefore, the tight structure $\xi$ on $M$ depends only on $\Gamma_B$, which in turn is determined by
the sign of the boundary-parallel component of $B$ along $T_1$, together with $i+2=\#\Gamma_B$.  If the sign is
$+$ ($-$), then $\xi=\xi_{i+1}^+$ ($\xi_{i+1}^-$).
\qed

\begin{lemma} \label{cc}  The $\xi_n^\pm$, $n\in \Z^+$ are distinct.
\end{lemma}

\proof  We distinguish among the four classes
$\xi_{2m-1}^+$, $\xi_{2m-1}^-$, $\xi_{2m}^+$, $\xi_{2m}^-$,
$m\in \Z$, according to whether attaching $N_{-{\pi\over 2}}$ to the front preserves
tightness (they do for $\xi^+_n$) and whether attaching $N_0$ to the back preserves tightness (they do
for $\xi^+_{2m}$ and $\xi^-_{2m-1}$).  In the cases when tightness is not preserved, we can
find horizontal annuli with a dividing curve bounding a disk.

In each case, $m$ determines the twisting.  For example, consider $\xi^+_{2m}$.  If
we glue the front and back via the identity map, we obtain the tight contact
structure $(T^3,\xi_m)$ described previously, and the $\xi_m$ are distinguished
by the minimal twisting number for closed curves isotopic to $S^1=I/\sim$.
(This is due to Kanda \cite{K}.)\qed

\begin{prop}     \label{hey}
A tight  $(M=T^2\times I,\xi)$ with $\#\Gamma_{T_i}=2$, $s_0=0$, $s_1=-{p\over q}$,
$p>q>0$, and non-minimal twisting is universally tight.  Moreover, there exists a splitting
$T^2\times I=(T^2\times [0,{2\over 3}])\cup (T^2\times [{2\over 3},1])$ where $T_{2\over 3}$ is
convex with $\#\Gamma_{2\over 3}=2$, $T^2\times [{2\over 3},1]$ is minimally twisting,
and $T^2\times [0,{2\over 3}]$ is isotopic to some $\xi^\pm_n$.
\end{prop}

\proof Given such $(T^2\times I,\xi)$, there exist enough bypasses to factor $M$ into
$M_1=T^2\times [0,{1\over 3}]$, $M_2=T^2\times [{1\over 3},{2\over 3}]$,
and $M_3= T^3\times [{2\over 3},1]$, where $s_0=0$, $s_{1\over 3}=-{p\over q}$, $s_{2\over 3}=0$,
$s_1=-{p\over q}$,  and $M_1$, $M_3$ are minimally twisting.
Notice that the tight structure on $M_2\cup M_3$ is one of the $\xi_n^\pm$ as in Lemma
\ref{cc}, and is universally tight.
By Proposition \ref{universal}, this reduces the possibilities on $M_3$  to two.
A consideration of the signs will reveal that $\xi$  on $M$ is universally tight, and
can be split into $M_3$ with minimal twisting, and $M_1\cup M_2$ with some $\xi^\pm_n$. There
will be two such, according to whether the horizontal bypasses on $M_1$ are all positive or all negative.
This $n$ is unique --- this is proved in the same way as Lemma \ref{cc}.  \qed

\begin{prop}     \label{hiho}
The $I$--twisting $\beta_I$ of a tight contact structure $\xi$ on $T^2\times I$ is well-defined and finite.
In particular,  $\beta_I$ is independent of the factorization $T^2\times I=\cup_{k=0}^{l-1}(T^2\times
[{k\over l},{k+1\over l}])$ into minimally twisting slices.
\end{prop}

\proof  The finiteness follows (provided $\beta_I$ is well-defined)
from the factorization in Proposition \ref{hey} into a minimally twisting $T^2\times I$ and
a slice with some $\xi^\pm_n$, which can be factored into basic slices as above.
It remains to prove that $\beta_I$ remains invariant under subdivisions.  But this follows from observing
that
$\beta_I$ is well-defined on a minimally twisting $T^2\times I$, since any factorization will satisfy
$s_1<s_{l-1\over l}<s_{l-2\over l}<\cdots <s_0$ where $a<b$ if there exists a counterclockwise
subarc of $[s_1,s_0]\subset \bdry \H^2$ from $a$ to $b$.\qed

\subsection{Non-minimal boundary}

\subsubsection{Model for increasing the torus division number}\label{section:increase}
\begin{figure}[ht!]
\centerline{
\psfrag {T0}{\smash{\rlap{\kern -2pt \raise 0pt\hbox{$T_0$}}}}
\psfrag {T'}{$T'$}
\psfrag {A}{$A$}
\epsfysize=2in\epsfbox{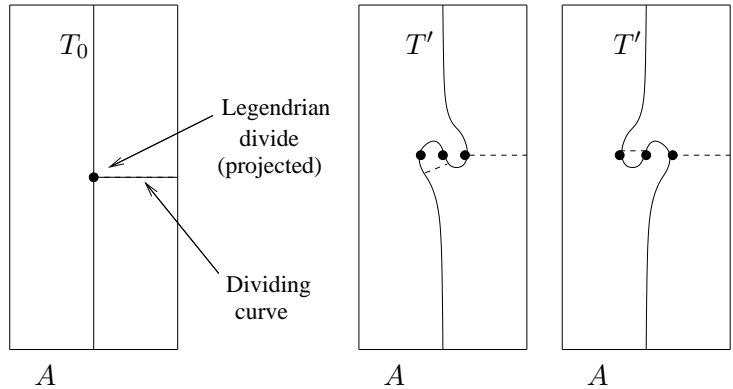}}
\caption{Perturbation to increase the division number}
\label{increase}
\end{figure}

Let $T^2$ be a convex
torus in standard form with $s=\infty$, $r=0$, and $\#\Gamma=2n$.  Since $T^2$ is
convex, there is a universally tight, $I$--invariant neighborhood $T^2\times
[-\varepsilon,\varepsilon]$ of $T^2= T_0$.
The horizontal annulus $A=S^1\times \{0\}\times [-\varepsilon,\varepsilon]$
has parallel dividing curves from $T_{-\varepsilon}$ to $T_\varepsilon$.  We will find
$T'$ $C^0$--close to $T_0$ so that the division number of $T'$ is $n+1$.  Modify $T_0$ near one
of its Legendrian divides to increase $\#\Gamma$ by $2$, as in Figure \ref{increase}.  Here,
$T_0$, $T'$ are invariant in the $y$--direction,
and their projections to $A$ are as shown.
One of the modifications will increase $\# R_+$ by $1$, and the other will
increase $\# R_-$ by $1$.  Now perturb $T'$ so it is standard.  The region bounded by $T'$
and $T_{-\varepsilon}$ will be universally tight.  Note that we can insert a bypass to create
any possible configuration for $T^2\times I$ with no twisting, $n_1=n+1$, $n_0=n$,
$s_1=s_0=\infty$. Here $n_i$ is the torus division number for $T_i$.

By iterating this procedure, we find that any $(N=T^2\times I,\xi)$ with
$n_1\geq  n_0$, $s_1
=s_0$, and bypasses on a horizontal annulus only along $T_1$, can be
obtained as a universally
tight structure inside a translation invariant one on $T^2\times I$.
Moreover, for any $(M,\xi)$ tight and $\bdry M$ a union of tori in standard
form, attaching
layers of the same type as $N$ is an operation which preserves tightness,
since the
resulting manifold and contact structure can be found inside $(M,\xi)$ due
to convexity.

\subsubsection{Template matching}   In the next section we shall reduce the problem of
classifying nonrotative contact structures to a 2--dimensional problem which we  treat first.
Given an oriented compact surface $\Sigma$ with boundary, and a finite subset $\sigma\subset
\bdry \Sigma$ which we call the {\it markings}, define
$\mathcal{C}(\Sigma)$ to be the set of configurations  $\Gamma$, where each configuration is a
set of arcs with endpoints on $\sigma$, and every point of $\sigma$ is used as an endpoint exactly once.  (In
particular, $|\sigma|$ must be even.)  If $\Sigma$ and $\Sigma'$ are glued along a boundary component $C$,
then there exists a natural map:
$$G\co \mathcal{C}(\Sigma)\times \mathcal{C}(\Sigma')\rightarrow \mathcal{C}(\Sigma\cup_C\Sigma'),$$
$$(\Gamma, \Gamma')\mapsto (\Gamma\cup\Gamma')_0,$$
where $()_0$ means throw away any closed curves.
If $\Sigma$ is an annulus, and if the number of
markings on the two boundary components are $m\leq n$,  then we define
$\mathcal{C}_0(\Sigma)\subset \mathcal{C}(\Sigma)$ to be the set of configurations where $m$ of
the markings on on one boundary component are connected to $m$ markings on the other boundary component
via an arc.

Consider annuli $A=S^1\times[0,1]$ and $B=S^1\times [1,2]$.  Fix markings
$p_1,\cdots, p_{2m}$ (in cyclical order) each on $S^1\times\{0\}$ and $S^1\times\{2\}$,  and markings $q_1,\cdots,
q_{2n}$ (in cyclical order) on $S^1\times \{1\}$, where $m<n$. If $\mathcal{C}\subset \mathcal{C}(A)$,
then we define the {\it dual of $\mathcal{C}$} to be $\mathcal{C}^*=\{\Gamma\in \mathcal{C}(B)|
\Gamma'\cup \Gamma =id, \forall \Gamma'\in \mathcal{C} \}$, where $id\in \mathcal{C}(A\cup B)$ is the unique element (up to changes in
holonomy).

\begin{lemma} [Reflexive Property] \label{reflexive}
Consider $\mathcal{A}=\{\Gamma_0\}\subset \mathcal{C}_0(A)$, for any element $\Gamma_0$.
Then $(\mathcal{A}^*)^*=\mathcal{A}$.
\end{lemma}

What this lemma says is that $\Gamma_0$ can be detected externally by considering the space of
{\it templates} on $B$ which give $2m$ parallel curves (and no closed homotopically trivial curves)
when glued to $\Gamma_0$.

\proof  By induction on $n-m$.  Assume first $n-m=1$.  Then $\Gamma_0$ will consist of $2m$ curves
which cross from $S^1\times \{0\}$ to $S^1\times \{1\}$, and one boundary-parallel curve
from $q_k$ to $q_{k+1}$.  $\mathcal{A}^*$ will have two configurations, $\Gamma_0^\pm$,
both with $2m$ curves crossing
from $S^1\times\{1\}$ to $S^1\times \{2\}$.  $\Gamma_0^+$ ($\Gamma_0^-$)
will have a boundary-parallel curve from $q_{k+1}$ to $q_{k+2}$ (resp. $q_{k-1}$ to $q_{k}$).
$(\mathcal{A}^*)^*=\{\Gamma_0^+\}^*\cap \{\Gamma_0^-\}^*=\mathcal{A}$.

Suppose the lemma is true for all $\Gamma_0$ with $n-m=l$.  Now assume $\Gamma_0$ has
$n-m=l+1$.  We claim any $\Gamma\in(\mathcal{A}^*)^*$ will have a factorization
$A=(S^1\times [0,{1\over 2}]) \cup (S^1\times [{1\over 2},1])$, $\Gamma=\Gamma_{l-1}\cup
\Gamma_l$, where $\Gamma_l$ consists of $2(n-1)$ curves which go across and $1$ boundary-parallel
curve, and $\Gamma_{l-1}$ consists of $2m$ curves which go across.  Moreover, the boundary-parallel
curve on $\Gamma_l$ will coincide with a boundary-parallel curve on $\Gamma_0$.
This can be seen as follows:  Let $q^*_i$ be the point on $S^1\times \{1\}$ which is connected to
$p_i$ by an arc of $\Gamma_0$.  Look at two consecutive $q^*_i$, $q^*_{i+1}$.  If $q^*_{i+1}-q^*_i>1$,
then there exists a boundary-parallel arc $\delta$ of $\Gamma_0$ inbetween.  Every boundary-parallel
arc with endpoints on the interval $[q^*_i, q^*_{i+1}]$ can be incorporated into an element of
$\mathcal{A}^*$, except when the endpoints are exactly the endpoints of $\delta$.  If
$q^*_{i+1}-q^*_i=1$, and $m>1$, then we take parallel curves from $q^*_i$ and $q^*_{i+1}$ to $S^1\times \{2\}$,
and extend to an element of $\mathcal{A}^*$.    Now consider $\Gamma\in (\mathcal{A}^*)^*$.  The discussion
above restricts the possible positions of the boundary-parallel arcs of $\Gamma$.  If there is a boundary
parallel arc $\delta$ which is not a boundary-parallel arc for $\Gamma_0$ as well, then take
$q^*_i$, $q^*_{i+1}$ so that $\bdry \delta\subset [q^*_i,q^*_{i+1}]$.  If $q^*_{i+1}-q^*_i>1$, then there
is only one position where a closed homotopically trivial curve is not created by summing with
some $\Gamma'\in \mathcal{A}^*$.  If  $q^*_{i+1}-q^*_i=1$, and $m>1$,
then there exists $\Gamma'\in\mathcal{A}^*$ such that summing creates a boundary-parallel arc along
$S^1\times\{2\}$.  If $q^*_{i+1}-q^*_i=1$ and $m=1$, then the only $\Gamma$ which is not immediately
factorable is one where $\delta\subset \Gamma$ has endpoints $q^*_i$, $q^*_{i+1}$, and there are no
other boundary-parallel arcs (hence all the other arcs with endpoints on $S^1\times \{1\}$ are
concentric arcs).  This can be eliminated by taking $\Gamma'$ which extends
the union of two arcs $\delta_1$ (with endpoints $q^*_{i-1}$, $q^*_i$) and $\delta_2$ (with
endpoints $q^*_{i+1}$, $q^*_{i+2}$).  Therefore, $\Gamma$ can be factored as claimed above and we are
done by induction.\qed

\subsubsection{Factorization}
For $(T^2\times I,\xi)$ with convex boundary, we set $T_i=T^2\times \{i\}$,
$\Gamma_i=\Gamma_{T_i}$, $s_i=s(T_i)$ (slopes of the dividing sets), $r_i=r(T_i)$ (slopes of
the Legendrian rulings), and $n_i={1\over 2}(\#\Gamma_{T_i})$ (torus division number).

\begin{lemma} Let $\mbox{Tight}^0(T^2\times I,\Gamma)$ be the space of nonrotative tight
contact structures with fixed boundary condition $\Gamma=\Gamma_0\cup \Gamma_1$, $n_0\leq n_1$,
$s_0=s_1=\infty$, and $r_0=r_1=0$, and $\mathcal{G}$ is the set of dividing sets
$\Gamma_A$ on an annulus $A$ with a fixed number of endpoints on each component of $\bdry A$,
subject to the condition that such that at least two dividing curves go across from $T_0$ to
$T_1$.  There exists a bijection
$$\Psi\co \pi_0(\mbox{Tight}^0(T^2\times I,\Gamma))
\stackrel{\sim}{\rightarrow} \mathcal{G}.$$
\end{lemma}

\proof Let $\xi\in \mbox{Tight}^0(T^2\times I,\Gamma)$.  Let $A_{[0,1]}
=S^1\times\{0\}\times[0,1]$ be a horizontal convex
annulus with Legendrian boundary on $T^2\times [0,1]$. Since $\xi$ is nonrotative,
$\Gamma_{A_{[0,1]}}$ must have at least two dividing curves which go across.
Then $\Gamma_{A_{[0,1]}}$ completely
determines the isotopy type of $\xi$, since $(T^2\times I)\backslash A_{[0,1]}$ is a solid
torus which has boundary slope $-{1\over k}$ after rounding (here $2k$ is the number of
dividing curves which go across) . In particular, a tight contact structure $\xi(\Gamma_{A_{[0,1]}})$
which has dividing set $\Gamma_{A_{[0,1]}}$ is isotopic to an $S^1$--invariant tight contact structure
on $S^1\times A_{[0,1]}$, all of whose cross sections $\{pt\}\times A$ have the same dividing set
$\Gamma_{A_{[0,1]}}$.

\begin{figure}[ht!]
\centerline{
\psfrag {T0}{\smash{\rlap{\kern -2pt \raise 0pt\hbox{$T_0$}}}}
\psfrag {T1}{\smash{\rlap{\kern -2pt \raise 0pt\hbox{$T_1$}}}}
\psfrag {T2}{\smash{\rlap{\kern -2pt \raise 0pt\hbox{$T_2$}}}}
\epsfysize=2.5in\epsfbox{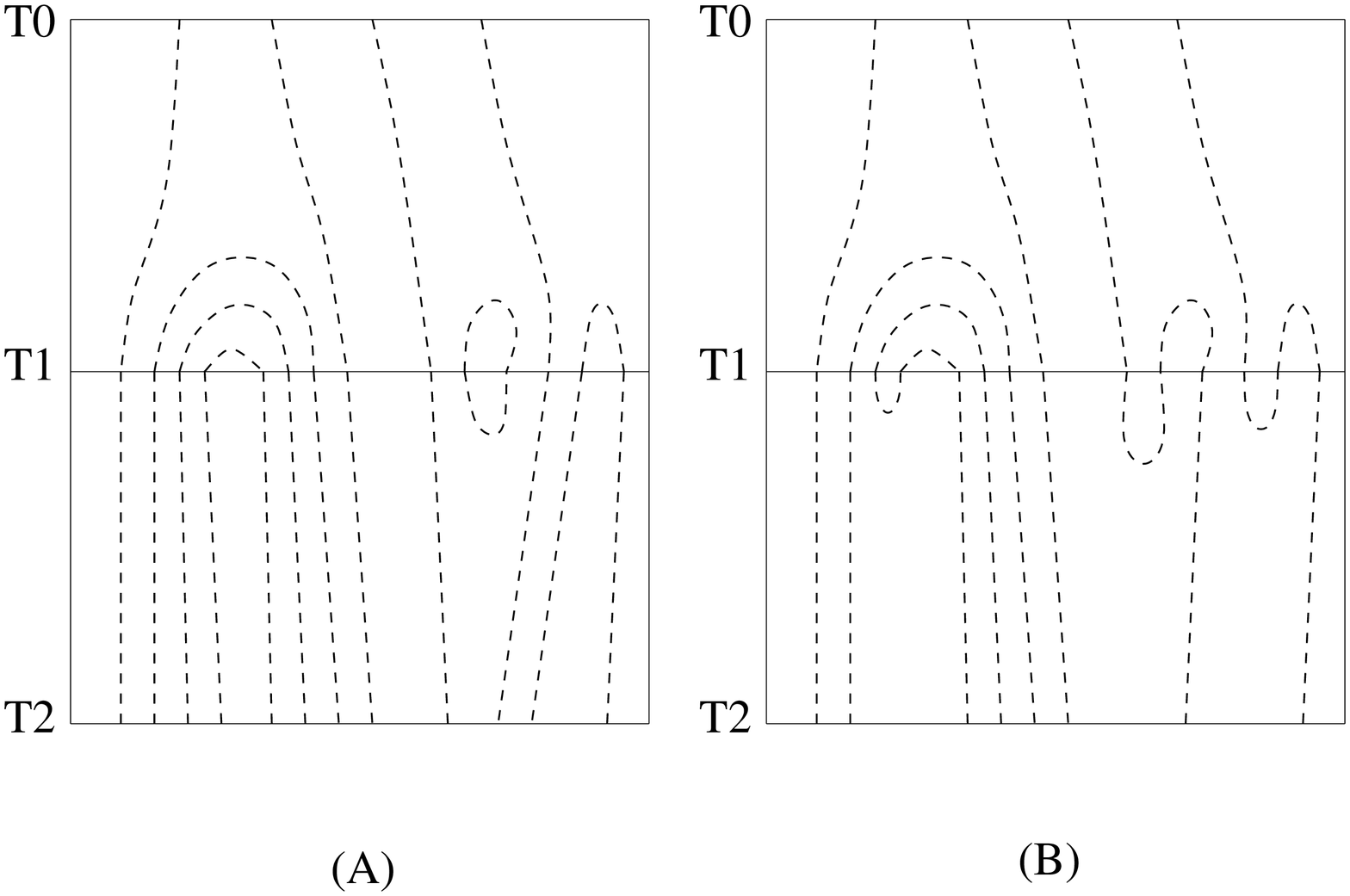}}
\caption{Gluing}
\label{matching}
\end{figure}

It remains to show that $\Gamma_{A_{[0,1]}}$ is uniquely determined by
$\xi\in \mbox{Tight}^0(T^2\times I,\Gamma)$.
Assume first that there exist no boundary-parallel components on $\Gamma_{A_{[0,1]}}$ along $T_0$.
We prove that there cannot exist $A'_{[0,1]}$ with a different $\Gamma_{A'_{[0,1]}}$.
The idea is to take advantage of the fact
that $\xi$ is $S^1$--invariant and apply dimensional reduction. We attach various
$T^2\times[1,2]$ with $n_{2}<n_1$, $s_{2}=s_1=\infty$, and no twisting, onto $T^2\times[0,1]$.
Equivalently, set  $\mathcal{A}=\Gamma_{A_{[0,1]}}$ and
consider all possible gluings to $A_{[1,2]}$ with dividing set $\Gamma'\in \mathcal{A}^*$.
The elements $\Gamma'$ correspond to all the gluings which (1) do not produce an overtwisted
disk after gluing and (2) do not produce a bypass along $T_0$ after gluing.
See Figure \ref{matching}(A) for an illustration.
Any gluing which does not produce a dividing curve on
$A_{[0,1]} \cup A_{[1,2]}$ bounding a disk will yield a universally tight structure --- this
follows from observing that both contact structures are invariant in the $S^1$-direction and
using Giroux's criterion for tightness of $I$--invariant neighborhoods of $\Sigma$.
Now apply Lemma \ref{reflexive} and obtain that $\Gamma_{A_{[0,1]}}$ is
completely determined by the isotopy type of $\xi$ (modulo holonomy).
To take show that the holonomy is the same, we use the same technique as in Proposition
\ref{nonrotative}.

In the general case we first factor $T^2\times [0,1]=(T^2\times [0,{1\over 2}])\cup
(T^2\times [{1\over 2},1])$ and $A_{[0,1]}=A_{[0,{1\over 2}]}\cup A_{[{1\over 2},1]}$ so that
$\Gamma_{A_{[0,{1\over 2}]}}$, $\Gamma_{A_{[{1\over 2},1]}}$
have no boundary-parallel components along $S^1\times\{0\}\times  \{{1\over 2}\}$. Suppose
there exists another $A_{[0,1]}'$ with the same boundary.
By passing to a large enough cover $S^1\times [-N,N]\times I$ as in
Proposition \ref{nonrotative}, we take disjoint copies $\widetilde{A_{[0,1]}}$ and $\widetilde{A_{[0,1]}'}$.
If the minimal geometric intersection number $2m=\#(\Gamma_{A_{[0,1]}}\cap (S^1\times\{{1\over 2}\})$
is not equal to $2m'=\#(\Gamma_{A_{[0,1]}'}\cap (S^1\times\{{1\over 2}\})$ (say $m>m'>0$), then we can cut along
$\widetilde{A_{[0,1]}}$, round the edges, and obtain a Legendrian curve of twisting number $-m'$ inside
a standard neighborhood of a Legendrian curve of twisting number $-m'$.  Therefore, $m=m'$.
Next, use Legendrian realization to take curves $\gamma$, $\gamma'$ on $\widetilde{A_{[0,1]}}$,
$\widetilde{A_{[0,1]}'}$  with twist number $m$, take a convex annulus interpolating from $\gamma$ to
$\gamma'$.   Since there cannot exist any boundary-parallel dividing curves (this would imply a bypass),
the dividing curves must connect between $\gamma$ and $\gamma'$.  This implies that there exists a splitting
$T^2\times [0,1]=(T^2\times [0,{1\over 2}])\cup (T^2\times [{1\over 2},1])$ which simultaneously splits
$A_{[0,1]}'=A_{[0,{1\over 2}]}\cup_\gamma A_{[{1\over 2},1]}$ and
$A_{[0,1]}'=A_{[0,{1\over 2}]}'\cup_{\gamma'} A_{[{1\over 2},1]}'$  where there are no boundary-parallel
components along $\gamma$, $\gamma'$, and $\gamma$, $\gamma'$ can be identified without loss of
generality.  Finally, apply template matching (Lemma \ref{reflexive}) to both components of the splitting. \qed

\medskip
This proves Theorem \ref{3}, Part (4).

\begin{prop} \label{factor}
There exists a unique factorization of a tight contact manifold $(M,\xi)$
with convex $\bdry M=\cup_{i=1}^n T^2_i$ so that $M=(\cup_i N_i)\cup M_0$, where
(1) $N_i\simeq T^2_i\times I$ with identical boundary slopes on both boundary components
and no twisting, (2) $T^2_i\times\{1\}=T^2_i$, (3) $T^2_i\times\{0\}$ has the minimum
possible torus division number, where the minimum is taken over all
$T^2_i\times I\subset (M\backslash \cup_{k=1}^{i-1} N_k)$
satisfying (1) and (2), and (4) $M_0=M\backslash (\cup_i N_i)$.
\end{prop}

\proof Let us show that the first factorization is unique.  Suppose  $M=N_1\cup M_0=N_1'
\cup M_0'$.  Then the cross-sectional annuli for $N_1$ and $N_1'$ must have identical dividing sets, by
using the template technique.  Therefore, $N_1$ and $N_1'$ can be matched up using an isotopy.
It then remains to show that $M_0$ and $M_0'$ are isotopic.  This follows from attaching a template
$T^2_1\times [1,2]$ such that $T^2_1\times [0,2]$ is now an $I$--invariant neighborhood of $T^2_i$.
$N_1$ and $N_1'$ are therefore isotopic.
\qed

\medskip

{\bf Proof of Theorem \ref{2}(1)}\qua
Consider $M=T^2\times I$ with convex boundary and boundary slopes $-{p\over q}$ and $-1$.
If $-{p\over q}<-1$ or $-{p\over q}=-1$ and $\phi_I>0$, then there exist nonrotative outer layers
$T^2\times [0,{1\over 3}]$ and $T^2\times [{2\over 3},1]$, where $T_{i\over 3}$, $i=0,1,2,3$, are
convex and $\#\Gamma_{1\over 3}=\#\Gamma_{2\over 3}=2$.  Moreover, Proposition \ref{factor}
indicates that the factorization is unique up to isotopy rel boundary.

If $p=q=1$ with
no twisting,  then $M$ will have an inner layer $T^2\times [{1\over 3},{2\over 3}]$ with boundary slopes $-1$ and torus
division number $n$, together with a horizontal convex annulus, all of whose dividing
curves go across.  This is the only time for $T^2\times I$ that the minimal
possible torus division number is not necessarily $1$.  In both cases, Proposition
\ref{factor} allows us to factor $M$ into an essential inner layer, together with
universally tight outer layers which can be thought of as decoration.      \qed

\section{Remarks and questions}

The results in this paper are best thought of as {\it building blocks} for a more topological
(cut-and-paste) theory of tight contact structures on 3--manifolds.  Using the techniques
presented here, we completely classify tight contact structures on the following classes of
3--manifolds in subsequent papers:
\begin{itemize}
\item Torus bundles which fiber over the circle \cite{H2}.
\item Circle bundles which fiber over closed Riemann surfaces \cite{H2}.
\item Some Seifert fibered spaces over $S^2$, such as the Poincar\'e homology sphere \cite{EH99}.
\end{itemize}

We also list some classes of 3--manifolds which are more stubborn, for which only (very weak)
partial results are known.

\begin{itemize}
\item Genus $g$ handlebodies where $g>1$.
\item Circle bundles which fiber over surfaces with boundary (even for the 3--holed sphere).
\item Seifert fibered spaces.
\item $\Sigma\times I$, where $\Sigma$ is a closed surface of genus $g>1$.
\item Surface bundles over the circle with pseudo-Anosov monodromy.
\end{itemize}

Here are a few facts and questions.

\begin{prop}
Let $M$ be a genus $g>1$ handlebody and $\Gamma$ be a dividing set for $\bdry M$.
Then $|\pi_0(\mbox{Tight}(M, \Gamma))|$ is finite.
\end{prop}

This follows from the fact that there exist $g$ compressing disks $D_1,\cdots, D_g$ so that
$M\backslash (D_1\cup \cdots \cup D_g)$ is a 3--ball, and that the number of possible dividing sets on
each $D_i$ is finite.

\begin{quest}
Can every tight $(M,\xi)$, where $M$ is a genus $g$ handlebody, be embedded inside
a symplectically semi-fillable $(M',\xi')$?
\end{quest}

By Theorem \ref{3}, when $g=1$ every tight $(M,\xi)$ can be embedded inside a lens space $L(p,q)$ with a
tight contact structure which is holomorphically fillable.  \footnote{Recently the author showed that not
every tight handlebody can be embedded inside
a symplectically semi-fillable $(M',\xi')$ \cite{H3}.}

\begin{quest}
Can every tight $(M,\xi)$, where $M=S^1\times \Sigma$ and $\Sigma$ is a 3--holed sphere, be
embedded inside a symplectically semi-fillable $(M',\xi')$?
\end{quest}

The author believes the answer is no.

\medskip
{\bf Acknowledgements}\qua I would like to thank John Etnyre for
his interest and for illuminating some of the potential pitfalls.  Thanks also go to
Francis Bonahon for a careful reading of this paper and for suggesting numerous improvements
and corrections.

\newpage

\volumenumber{5}
\volumeyear{2001}
\pagenumbers{1501}{1514}

\count0=\startpage

\addtocounter{section}{-6}

\gt\hfill      
\beginpicture
\setcoordinatesystem units <0.33truein, 0.33truein> point at 2.2 0.9
\setplotsymbol ({$\cal G$})
\plotsymbolspacing=9truept
\circulararc 315 degrees from 0 1 center at 0 0
\setplotsymbol ({$\cal T$})
\circulararc 315 degrees from 1 -1 center at 1 0
\endpicture
%
\break
{\small Volume 5 (2001) \startpage--\finishpage\ (temporary page numbers)\nl
Erratum 1\nl
Published:  24 October 2001}

\vglue 0.4truein
{\parskip=0pt\leftskip 0pt plus 1fil\def\\{\par\smallskip}{\Large
\bf Factoring nonrotative $T^2\times I$ layers}\par\medskip}   

\vglue 0.05truein 

{\parskip=0pt\leftskip 0pt plus 1fil\def\\{\par}{\sc Ko Honda}
\par\medskip}%
 
\vglue 0.03truein 

%
{\small\parskip=0pt
{\leftskip 0pt plus 1fil\def\\{\par}{\sl University of Southern California, Los Angeles, CA 90089}\par}

\vglue 5pt
\cl{Email:\stdspace\tt khonda@math.usc.edu}

\vglue 5pt
\cl{URL:\stdspace\tt http://math.usc.edu/\char126 khonda}

\vglue 10pt 

{\small\leftskip 25pt\rightskip 25pt{\bf Abstract}\stdspace
In this note we seek to remedy errors which appeared in \cite{H1} and were 
	propagated in subsequent papers.	

\vglue 7pt
{\bf AMS Classification}\stdspace 57M50; 53C15\par
\vglue 7pt
{\bf Keywords}\stdspace Tight, contact structure\par}}
\vglue 7pt


\section{Introduction}

The goal of this note is to highlight two errors which appeared in \cite{H1} and 
to provide substitutes for them.  The two incorrect statements are 
Proposition~5.8 and Part~1 of Theorem~2.2, which is a corollary of 
Proposition~5.8. The incorrect proofs of both statements appear on pages 365--366 
of \cite{H1}.  (The rest of Theorem~2.2 is unaffected by this mistake and is 
still valid.)  After making a few preliminary definitions, we will explain what 
the incorrect statements are, why they are wrong, and what can be salvaged.

In this note we assume the ambient manifold $M$ is an oriented, compact  
3-manifold and the contact structure $\xi$ on $M$ is oriented and positive, 
unless otherwise stated.   We denote the dividing set of a convex surface 
$\Sigma$ by $\Gamma_\Sigma$, and the number of connected components of 
$\Gamma_\Sigma$ by $\#\Gamma_\Sigma$. 

\subsection{}  First we recall the classification of {\it nonrotative} tight  
contact structures on $T^2\times [0,1]$.  Fix an oriented identification 
$T^2\simeq \R^2/\Z^2$, so we may talk about {\it slopes} of essential curves on 
$T^2$.  We will denote $T_t=T^2\times \{t\}$ and the slope of $\Gamma_{T_t}$ by 
$s_t$.   Let $\xi$ be a tight contact structure on $T^2\times [0,1]$ with convex 
boundary. Then $\xi$ is said to be {\it nonrotative} if all convex surfaces 
parallel to $T_0$ (or $T_1$) have dividing curves of the same slope; otherwise 
$\xi$ is said to be {\it rotative}.   An annulus $A$ in a nonrotative 
$(T^2\times [0,1],\xi)$ is {\it horizontal} if it is convex with Legendrian 
boundary, and each component of $\Gamma_{T_0}\sqcup \Gamma_{T_1}$ intersects 
$\bdry A$ exactly once. Note we may need to modify $\xi|_{T_0\sqcup T_1}$ using 
Giroux's Flexibility Theorem  (see  \cite{Gi91}) --- such modifications will 
usually be made in this note without explicit mention of the Flexibility 
Theorem.  

Recall the following, which is Lemma~5.7 of \cite{H1}.   

\begin{prop}    \label{nonrotative2}
The set of isotopy classes, rel boundary, of nonrotative tight contact 
structures on $T^2\times I$ with a fixed convex boundary, where 
$s_0=s_1=\infty$, $\#\Gamma_{T_0}=2n_0$, $\#\Gamma_{T_1}=2n_1$, and the 
characteristic foliation consists of horizontal Legendrian rulings, is in 1-1 
correspondence with isotopy classes of dividing curves $\Gamma_A$ on the 
horizontal annulus $A$, rel $\bdry A$, which consist of $n_0+n_1$ arcs 
which connect among the $2(n_0+n_1)$ fixed endpoints on $\bdry A$ ($2n_0$ along $T_0$ and $2n_1$ along 
$T_1$), at least two of which are nonseparating. \end{prop}

A connected component $\delta$ of $\Gamma_A$ is {\it nonseparating} if 
$A\setminus \delta$ is connected.

\subsection{}
Let $(M,\xi)$ be a tight contact manifold.  We define a {\it nonrotative outer 
layer} of $(M,\xi)$ to be a toric annulus $T^2\times [0,1]\subset M$ for which:

\begin{itemize}
\item $T_1$ is a boundary component of $M$, 
\item $(T^2\times [0,1], \xi|_{T^2\times [0,1]})$ is nonrotative, and 
\item $\#\Gamma_{T_0}=2$, $\#\Gamma_{T_1}=2n\geq 2$. 
\end{itemize}

Assume $(M,\xi)$ admits a factorization 
$M=(T^2\times[0,1])\cup M_0$, where $T^2\times [0,1]$ is a nonrotative outer 
layer. It was claimed (Proposition~5.8 of \cite{H1}) that such a factorization 
is unique up to isotopy, but this is hardly the case.  There is a small amount 
of flexibility in the factorization process, arising out of one case which was 
forgotten in the ``proof'' of Proposition~5.8 of \cite{H1}.  Also, in Part~1 of 
Theorem~2.2 of \cite{H1}, it was claimed that if $(T^2\times [0,1],\xi)$ is a 
tight contact manifold with convex boundary and $s_0\not=s_1$, then there exists 
a {\em unique} factorization $T^2\times [0,1]=(T^2\times[0,{1\over 3}])\cup 
(T^2\times [{1\over 3},{2\over 3}])\cup (T^2\times [{2\over 3},1])$, where 
$T_{i\over 3}$, $i=0,1,2,3$ are convex, $T^2\times[0,{1\over 3}]$ and $T^2\times 
[{2\over 3},1]$ are nonrotative, and $\#\Gamma_{T_{1/3}}=\#\Gamma_{T_{2/3}}=2$.  
The existence of such a factorization is still valid, but the uniqueness 
(purportedly a corollary of Proposition~5.8) does not hold.  Potential sources 
of this nonuniqueness will be explained in Sections~\ref{nonunique} and 
\ref{shufflable}. 

In general, it appears that the mechanism of 
factoring the nonrotative outer layer is a rather subtle one, and the 
following problem does not have a complete solution at this moment:

\begin{p}  \label{problem1}
Classify tight contact structures on $T^2\times [0,1]$ with convex boundary, in 
the case $\#\Gamma_{T_0}$ and $\#\Gamma_{T_1}$ are greater than $2$. 
\end{p} 

In this paper, we will provide partial results towards the mechanism of 
factorization.  In Section~\ref{section:general}, we introduce the notion of 
{\it disk-equivalence} and prove the following theorem:

\begin{thm}  \label{disk-equiv}
Any two nonrotative outer layers of a tight contact manifold $(M,\xi)$ 
corresponding to the same torus boundary component of $M$ are disk-equivalent.  
\end{thm}  

Theorem~\ref{disk-equiv} has the advantage that it is a general theorem which is 
sufficient for many purposes.  For example, the proofs of gluing theorems in 
\cite{H22}, which mistakenly used Proposition~5.8 of \cite{H1}, can be easily 
patched by using Theorem~\ref{disk-equiv} --- we did not need the full strength 
of the (incorrect) Proposition~5.8.  This will be explained in 
Section~\ref{section:remedy}.  

The drawback of Theorem~\ref{disk-equiv} is that the full set of nonrotative 
outer layers $T^2\times I$ for a tight contact manifold $(M,\xi)$ may not be all 
the toric annuli disk-equivalent to the initial outer layer.  In 
Section~\ref{special-cases} we exhibit two extreme cases: the {\it shufflable 
case}, where all the disk-equivalent toric annuli are represented, and the {\it 
universally tight case}, where the full set of nonrotative outer layers is 
substantially smaller.   

There are two general strategies for analyzing the factorization process.   The 
easier strategy is to probe the tight contact structure on $(M,\xi)$ {\it 
externally}.  This involves attaching nonrotative $T^2\times I$ layers from 
outside (called {\it templates}), and weighing their effect on the resulting 
glued-up contact manifold.  The key is to keep track of the layers which glue to 
give tight contact manifolds, as well as those which glue to give overtwisted 
contact manifolds.  The other strategy is an internal probe, called {\it state 
traversal}, explained in \cite{H32}.  This internal probe, although 
usually more difficult to implement in practice, yields more complete 
information than that of {\it template attaching}.  In this note, we shall 
restrict ourselves to the (much easier) template method.  State traversal 
should yield a complete solution to Problem~\ref{problem1}, but the 
combinatorics seem highly nontrivial.

\section{General case}   \label{section:general}

In this section, we prove the general result on nonrotative outer layers, namely 
Theorem~\ref{disk-equiv}.  Theorem~\ref{disk-equiv} has the advantage that it 
has a nice formulation in terms of {\it disk-equivalence} which is useful in 
practice.  It also admits a relatively elementary proof using template 
attaching.  

\subsection{}  Consider two nonrotative outer layers $N=T^2\times [0,1]$ and  
$N'=(T^2\times [0,1])'$ of $(M,\xi)$, where $T_1=T_1'$ is a boundary component 
of $M$.  Let $A$ and $A'$ be the corresponding horizontal annuli with $\bdry 
A=\delta_0\sqcup \delta_1$ and $\bdry A'=\delta'_0\sqcup \delta'_1$.  After 
sliding $\delta'_1$ along $T_1=T_1'$ if necessary, we may assume that 
$\delta'_1=\delta_1$ and $\Gamma_{A}\cap \delta_1 = \Gamma_{A'}\cap \delta_1$.   
Now, we say $A$ and $A'$ (or $N$ and $N'$) are {\it disk-equivalent} if there 
exist embeddings $\phi:A\hookrightarrow D^2$ and $\phi':A'\hookrightarrow D^2$ 
where $\phi(\delta_1)=\phi'(\delta_1)=\bdry D^2$ and 
$\phi|_{\delta_1}=\phi'|_{\delta'_1}$,  such that $\Gamma_{D^2}$ on $D^2$, 
obtained from $\phi(\Gamma_{A})$ by connecting the two endpoints of 
$\phi(\Gamma_{A}\cap \delta_0)$ via an arc in $D^2\setminus \phi(A)$, and 
$\Gamma'_{D^2}$, obtained similarly from $\phi'(\Gamma_{A'})$, are isotopic rel 
$\bdry D^2$.  

We are now ready to prove Theorem~\ref{disk-equiv}.

\begin{proof}[Proof of Theorem~\ref{disk-equiv}]
Consider the factorization $M=N\cup M_0$, where $N=T^2\times [0,1]$ 
is a nonrotative outer layer and  $A_{[0,1]}$ is its horizontal annulus.  We 
prove that $A_{[0,1]}'$ corresponding to another nonrotative outer layer $N'$ is 
disk-equivalent to $A_{[0,1]}$.  Write $\bdry A_{[a,b]}=\delta_a\sqcup 
\delta_b$.

Let $\mathcal{T}_{A_{[0,1]}}$ (resp.\ $\mathcal{T}$) be the set of isotopy 
classes of nonrotative tight contact structures $(T^2\times [1,2],\zeta)$ with 
a fixed boundary characteristic foliation and $\#\Gamma_{T_2}=2$, which glue to 
$(N=T^2\times [0,1],\xi|_N)$ to yield a tight contact structure on 
$T^2\times[0,2]$ which is $I$-invariant (resp.\ a tight contact structure on 
$M\cup (T^2\times [1,2])$).   Here, the {\it $I$-invariant} tight contact 
manifold is isomorphic to an invariant neighborhood of a convex surface $T_2$ 
(or $T_0$).   By Proposition~\ref{nonrotative2}, a nonrotative 
$(T^2\times [1,2],\zeta)$ is characterized by the dividing set of its horizontal 
annulus $A_{[1,2]}$.   Any $\Gamma_{A_{[1,2]}}$ will have exactly two endpoints 
along $\delta_2$ and exactly two nonseparating arcs.   Associate to 
$\mathcal{T}_{A_{[0,1]}}$ (resp.\ $\mathcal{T}$) the corresponding set of 
isotopy classes $\mathcal{A}_{A_{[0,1]}}$ (resp.\ $\mathcal{A}$)  of 
$\Gamma_{A_{[1,2]}}$.    Let $A_{[0,2]}=A_{[0,1]}\cup A_{[1,2]}$ be the 
horizontal annulus for $T^2\times [0,2]$, where we assume that $A_{[0,1]}$ and 
$A_{[1,2]}$ have common boundary $\delta_1$. Now, $\Gamma_{A_{[1,2]}}\in 
\mathcal{A}_{A_{[0,1]}}$ if and only if $\Gamma_{A_{[0,2]}}$ consists of exactly 
two parallel nonseparating arcs. Clearly, $\mathcal{A}_{A_{[0,1]}}\subset 
\mathcal{A}$, since the $I$-invariance of $T^2\times [0,2]$ implies there is a 
contact diffeomorphism $(M,\xi)\cup ((T^2\times[1,2]),\zeta)\simeq 
(M_0,\xi|_{M_0})$.  Of course, $\mathcal{A}$, unlike $\mathcal{A}_{A_{[0,1]}}$, 
depends on the ambient $(M,\xi)$, and $\mathcal{A} - \mathcal{A}_{A_{[0,1]}}$ 
may or may not contain certain $\Gamma_{A_{[1,2]}}$ for which 
$\Gamma_{A_{[0,2]}}$ contains (necessarily homotopically essential) closed 
curves.  See Figure~\ref{fig7} for  various possibilities of 
$\Gamma_{A_{[1,2]}}$.    

\begin{figure}[ht!] 			
{\epsfysize=1.5in\centerline{\epsfbox{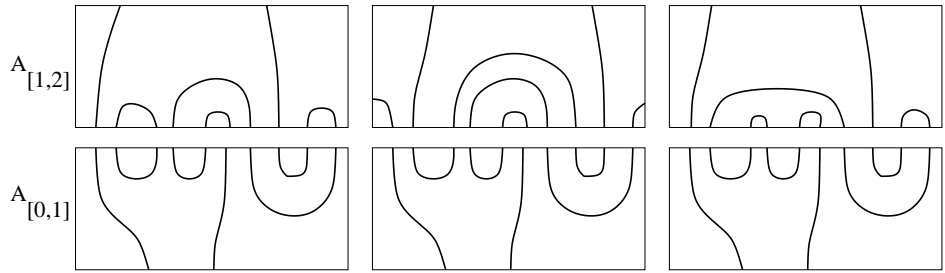}}} 		
\caption{In all the figures, the sides are identified.  The right-hand $\Gamma_{A_{[1,2]}}$ is in 
$\mathcal{T}_{A_{[0,1]}}$, the left-hand diagram is not in $\mathcal{T}$, and it 
cannot be determined simply by looking at $A_{[0,2]}$ whether the middle is in 
$\mathcal{T}-\mathcal{T}_{A_{[0,1]}}$.} 
\label{fig7} 
\end{figure}

The induction is done by fixing $(M_0,\xi|_{M_0})$ and inducting on 
$\#\Gamma_{T_1}=2n$ over the space of all nonrotative outer layers $N=T^2\times 
[0,1]$ with $\#\Gamma_{T_0}=2$.  Note that all nonrotative $N=T^2\times [0,1]$ 
with $\#\Gamma_{T_0}=2$ can be embedded inside an $I$-invariant neighborhood of 
$T_0$ by {\it folding} (see Section 5.3 of \cite{H1}), so all contact structures 
on $M_0\cup N$ constructed this way are tight.  When $n=1$, the nonrotative 
outer layer is clearly unique.  Therefore, assume the theorem is true for 
$\#\Gamma_{T_1}=2n$, and we prove it for $\#\Gamma_{T_1}=2(n+1)$.  There are two 
cases: either $\Gamma_{A_{[0,1]}}$ has at least two $\bdry$-parallel curves or  
there is only one $\bdry$-parallel curve.

Suppose first that there are at least two $\bdry$-parallel curves on 
$A_{[0,1]}$.   Let $\gamma$ be an arc on $A_{[1,2]}$ whose 
endpoints are consecutive points of  $\Gamma_{A_{[0,1]}}\cap \delta_1$, ie, 
$\gamma$ is {\it $\bdry$-parallel}. If the endpoints of $\gamma$ coincide with 
the endpoints of a $\bdry$-parallel curve of $A_{[0,1]}$, then, for any 
completion of $\gamma$ to a dividing set $\Gamma_{A_{[1,2]}}\supset \gamma$, 
the gluing $A_{[0,1]}\cup A_{[1,2]}$ corresponds to an overtwisted contact 
structure.  On the other hand, if the endpoints of $\gamma$ are not (i) the two 
endpoints of the nonseparating curves of $\Gamma_{A_{[0,1]}}$ and not (ii) the 
two endpoints of a $\bdry$-parallel curve of $\Gamma_{A_{[0,1]}}$, then $\gamma$ 
can be completed into some $\Gamma_{A_{[1,2]}}\in \mathcal{A}$.  We now 
summarize the completability of $\gamma$ to an element in $\mathcal{A}$: unknown 
if endpoints are (i),  no if endpoints are (ii), and yes otherwise. (Here 
``unknown'' means that it depends on whether adding an extra $\pi$-twisting 
$T^2\times I$ layer to $M_0$ yields a tight contact structure or an overtwisted 
contact structure.)   Now, since there are at least two $\bdry$-parallel curves 
of $\Gamma_{A_{[0,1]}}$, there are at least two $\bdry$-parallel $\gamma$ which 
cannot be completed to an element of $\mathcal{A}$, and at least one of them 
must have the same endpoints as a $\bdry$-parallel curve of 
$\Gamma_{A'_{[0,1]}}$.  (This follows from repeating the same argument for 
$A'_{[0,1]}$ instead of $A_{[0,1]}$.) Thus, there is a common $\bdry$-parallel 
position for both $A_{[0,1]}$ and $A'_{[0,1]}$.  Now, attach a horizontal 
annulus with $2n$ nonseparating curves and one $\bdry$-parallel dividing curve 
$\gamma$ right next to the common $\bdry$-parallel position of $A_{[0,1]}$ and 
$A'_{[0,1]}$ as in Figure \ref{fig8}, and use the inductive step.     

\begin{figure} [ht!]			
	{\epsfysize=2in\centerline{\epsfbox{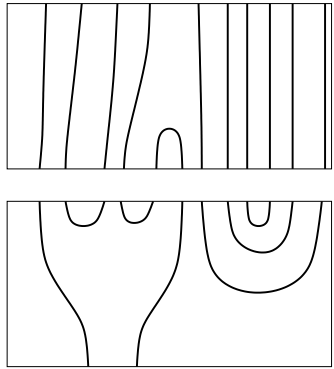}}} 		
	\caption{Inductive step} 			
	\label{fig8}
\end{figure}

\begin{figure} [ht!]		
	{\epsfysize=2in\centerline{\epsfbox{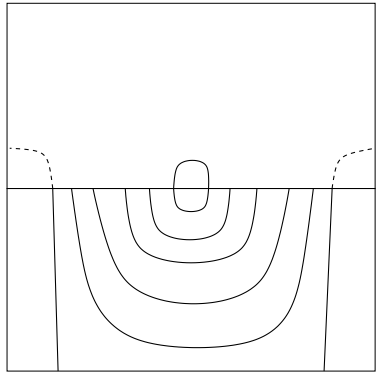}}} 	
	\caption{Only one $\bdry$-parallel dividing curve. The bottom annulus is 
	$A_{[0,1]}$ and the top one is $A_{[1,2]}$.} 			
	\label{fig6}
\end{figure}

Suppose now that there exists only one $\bdry$-parallel arc of 
$\Gamma_{A_{[0,1]}}$.   Then the two nonseparating curves must be consecutive 
(ie, one of the regions of $A_{[0,1]}$ divided by these two curves does not 
contain any other dividing curves), and all the separating curves must be nested 
concentrically around the one $\bdry$-parallel dividing curve.  See Figure 
\ref{fig6}. The $\bdry$-parallel arc $\gamma$ on $A_{[1,2]}$ satisfying (ii) is 
at the center (solid line), and $\gamma$ satisfying (i) is given by dotted 
lines.    Then $\Gamma_{A'_{[0,1]}}$ satisfies one of the following:
\begin{itemize} 
\item $\Gamma_{A'_{[0,1]}}=\Gamma_{A_{[0,1]}}$.
\item The positions of (i) and (ii) are reversed.
\item Positions (i), (ii) for $\Gamma_{A_{[0,1]}}$ are both (ii) for 
$\Gamma_{A'_{[0,1]}}$.
\end{itemize}
In each case, $\Gamma_{A_{[0,1]}}$ and $\Gamma_{A'_{[0,1]}}$ are 
disk-equivalent. \end{proof}

Note that Theorem~\ref{disk-equiv} does not completely address exactly which 
nonrotative outer layers exist for a given $(M,\xi)$.  

\begin{cor}     \label{iso-core}
Given two factorizations $M=N\cup M_0$ and $M=N'\cup M_0'$ of a tight $(M,\xi)$, 
where $N$, $N'$ are nonrotative outer layers corresponding to the same torus 
boundary component of $M$, there exists an isomorphism $(M_0,\xi|_{M_0})\simeq
(M_0',\xi|_{M_0'})$. 
\end{cor}
                                                    
\begin{proof}  The actual isomorphism is not an arbitrary isomorphism, but an 
isotopy in the following sense.  Let $(T^2\times[1,2],\zeta)$ be an element of 
$\mathcal{T}_{A_{[0,1]}}$ as in the proof of Theorem~\ref{disk-equiv}.  Then there 
exists a contact isotopy of $(M_0,\xi|_{M_0})$ to $(M,\xi)\cup 
(T^2\times[1,2],\zeta)$ inside $M\cup (T^2\times [1,2])$.  This is clear from 
the $I$-invariance of $N\cup (T^2\times [1,2])$.  Now we claim that the 
disk-equivalence of $A_{[0,1]}$ and $A'_{[0,1]}$ implies that 
$N'\cup (T^2\times [1,2])$ is $I$-invariant, thus proving the contact isotopy of 
$(M_0',\xi|_{M_0'})$ to $(M,\xi)\cup 
(T^2\times[1,2],\zeta)$ inside $M\cup (T^2\times [1,2])$.  Write 
$A'_{[0,2]}=A'_{[0,1]}\cup A_{[1,2]}$,  $\bdry 
A_{[0,2]} =\delta_0\sqcup \delta_2$, and $\bdry A'_{[0,2]}=\delta'_0\sqcup 
\delta'_2$.   We then complete $A_{[0,2]}$ (resp.\ $A'_{[0,2]}$) by attaching a 
disk $D$ and (resp.\ $D'$) along $\delta_0$ (resp.\ $\delta'_0$).  By the 
disk-equivalence, the dividing sets on $A_{[0,2]}\cup D$ and $A'_{[0,2]}\cup D'$ 
are identical and consist of exactly one $\bdry$-parallel arc along 
$\delta_2=\delta'_2$.    This in turn implies that, after removing $D'$ from 
$A'_{[0,2]}\cup D'$, $\Gamma_{A'_{[0,2]}}$ consists of exactly two nonseparating 
arcs.   This completes the proof. \end{proof}

\subsection{} \label{section:remedy}
In this section we seek to remedy some tightness proofs in \cite{H22} which were 
affected by the misuse of unique factorizations for nonrotative outer layers.  
The situation we are interested in is the following.  Let $(M,\xi)$ be a contact 
manifold and $T\subset M$ an incompressible torus.  Using state traversal in 
\cite{H22} and \cite{H32}, we want to determine whether $(M,\xi)$ is tight.  When 
we use this method, we start with $T$ convex and for which it is easy to 
determine that $(M\backslash T,\xi|_{M\backslash T})$ is tight.  Successively we 
find $T'$ isotopic to and disjoint from $T$, and ask whether $(M\backslash T', 
\xi|_{M\backslash T'})$  is tight.  If yes, then we let $T'$ be the new $T$, and 
continue.  If tightness is preserved for all possible $T'$, then $(M,\xi)$ is 
tight.  Usually, the initial state consists of $\#\Gamma_T=2$, but, during the 
course of the state transitions, $\#\Gamma_{T'}$ may become large.  The 
following theorem allows us to avoid these more complicated states.

\begin{thm}  \label{patchup}
It is sufficient to verify the following in order to prove the tightness of the 
contact manifold $(M,\xi)$ using state tranversal:
\be
\item $\xi|_{M\backslash T'}$ is tight for every convex $T'$ 
with $\#\Gamma_{T'}=2$, obtained from $T$ via a sequence of bypass moves, each 
of which leaves $\#\Gamma=2$.  
\item Let $T''$ be a convex torus isotopic to $T$ with tight $\xi|_{M\backslash 
T''}$.  Let $T^2\times [-0.5,0.5] \hookrightarrow M$ be a toric annulus with 
$T_0=T''$ and nonrotative $T^2\times[-0.5,0]$ and $T^2\times[0,0.5]$.  Then 
there exists an extension to $T^2\times [-1,1]\hookrightarrow M$ where 
$T^2\times [-1,0]$ and $T^2\times[0,1]$ are nonrotative outer layers in 
$M\backslash T''$.  In particular, $\#\Gamma_{T_{-1}}=\#\Gamma_{T_1}=2$. 
\ee 
\end{thm}

\begin{proof}  The smallest state transition unit $T\sa T''$ consists of 
attaching a bypass along $T$ to obtain $T''$.  Hence, every pair 
$T$, $T''$ of isotopic tori is related by a sequence of bypass attachments.  
Suppose that  $(M\backslash T', \xi|_{M\backslash T'})$ is tight for every 
convex $T'$ with $\#\Gamma_{T'}=2$, obtained from $T$ via a sequence of bypass 
moves which do not change $\#\Gamma$.   Observe that if $T=\Sigma_0\sa 
\Sigma_1\sa \cdots \sa \Sigma_k$ is the sequence of bypass moves which 
extricates the original $T$ from a candidate overtwisted disk, then there will 
exist intervals $\Sigma_i\sa \cdots \sa \Sigma_j$ where 
$\#\Gamma_{\Sigma_i}=\#\Gamma_{\Sigma_j}=2$ and $\#\Gamma_{\Sigma_l}>2$ 
inbetween, or half-intervals $\Sigma_i\sa \cdots \sa \Sigma_k$, where 
$\#\Gamma_{\Sigma_i}=2$ and $\#\Gamma_{\Sigma_l}>2$ thereafter.  We will prove 
that the state transitions when $\#\Gamma>2$ are rather superficial, and that 
$(M\setminus \Sigma_i,\xi|_{M\setminus \Sigma_i})\simeq (M\setminus 
\Sigma_j,\xi|_{M\setminus \Sigma_j})$.  

We inductively assume the following:

\be
\item[(A)] $T''$ is one of the $\Sigma_l$ between $\Sigma_i$ and $\Sigma_j$ (or 
$\Sigma_k$). 
\item[(B)] $(M\setminus T'',\xi|_{M\setminus T''})$ is tight.  
\item[(C)]  There exist nonrotative layers $T^2\times [-1,0]$, $T^2\times [0,1]$ 
with $T_0=T''$ and $\#\Gamma_{T_{-1}}=\#\Gamma_{T_1}=2$, and such that 
$T^2\times [-1,1]$ is $I$-invariant.  
\item[(D)] There is an isomorphism $$(M\setminus \Sigma_i,\xi|_{M\setminus 
\Sigma_i}) \simeq (M\setminus (T^2\times [-1,1]), \xi|_{M\setminus (T^2\times 
[-1,1])}).$$ 
\ee
Let $A_{[-1,0]}$ and $A_{[0,1]}$ be the horizontal annuli corresponding 
to $T^2\times [-1,0]$ and $T^2\times [0,1]$.

Let $(T^2\times [-0.5,0])'$ be the layer between $\Sigma_{l}=T''$ 
and $\Sigma_{l+1}$.  It is nonrotative because $\#\Gamma_{\Sigma_l}>2$ and we 
are considering a single bypass move from $\Sigma_l$ to $\Sigma_{l+1}$. The 
hypotheses of the theorem guarantee an extension to $(T^2\times [-1,0])'$, a 
nonrotative outer layer of $M\backslash T''$.  There also exists a nonrotative 
outer $(T^2\times [0,1])'$ on the other side of $T''$.  Call the corresponding 
new horizontal annuli $A'_{[-1,0]}$ and $A'_{[0,1]}$. (Also let 
$A'_{[-1,1]}=A'_{[-1,0]}\cup A'_{[0,1]}$.)  

The key is to prove that the new layer $(T^2\times [-1,1])'$ containing 
$\Sigma_{l+1}$ is $I$-invariant.  This is done by completing $A_{[-1,0]}$ to a 
disk $D_1$, $A_{[0,1]}$ to a disk $D_2$, and likewise forming $D_1'$ and $D_2'$ 
from $A_{[-1,0]}'$ and $A_{[0,1]}'$.  
If we put $D_1$ and $D_2$ together to form $S^2$ so the dividing curves match 
up, then there is exactly one dividing curve, since $\Gamma_{A_{[-1,1]}}$ 
consists of two parallel nonseparating curves.  (The corresponding toric annulus 
is $I$-invariant.) Now use Theorem~\ref{disk-equiv} to see that $D_1'\cup D_2'$ 
must also consist of exactly one dividing curve, due to disk-equivalence.  Now, 
$\Gamma_{A'_{[-1,1]}}$ is obtained by removing two small disks from $D_1'\cup 
D_2'$, each containing a short arc of the dividing set.  Therefore, 
$\Gamma_{A'_{[-1,1]}}$ must consist of parallel nonseparating curves.    This 
proves that Condition C of the inductive step also holds for $\Sigma_{l+1}$.  
Next, Condition D is satisfied, since $$(M\setminus \Sigma_i,\xi|_{M\setminus 
\Sigma_i}) \simeq (M\setminus (T^2\times [-1,1]), \xi|_{M\setminus (T^2\times 
[-1,1])}),$$ and $$(M\setminus (T^2\times [-1,1]), \xi|_{M\setminus (T^2\times 
[-1,1]}))\simeq (M\setminus (T^2\times [-1,1])', \xi|_{M\setminus (T^2\times 
[-1,1])'}),$$ due to Corollary~\ref{iso-core}.   Condition B is now obvious, 
since $(M\setminus \Sigma_{l+1}, \xi|_{M\setminus \Sigma_{l+1}})$ is obtained 
from  $(M\setminus \Sigma_i,\xi|_{M\setminus \Sigma_i})$ by folding.
\end{proof}

The following suffices for the purposes of gluing in \cite{H22}.

\begin{cor}
Let $M=(T^2\times[0,1])/\sim$ be a $T^2$-bundle over $S^1$, obtained by 
identifying $T_0\sim T_1$, and let $\xi$ be a contact structure on $M$.  
If $\xi|_{T^2\times [0,1]}$ is a rotative tight contact structure, then 
$\xi|_M$ is tight if Condition~1 of Theorem~\ref{patchup} is satisfied. 
\end{cor}

\begin{proof}
Let $T= T_0=T_1$.  Then $\xi|_{M\setminus T}$ is rotative and any pair of  
nonrotative layers  $(T^2\times [0,0.1])\sqcup (T^2\times[0.9,1])$ can be 
extended to a pair of nonrotative outer layers $(T^2\times [0,0.2])\sqcup 
(T^2\times[0.8,1])$ using bypasses and the Imbalance Principle \cite{H1}.  
Moreover, for each state transition $T\sa T''$, if $\xi|_{M\setminus T}$ is 
rotative, then so is $\xi|_{M\setminus T''}$. \end{proof}

\section{Special cases}  \label{special-cases}

In this section we assume the following:

\medskip

{\bf Extendability Condition}\qua  Let $(M,\xi)$ be a tight contact manifold 
with convex boundary $\bdry M$, one component of which is a torus $T$.  
Assume there exists a  factorization $M=(T^2\times[-1,1])\cup M_0$, where 
$T_1=T$,  $s_{-1}=0$, $s_1=-\infty$, $\Gamma_{T_{-1}}=2$, 
$\Gamma_{T_1}>2$, and every convex torus in $T^2\times[-1,1]$ parallel to 
$T_{-1}$ (or $T_1$) has slope $s$ satisfying $-\infty \leq s\leq 0$.   

\medskip

Let us call such a $T^2\times [-1,1]$ a {\it rotative outermost layer}.  Note 
that the Extendability Condition is very similar to the ``quasi-pre-Lagrangian'' 
condition in Colin \cite{Co99a}.

\subsection{} \label{nonunique} Here we present the first sources of 
nonuniqueness of nonrotative outer layers.  Suppose $(M,\xi)$ is universally 
tight and satisfies the Extendability Condition.  Consider a rotative 
outermost layer $T^2\times [-1,1]\subset M$, where
$s_1=\infty$ and $s_{-1}=0$. Consider a factorization of $T^2\times [-1,1]$ into 
$T^2\times [-1,0]$ and $T^2\times [0,1]$, where the first is a {\it basic slice} 
(ie, contactomorphic to $(T^2\times [-1,0],\xi)$ with convex boundary, 
$\#\Gamma_{T_{-1}}=\#\Gamma_{T_0}=2$, $s_{-1}=0$, $s_0=-\infty$, and every 
convex surface parallel to $T_0$ has dividing curves of slope $s$ satisfying 
$-\infty\leq s\leq 0$) and the second is a nonrotative outer layer.    Let 
$A_{[0,1]}$ be the horizontal annulus for $T^2\times [0,1]$ and $A_{[-1,0]}$ be 
the ``horizontal annulus'' for $T^2\times [-1,0]$ in the sense that $A$ is 
convex with efficient Legendrian $\bdry A_{[-1,0]}=\delta_{-1}\sqcup \delta_0$ 
of slope $0$ on $T_{-1}$ and $T_0$.     Here, a closed curve $\gamma$ on a 
convex surface $\Sigma$ is {\it efficient} if $\gamma\pitchfork \Gamma_\Sigma$ 
and the geometric intersection number $|\gamma\cap\Gamma_\Sigma|$ equals the 
actual number of intersection points.  Let $\eta_1,\cdots,\eta_k$ be the 
`innermost' dividing curves on $A_{[0,1]}\cup A_{[-1,0]}$, ie, there exists an 
arc from $\eta_i$ to $\delta_{-1}$ which intersects no other dividing curve 
except perhaps for closed essential dividing curves on $A_{[-1,0]}$ (if they 
exist).  Then the various nonrotative outer layers are obtained by truncating 
some $\eta_i$. 

\begin{figure} [ht!]		
	{\epsfysize=2in\centerline{\epsfbox{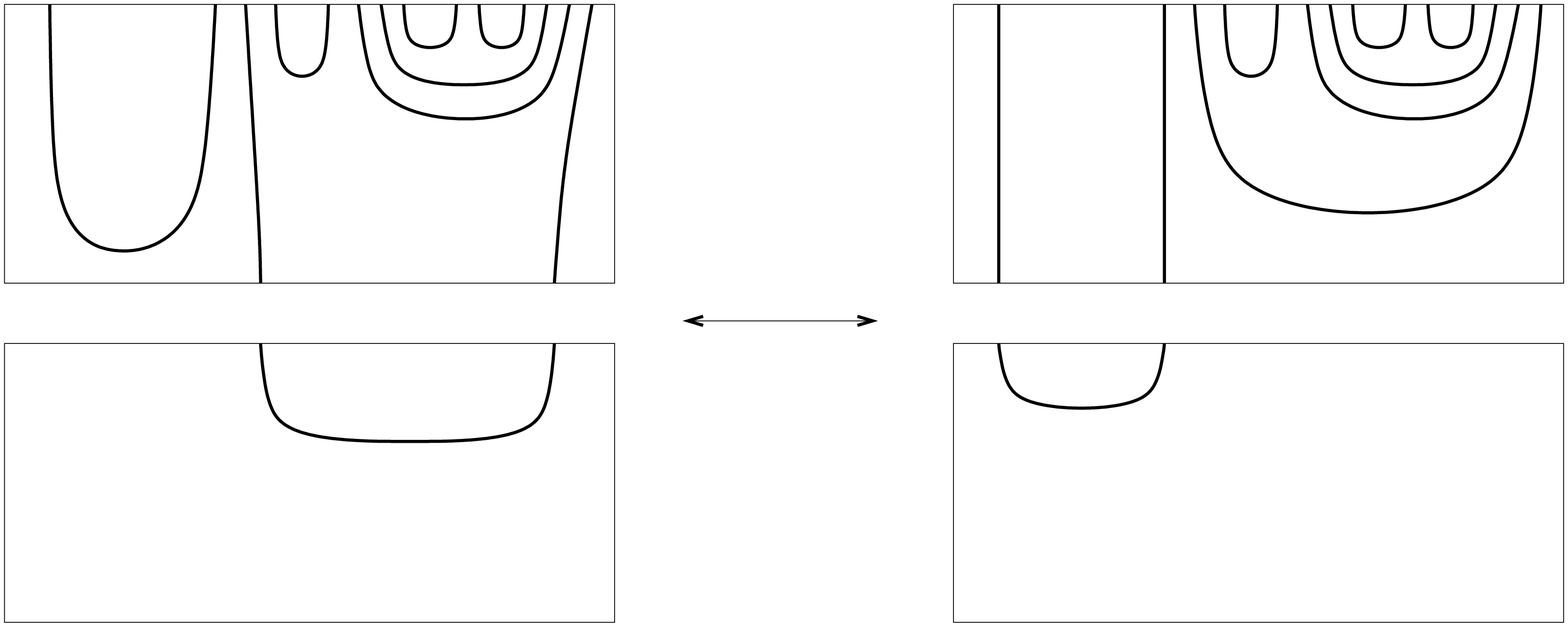}}} 	
	\caption{Equivalence in the universally tight case. The top annulus
	is $A_{[0,1]}$ and the bottom annulus is $A_{[-1,0]}$. } 		
	\label{fig4}
\end{figure}

\subsection{} \label{shufflable} 

Next we consider the following situation, which we call the {\it shufflable 
case}.  

\medskip

{\bf Assumption}\qua  Let $(M,\xi)$ be a tight contact manifold with convex 
boundary and $T$ a torus component of $\bdry M$. Suppose there exists a layer 
$T^2\times [-2,1]\subset M$ with $T_1=T$, for which $s_{-2}={1\over 2}$, 
$s_{-1}=0$, $s_0=s_1=\infty$, $\#\Gamma_{T_{-2}} = \#\Gamma_{T_{-1}} = 
\#\Gamma_{T_0}=2$, and $\#\Gamma_{T_1}=2n$.  Let
$T^2\times[-2,-1]$ and $T^2\times[-1,0]$ be basic slices, and let $T^2\times 
[0,1]$ be a nonrotative outer layer.  Moreover, assume that the relative Euler 
classes of $T^2\times [-2,-1]$ and $T^2\times [-1,0]$ are $\pm(1,1)$, 
$\mp(1,1)$, respectively.  These two basic layers can be switched via a contact 
isotopy, which is called {\it shuffling} in \cite{H1}.  Therefore, if we have 
such a $T^2\times [-2,1]$-layer, we say we are in the {\it shufflable case.} 
\medskip

\begin{figure} [ht!]		
	{\epsfysize=2in\centerline{\epsfbox{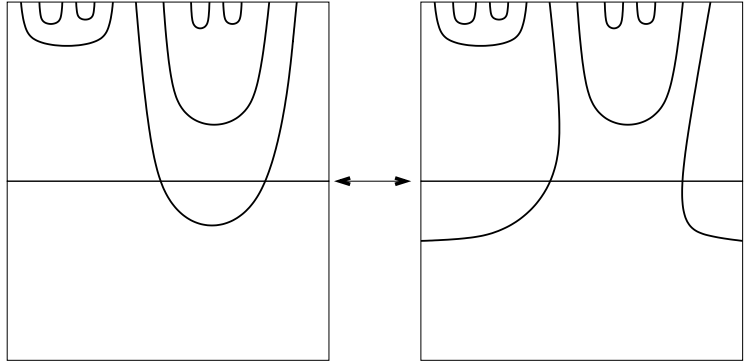}}} 	
	\caption{Equivalences in the shufflable case} 		
	\label{fig5}
\end{figure}

In the shufflable case, the rotative outermost layer is certainly not unique, as 
can be seen from Figure~\ref{fig5}. In other words, there is a clear equivalence 
relation, where the dividing curve configuration for $A_{[-1,0]}$ is substituted 
by the other possibility (ie, coming from $A_{[-2,-1]}$ after shuffling).  

If we combine moves described in Section~\ref{nonunique} with the moves 
described in Figure~\ref{fig5}, it is clear that {\it all} the configurations of 
$A_{[0,1]}$ disk-equivalent to the initial one are realized.  Combining this 
with Theorem~\ref{disk-equiv}, we obtain the following:

\begin{prop}  Let $(M,\xi)$ be a tight contact manifold with convex boundary 
$\bdry M$ and let $T$ be a torus component of $\bdry M$.  Suppose $M$ is 
shufflable along $T$.  If we fix a nonrotative outer layer $N=T^2\times [0,1]$ 
with $T_1=T$ and let $A_{[0,1]}$ be its horizontal annulus, then the set of 
isotopy classes of nonrotative outer layers (rel boundary) for $(M,\xi)$ along 
$T$ is in 1-1 correspondence with the set of isotopy classes of dividing 
multicurves (rel boundary) disk-equivalent to  $\Gamma_{A_{[0,1]}}$. \end{prop}

\subsection{}  \label{basic2}
The following is the analog of Proposition~\ref{nonrotative2} for rotative 
outermost layers.

\begin{lemma}   \label{rotative-outermost}
Let $(M=T^2\times[-1,1],\xi)$ be a rotative outermost layer. Then there 
exists a unique dividing set $\Gamma_{A_{[-1,1]}}$, modulo closed curves which 
are parallel to the boundary. \end{lemma}                                                                    

\begin{proof} 
We take $s_{-1}=0$, $s_1=\infty$, and $\#\Gamma_{T_1}>2$.  As in the proof of 
Theorem~\ref{disk-equiv}, consider the set $\mathcal{T}$ of nonrotative tight 
contact structures $(T^2\times [1,2],\zeta)$ with $\#\Gamma_{T_2}=2$, which glue 
to $(M=T^2\times [-1,1],\xi)$ to yield a tight contact structure on 
$T^2\times[-1,2]$.  The key difference between this case and 
Theorem~\ref{disk-equiv} is that it is possible to determine $\mathcal{T}$ and 
its corresponding $\mathcal{A}$ precisely.  That is, $\mathcal{A}$ consists of 
all $\Gamma_{A_{[1,2]}}$ for which $\Gamma_{A_{[-1,2]}}$ does not have any 
homotopically trivial dividing curves. --- in other words, the ``unknown'' 
gluings which produced the middle configuration in Figure~\ref{fig7} are now 
known to be tight gluings.   Elements of $\mathcal{A}$ correspond to $(T^2\times 
[1,2],\zeta)$, whose attachment makes $T^2\times[-1,2]$ either into a {\it basic 
slice} or adds extra twisting by a multiple of $\pi$.   

Now, we want to prove that if $A_{[-1,1]}$ and $A'_{[-1,1]}$ are two horizontal 
annuli for $T^2\times [-1,1]$, then $A_{[-1,1]}=A'_{[-1,1]}$ modulo parallel 
closed essential curves.   This is proved by induction on $\#\Gamma_{T_1}$.  If 
$\#\Gamma_{T_1}=2$, then there are two possibilities for $\Gamma_{A'_{[-1,1]}}$ 
modulo parallel closed essential curves, corresponding to the two possible 
positions for $\bdry$-parallel dividing curves.  In this step only, we attach 
templates which are basic slices $(T^2\times [1,2],\zeta')$ (not nonrotative 
layers) with $s_1=\infty$ and $s_2=0$, and corresponding ``horizontal'' annuli 
$A_{[1,2]}$.  The two basic slices are also distinguished by the positions of 
the $\bdry$-parallel dividing curves along $\delta_1$.  (As before, we are 
assuming that $\bdry A_{[-1,1]}=\delta_{-1}\sqcup \delta_1$ and $\bdry 
A_{[1,2]}=\delta_1\sqcup \delta_2$.  Note they have a common boundary 
$\delta_1$.)  Since the gluing is tight if and only if a closed homotopically 
trivial curve does not appear on $A_{[-1,2]}$, the two possible 
$\Gamma_{A'_{[-1,1]}}$ can be distinguished using templates.

Next, assume inductively that the claim holds for $\#\Gamma_{T_1}=2n$. 
Let $\#\Gamma_{T_1}=2(n+1)$.  Now any arc $\gamma$ on $A_{[1,2]}$ with 
consecutive endpoints on $\delta_1\cap \Gamma_{A_{[-1,1]}}$ can be extended to 
some $\Gamma_{A_{[1,2]}}\in \mathcal{A}$, if and only if the endponts of 
$\gamma$ are not the endpoints of a $\bdry$-parallel dividing curve of 
$A_{[-1,1]}$.  This implies that the set of $\bdry$-parallel curves must be the 
same for $A_{[-1,1]}$ and $A_{[-1,1]}'$.  We then reduce to the case 
$\#\Gamma_{T_1}=2n$ in the same manner as in the proof of 
Theorem~\ref{disk-equiv}.  \end{proof}

\subsection{}  The argument in Section~\ref{basic2} generalizes to the case where 
$(M,\xi)$ is universally tight.

\begin{prop}   \label{univ-tight}
If $(M,\xi)$ is universally tight and satisfies the Extendability Condition, and 
$\bdry M$ is an incompressible torus, then any two rotative outermost layers 
are contact diffeomorphic. \end{prop}

\begin{proof}    
In this case, we can apply the same 
template matching as in Lemma~\ref{rotative-outermost}.  Let $N=T^2\times[-1,1]$ 
be an outermost rotative layer with $T_1=\bdry M$, and $A_{[-1,1]}$ the 
corresponding horizontal annulus.  Let $\mathcal{A}$ be the set of 
configurations on $A_{[1,2]}$, corresponding to nonrotative $T^2\times[1,2]$ 
for which $M\cup (T^2\times [1,2])$ remains tight.  We claim that $\mathcal{A}$ 
once again is the set of $\Gamma_{A_{[1,2]}}$ for which no homotopically trivial 
dividing curves appear after merging with $A_{[-1,1]}$.  Note that there might be 
some attachments of $T^2\times [1,2]$ for which the twisting increases by a 
multiple of $\pi$ when we compare $T^2\times [-1,1]$ and $T^2\times [-1,2]$.  This 
happens when homotopically nontrivial closed curves are created on 
$A_{[-1,2]}$.  The tightness is guaranteed by Colin's gluing theorem for 
universally tight contact structures along incompressible tori (see 
\cite{Co99a}).  Finally, $\mathcal{A}$ is sufficient to recover $A_{[-1,1]}$.  
This proves that any two rotative outermost layers are contact diffeomorphic.  
\end{proof}

\subsection{} We make some remarks.  Although we were able to corral in the 
nonrotative outer layers up to disk-equivalence using Theorem~\ref{disk-equiv}, 
the {\it exact set} of allowable nonrotative outer layers for a fixed $(M,\xi)$ 
with torus boundary is much more difficult to determine.  

One of the difficulties (though by no means the only one) is our inability to 
answer the following question:

\medskip

{\bf Question}\qua Let $(M,\xi)$ be a tight contact manifold with a fixed 
convex torus boundary component $T$, and let $T^2\times 
I\subset M$ be a nonrotative outer layer with $T_1=T$. If $T^2\times I$ can be 
extended to a rotative toric annulus inside $M$, then can any other nonrotative 
outer layer $(T^2\times I)'$ in $M$ with $T_1'=T$ be extended to a rotative 
toric annulus inside $M$? 

\medskip

If such a statement is true, it can be proved only by probing deeper into the 
manifold.  In other words, nonrotative outer layers do not always exhibit purely 
superficial data.

\medskip

{\bf Acknowledgements}\qua   The author gratefully acknowledges the American 
Institute of Mathematics and IHES for their hospitality.   He also thanks Tanya 
Cofer, whose  computer implementation of the Gluing Theorem indicated that 
Proposition 5.8 in \cite{H1} was incorrect.

The author is supported by NSF grant DMS-0072853 and the American Institute of 
Mathematics.


\begin{thebibliography}


\bibitem{B}{\bf D Bennequin}, {\it Entrelacements et \'equations de Pfaff}, Ast\'erisque, {107--108}
(1983) 87--161

\bibitem{E89}{\bf Y Eliashberg}, {\it Classification of overtwisted contact structures on 3--manifolds},
Invent. Math. {98} (1989) 623--637

\bibitem{E90}{\bf Y Eliashberg}, {\it Topological characterization of
Stein manifolds of dimension $>$ 2}, Intern. Journal of Math. {1}
(1990) 29--46

\bibitem{E91}{\bf Y Eliashberg}, {\it Filling by holomorphic discs and
its applications}, London Math. Soc. Lecture Note Series, {151} (1991)
45--67


\bibitem{E}{\bf Y Eliashberg}, {\it Contact 3--manifolds, twenty years since J
Martinet's work}, Ann. Inst. Fourier {42} (1992) 165--192

\bibitem{Et}{\bf J Etnyre}, {\it Tight contact structures on lens
spaces}, to appear in Communications in Contemp. Math.

\bibitem{Et2}{\bf J Etnyre}, {\it Transversal torus knots},
Geom. Topol. {3} (1999) 253--268

\bibitem{EH99} {\bf J Etnyre}, {\bf K Honda}, {\it On the
non-existence of tight contact structures}, preprint (1999) available
at {\tt http://www.math.uga.edu/\char126 honda} and {\tt 
arXiv:math.GT/9910115}


\bibitem{ET}{\bf Y Eliashberg}, {\bf W Thurston}, {\it Confoliations},
University Lecture Series, {13}, Amer. Math. Soc. Providence (1998)

\bibitem{F}{\bf M Fraser}, {\it Classifying Legendrian knots in tight
contact 3--manifolds}, PhD  thesis (1994)

\bibitem{Giroux}{\bf E Giroux}, {\it Une structure de contact, m\^eme tendue,
est plus ou moins tordue}, Ann. Scient. Ec. Norm. Sup. {27} (1994)
697--705

\bibitem{Giroux2}{\bf E Giroux}, {\it Convexit\'e en topologie de contact},
Comm. Math. Helv.  {66} (1991) 637--677

\bibitem{Giroux3}{\bf E Giroux}, {\it Structures de contact en
dimension trois and bifurcations des feuilletages de surfaces},
Invent. Math. 141 (2000) 615--689


\bibitem{Gompf}{\bf R Gompf}, {\it Handlebody construction of Stein
surfaces}, Ann. of Math. 148 (1998) 619--693

\bibitem{Gromov}{\bf M Gromov}, {\it Pseudo-holomorphic curves in
symplectic manifolds}, Invent. Math. {82} (1985) 307--347


\bibitem{Hatcher}{\bf A Hatcher}, {\it Notes on basic 3--manifold topology},
preliminary notes

\bibitem{H2}{\bf K Honda}, {\it On the classification of tight contact
strutures II}, to appear in Jour. Diff. Geom.,   available at
{\tt http://www.math.uga.edu/\char126 honda}

\bibitem{H3}{\bf K Honda}, {\it Gluing tight contact structures},
preprint (2000) available at {\tt http://www.math.uga.edu/\char126
honda}

\bibitem{K}{\bf Y Kanda}, {\it The classification of tight contact
structures on the 3--torus}, Comm. in Anal. and Geom. {5} (1997)
413--438

\bibitem{K98} {\bf Y Kanda}, {\it On the Thurston--Bennequin invariant
of Legendrian knots and non exactness of Bennequin's inequality},
Invent. Math. {133} (1998) 227--242


\bibitem{LM}{\bf P Lisca}, {\bf G Mati\'c}, {\it Tight contact
structures and Seiberg--Witten invariants,} Invent. Math. {129} (1997)
509--525


\bibitem{ML}{\bf S Makar-Limanov}, {\it Tight contact structures on
solid tori}, Trans. Amer. Math. Soc. {350} (1998) 1045--1078



\end{thebibliography}

\begin{thebibliography}

\bibitem{Co99a}
{\bf V Colin}, {\it Recollement de vari\'et\'es de contact tendues}, Bull.\ Soc.\ 
Math.\ France {127} (1999) 43--69

\bibitem{Gi91}
{\bf E Giroux}, \textit{Convexit\'e en topologie de contact}, Comment.\ Math.\ 
Helv.\ {66} (1991) 637--677

\bibitem{Gi99a}
{\bf E Giroux}, {\it Une infinit\'e de structures de contact tendues sur une 
infinit\'e de vari\'et\'es}, Invent.\ Math.\ {135} (1999) 789--802

\bibitem{H1} {\bf K Honda}, \textit{On the classification of tight
contact structures I}, Geom.\ Topol.\ {4} (2000) 309--368

\bibitem{H22} {\bf K Honda}, {\it On the classification of tight
contact structures II}, J.\ Differential Geom.\ {55} (2000) 83--143

\bibitem{H32} {\bf K Honda}, {\it Gluing tight contact structures}, to
appear in Duke Math.\ J.\

\end{thebibliography}
\end{document}